\documentclass[11pt,a4paper]{article}

\usepackage{titling}
\title{A Li--Yau and Aronson--Bénilan approach for the Keller--Segel system with critical exponent}

\newcommand{\authorPDF}{Elbar, Fernandez-Jimenez, Santambrogio.}
\newcommand{\subjectPDF}{%
    35Q92,  	
    35K57, 
    35B33, 
    35A15, 
    49J20. 
}
\newcommand{\keywordsPDF}{%
    Keller--Segel equation, Li--Yau, Aronson--Bénilan, Liouville equation, Lane--Emden equation, optimal-control theory.
}

\usepackage{authblk}

\author[1]{Charles Elbar\thanks{elbar@math.univ-lyon1.fr}}
\author[1,2]{Alejandro Fernández-Jiménez\thanks{a.fernandezjimenez@vu.nl}}
\author[1]{Filippo Santambrogio\thanks{santambrogio@math.univ-lyon1.fr}}
\affil[1]{Universite Claude Bernard Lyon 1, ICJ UMR5208, CNRS, Ecole Centrale de Lyon, INSA Lyon, Université Jean Monnet, 69622 Villeurbanne, France.}
\affil[2]{Vrije Universiteit Amsterdam,
    Department of Mathematics,
    Faculteit der Exacte Wetenschappen,
    De Boelelaan 1081a,
    1081 HV Amsterdam, The Netherlands.}

\usepackage{amsmath}
\usepackage{amsthm}
\usepackage{amsfonts}
\usepackage{amssymb}

\input{./tex/commands/document}

\usepackage[utf8]{inputenc}
\usepackage[T1]{fontenc}
\usepackage{dsfont}
\usepackage[UKenglish]{babel}

\newcommand{\R}{{\mathbb{R}}}

\newcommand{\Rd}{{\R^d}}

\renewcommand{\d}{\mathrm{d}}

\newcommand{\dive}{\mbox{div}}

\newcommand{\dd}{\mathop{}\!\d}

\newcommand{\dx}{\dd x}

\usepackage{scalerel}[2014/03/10]
\usepackage[usestackEOL]{stackengine}
\newcommand{\dashint}{\,\ThisStyle{\ensurestackMath{%
  \stackinset{c}{.2\LMpt}{c}{.5\LMpt}{\SavedStyle-}{\SavedStyle\phantom{\int}}}%
  \setbox0=\hbox{$\SavedStyle\int\,$}\kern-\wd0}\int} 

\usepackage{mathtools}

\DeclarePairedDelimiter{\positive}{(}{)^{+}}
\DeclarePairedDelimiter{\negative}{(}{)^{-}}

\newcommand\pos\positive
\renewcommand\neg\negative

\WithSuffix\newcommand\pos*{\positive*}
\WithSuffix\newcommand\neg*{\negative*}

\usepackage[dvipsnames]{xcolor}

\definecolor{color1}{RGB}{0, 121, 178}
\definecolor{color2}{RGB}{255, 124, 37}
\definecolor{color3}{RGB}{37, 160, 55}
\definecolor{color4}{RGB}{220, 32, 44}
\definecolor{color5}{RGB}{147, 104, 186}
\definecolor{color6}{RGB}{143, 85, 76}
\definecolor{color7}{RGB}{230, 119, 192}
\definecolor{color8}{RGB}{127, 127, 127}
\definecolor{color9}{RGB}{192, 188, 55}
\definecolor{color10}{RGB}{0, 191, 206}

\newcommand{\plotHeight}{2.25in}
\newcommand{\plotAspectRatio}{1.35}

\usepackage{tikz}

\usepackage{pgfplots}
\pgfplotsset{compat=1.16}

\usepgfplotslibrary{groupplots}

\pgfplotsset{every axis/.append style={
			axis background/.style={fill=gray!10},
			axis line style={draw=none},
			grid=both,
			grid style={white, line width=.1pt},
			legend style={
					line width = 1pt,
					draw=none,
					/tikz/every even column/.append style={column sep=0.5cm}
				},
			line width=1pt,
			major grid style={white, line width=1.5pt},
			tick style={draw=none},
			xlabel = $x$,
			ylabel = $y$,
			width=1in,
			height=\plotHeight,
			width={\plotAspectRatio*\plotHeight},
		}}

\theoremstyle{plain}
\newtheorem{theorem}{Theorem}[section]
\newtheorem{lemma}[theorem]{Lemma}
\newtheorem{proposition}[theorem]{Proposition}
\newtheorem{corollary}[theorem]{Corollary}

\newtheorem{problem}[theorem]{Problem}

\theoremstyle{remark}
\newtheorem{remark}[theorem]{\bf Remark}
\newtheorem{definition}[theorem]{\bf Definition}

\newcounter{review}

\makeatletter

\newcommand\listreviewname{List of Reviews}
\newcommand\listofreviews{%
	\section*{\listreviewname}\@starttoc{tor}}
\makeatother

\newif\ifreview

\newcommand {\f}{\frac}
\newcommand{\diff}{\mathop{}\!\mathrm{d}}

\usepackage{enumitem}
\usepackage{esint}
\usepackage{float}
\usepackage{tikz}
\usepackage{bbm}
\usepackage[most]{tcolorbox}
\usepackage{subcaption}

\usepackage{cite}

\newcommand{\Lip}{\mathrm{Lip}}

\newlist{enumeratetheorem}{enumerate}{2}
\setlist[enumeratetheorem,1]{
	label=\rm \roman*), 
	ref=Theorem \thetheorem--Item \roman*,
	topsep=0ex
	}
\setlist[enumeratetheorem,2]{label=\rm \alph*), ref=\theenumeratetheoremi.\alph*}

\newlist{enumeratedefinition}{enumerate}{2}
\setlist[enumeratedefinition,1]{
	label=\rm \roman*), 
	ref=Definition \thedefinition--Item \roman*,
	topsep=0ex}
\setlist[enumeratedefinition,2]{label=\rm \alph*), ref=\theenumeratedefinitioni.\alph*}

\newlist{enumeratesteps}{enumerate}{2}
\setlist[enumeratesteps,1]{
  label=\textit{Step \arabic*\protect\thisstep},
  ref=\arabic*,
  align=left, leftmargin=0pt, labelindent=!,
  listparindent=\parindent, labelwidth=0pt, itemindent=!
}
\setlist[enumeratesteps,2]{
    label= \textit{Step \theenumeratestepsi.\alph*\protect\thisstep},
    ref=\theenumeratestepsi.\alph*,
    align=left, leftmargin=0pt, labelindent=\parindent,
  listparindent=\parindent, labelwidth=0pt, itemindent=!
}
\newcommand{\step}[1][]{%
  \if\relax\detokenize{#1}\relax
    \def\thisstep{.}%
  \else
    \def\thisstep{: \protect #1.}%
  \fi
  \item}

\Crefname{enumeratestepsi}{Step}{Steps}
\Crefname{enumeratestepsii}{Step}{Steps}

\newlist{enumeratecases}{enumerate}{2}
\setlist[enumeratecases,1]{
  label=\textit{Case \alph*\protect\thiscase},
  ref=\arabic*,
  align=left, leftmargin=0pt, labelindent=!,
  listparindent=\parindent, labelwidth=0pt, itemindent=!
}
\setlist[enumeratecases,2]{
    label= \textit{Case \theenumeratecasesi.\alph*\protect\thiscase},
    ref=\theenumeratecasesi.\roman*,
    align=left, leftmargin=0pt, labelindent=\parindent,
  listparindent=\parindent, labelwidth=0pt, itemindent=!
}
\newcommand{\case}[1][]{%
  \if\relax\detokenize{#1}\relax
    \def\thiscase{.}%
  \else
    \def\thiscase{: \protect #1.}%
  \fi
  \item}

\Crefname{enumeratecasesi}{Case}{Cases}
\Crefname{enumeratecasesii}{Case}{Cases}

\newcommand{\impl}{\Rightarrow}
\newcommand{\ind}{\mathbbm{1}}
\newcommand{\ve}{\varepsilon}

\date{}

\begin{document}

\begin{singlespace}\maketitle\end{singlespace}

\begin{abstract}
The Li--Yau and Aronson--Bénilan estimates are fundamental differential inequalities for heat and porous-medium type equations: they control the pressure pointwise, imply differential Harnack inequalities, and yield sharp smoothing and decay rates. We develop analogues for the critical parabolic--elliptic Keller--Segel system on $\mathbb{R}^d$, where the nonlinear diffusion competes with the aggregation on the same scale.

Our main tool is the effective pressure $v=p-u$, where $p=\log\rho$ in dimension $2$ and $p=\frac{m}{m-1}\rho^{m-1}$ for $d>2$, with $m=2-\frac{2}{d}$. In terms of $v$ the equation becomes
$
\partial_t\rho=\operatorname{div}(\rho\nabla v),
$ 
which yields Li--Yau estimates in dimension $2$ and Aronson--Bénilan estimates for $d>2$. In each case we obtain quantitative lower bounds on $\Delta v$,
showing that the classical pressure estimates survive, in a suitable form, in the presence of the nonlocal Keller--Segel interaction.

On one hand these estimates recover, in a unified way, results that have so far been obtained with different techniques; on the other, they yield new ones and they have the potential to open the door to a broader program for regularity theory in the Keller--Segel system and for aggregation--diffusion equations where the pressure becomes a systematic tool. For nonnegative finite-measure initial data of sufficiently small mass (assuming only a logarithmic moment in dimension $2$), we construct global solutions that are smooth for positive times and satisfy $\|\rho_t\|_{L^\infty}\leq C/t$ in general dimension. For initial data with finite free energy, we prove global $L^\infty$ bounds of order $1+t^{-1}$ in the subcritical and critical mass regime. We introduce in this setting a new notion of critical mass, which is the minimal mass $M_c^{\mathrm{sub}}(d)$ of subsolutions of the elliptic Lane--Emden equation (Liouvile in the two-dimensional case). We then prove that this minimal mass is the expected value $8\pi$ for $d=2$ and, more generally, coincides with the the Keller--Segel critical mass coming from the HLS inequality for $d\geq 3$. This requires reducing the computation of the minimal mass to a 1D optimal control problem.

Along the way we develop several tools of independent interest: $L^\infty$ bounds derived from lower estimates on $\Delta v$, critical log-HLS/HLS type estimates for families far from a Dirac mass, Hessian bounds for the Newtonian potential, and detailed Lipschitz and Hessian estimates for bounded solutions of porous-medium equations with drifts on the whole space, which are required to justify the comparison principle.
\end{abstract}

\subjectclassification{\subjectPDF}

\newpage

\tableofcontents

\section{Introduction}


This paper focuses on the Keller--Segel system in $\mathbb{R}^d$:
\begin{equation}\label{eq:KS}\tag{KS}
\begin{cases}
\displaystyle\partial_t\rho   = \Delta \rho^m  -  \mathrm{div}\left(\rho \,\nabla u \right) 
& \text{in } (0,\infty) \times \mathbb{R}^d,\\[6pt]
-\Delta u  =  \rho 
& \text{in } (0,\infty)\times \mathbb{R}^d,
\end{cases}
\end{equation}
which arises as a biological model where $\rho$ denotes the density of cells or bacteria and $u$ is a chemo-attractant that they produce and move toward. Furthermore, we take the initial condition $\rho_{0}\ge 0$, that is $\rho(0,\cdot)=\rho_{0}$. The dimension is $d\ge 2$, and we  consider the critical exponent $m = 2 - \frac{2}{d}$. This exponent is referred as \textit{critical} since provides a fair balance between the nonlinear diffusion $\Delta\rho^m$ and the chemotactic attraction given by $\mathrm{div}(\rho\,\nabla u)$. Because of this balance, and only for this very choice of the exponent $m$ (see \cite{Carrillo_Craig_Yao19}), the results strongly depend on the mass of the solution $M=\int_{\R^d}\rho(t,x)\diff x$, which is preserved in time. In the model, $u$ is chosen to be the solution of the Poisson equation:
\[
u(t,x)  =  (\Gamma\ast\rho_t)(x),
\]
where $\Gamma$ refers to the Newtonian potential in $\mathbb{R}^d$, i.e.
\[
\Gamma(x)  = 
\begin{cases}
  \displaystyle \f{1}{(d-2)\,\sigma_d}\,|x|^{2-d} & \text{if } d > 2, \\
  \displaystyle -\f{1}{2\pi}\,\log|x| & \text{if } d = 2,
\end{cases}
\]
and $\sigma_d$ is the $d-1$-dimensional measure of the unit sphere $\mathbb{S}^{d-1}\subset \mathbb{R}^d$ (so that $\sigma_d=d\omega_d$ where  $\omega_d$ is the volume of the unit ball in $\R^d$).

The model was proposed for the first time by Keller and Segel in \cite{Keller_Segel71} and Patlak in \cite{Patlak53} to explain a phenomenon known as \textit{chemotaxis}.
The Keller--Segel system has further applications for different biological phenomenon connected to the evolution of unicellular organism to come up with more complex structures. For instance, it explains pattern formation of cells through meiosis \cite{Calvez_Hawkins_Meunier_Voituriez12}; embryogenesis and angiogenesis \cite{Chaplain96, Levine_Sleeman_Nilsen-Hamilton01}; the Balo disease \cite{Khonsari_Calvez07}; and bio-convection \cite{Chertock_Fellner_Kurganov_Lorz_Markowich12}. In physics, the Keller--Segel system also models the movement of self-gravitating Brownian particles \cite{Chavanis03}.

With respect to the mathematical analysis, over the last few decades there has been a lot of work devoted to this problem, see for instance the surveys \cite{Horstmann03, Perthame04, Carrillo_Craig_Yao19}. As already mentioned above, this system presents a dichotomy between diffusion and aggregation. For the $2$-dimensional case this dichotomy depends on the evolution of the second moment. In particular, in \cite{Blanchet_Dolbeault_Perthame06}, Blanchet, Dolbeault and Perthame prove a remarkable identity on the evolution of the second moment $
\mathrm{m}_{2}[\rho]\coloneqq \int_{\R^d}|x|^2\rho(x)\diff x
$. They show that the evolution in time of the second moment is given by
\begin{equation*}
    \frac{\diff}{\diff t} \mathrm{m}_{2}[\rho_t] = 4M \left( 1 - \frac{M}{8 \pi} \right).
\end{equation*}
In particular, this identity shows that $M= 8 \pi$ is the critical mass and it gives three different regimes depending on the mass. If $M < 8 \pi$ (\textit{subcritical mass}) and the initial datum has finite free energy \eqref{eq:free_energy}, then there exists global in time solutions \cite{Jager_Luckhaus92, Nagai_Senba_Yoshida97, Biler_Karch_Laurencot_Nadzieja06, Blanchet_Dolbeault_Perthame06, Blanchet_Calvez_Carrillo08, Campos_Dolbeault14,2024arXiv241015095E}. Furthermore, in \cite{Bedrossian_Masmoudi14}, Bedrossian and Masmoudi show local well-posedness for measure-value initial data which do not include atoms of mass $8\pi$. If $M > 8 \pi$ (\textit{supercritical mass}), solutions blow-up in finite time. 
In \cite{Herrero_Velazquez96}, Herrero and Velázquez construct an example of a radial solution with finite time blow-up with a subsequent detailed stability analysis by Velázquez \cite{Velazquez02, Velazquez04a, Velazquez04b}. Afterwards, in \cite{Raphael_Schweyer14} Rapha\"el and Schweyer show the existence of singularities for the radial case with mass close to $8 \pi$. Finally, in a later work \cite{Collot_Ghoul_Masmoudi_Nguyen22, Collot_Ghoul_Masmoudi_Nguyen22b}, Collot, Ghould, Masmoudi and Nguyen refined the analysis in order to cover further cases.   
If the total mass is $M= 8 \pi$ (\textit{critical mass}) and the initial free energy \eqref{eq:free_energy} is finite, solutions exists globally in time but they blow-up at time infinity \cite{Biler_Karch_Laurencot_Nadzieja06, Blanchet_Carrillo_Masmoudi08}. A proof of the global existence in the critical mass case $8\pi$ based on the Wasserstein gradient flow structure and on the JKO sheme is contained in \cite{KarWol}.

The community has also studied its basin of attraction \cite{Blanchet_Carlen_Carrillo12, Lopez-Gomez_Nagai_Yamada13, Lopez-Gomez_Nagai_Yamada14}, which includes, among others, oscillating behaviour when the second moment is infinite. Moreover, when the second moment is finite there has been work focused on the study of the blow-up profiles \cite{Ghoul_Masmoudi18, Davila_delPino_Dolbeault_Musso_Wei24}.

In the case of higher dimension there exists also a dichotomy depending on the initial mass, \cite{Sugiyama07}. If $M < M_c$ the solutions remain bounded globally in time and there exists self-similar solutions decaying in time with the porous medium equation \cite{Blanchet_Carrillo_Laurencot09, Bedrossian11}. If $M > M_c$, in \cite{Bedrossian_Kim13} Bedrossian and Kim show that all radial solutions blow-up in finite time (but the finite-time blow-up of {\it all} solutions with supercritical mass is an open problem). If $M = M_c$ the solutions are globally well-posed and their $L^\infty$ norm is globally bounded in time \cite{Blanchet_Carrillo_Laurencot09}. Furthermore, Yao showed that radial solutions with compact support converge to some stationary solution in this family \cite{Yao14}. Further properties of the Keller--Segel equation with critical exponent in higher dimension have been studied, e.g. \cite{Kim_Yao12, Calvez_Carrillo_Hoffmann17, Carrillo_Santambrogio18, Carrillo_Delgadino_Patacchini19, Carrillo_Hittmeir_Volzone_Yao19}, but the literature on this case is not as extensive as for its two-dimensional counterpart. For further remarks on the choice of the exponent $m$ and relevant references we refer the reader to \Cref{rem:Different m}. 

Finally, we also mention results on the uniqueness of solutions, which is proven for $d=2$ in \cite{Carrillo_Lisini_Mainini14} and for $d > 2$ in \cite{Bedrossian_Rodriguez_Bertozzi11}, in both cases in the class of bounded solutions.

\medskip

Much effort has been spent on the properties of \eqref{eq:KS} over the last three decades. Nevertheless, in this manuscript, we propose an approach in order to understand the Keller--Segel equation from a new viewpoint. 
Thus, the main result of this work is to prove that the classical Li--Yau or Aronson--B\'enilan estimates,  originally proved for purely diffusive equations, hold when one adds the chemotactic drift term from the Keller--Segel system. Let us start by introducing these inequalities on the purely diffusive case.

\begin{itemize}
\item For the heat equation $\partial_t\rho = \Delta\rho$, Li and Yau~\cite{Li_Yau86, Yau94} showed
\[
\Delta \left(\log\rho_t \right)   \ge -\frac{d}{2t}.
\]
Note that in this case one can obtain a stronger result, i.e. $D^2\left(\log\rho_t \right)   \ge -\frac{1}{2t}I$, where $I$ is the identity matrix. This is known as the matrix Li–-Yau–-Hamilton estimate, see \cite{Li_Yau86, Hamilton93}.
\item For the porous medium equation $\partial_t\rho = \Delta\rho^m$, Aronson and B\'enilan~\cite{Aronson70, Aronson_Benilan79} derived
\[
\Delta  \left(\frac{m}{m-1}\,\rho^{m-1}_t \right)
     \ge 
   -\frac{\alpha}{\,t\,},
   \quad
   \alpha  = \frac{d}{d(m-1)+2}
\]
\end{itemize}
(note that in the case $m=2-\frac 2d$ we have $\alpha=1$). This lower bound on the Laplacian of the pressure is used in \cite{Aronson_Benilan79} in order to prove that there exists a
unique strong and continuous solution to the Cauchy problem with $L^1$-bounded initial data. See also \cite[Chapter 9]{Vazquez06} for a review on the Aronson--B\'enilan estimate for the porous medium equation. Afterwards, this estimate has been extended in order to cover further cases such as the filtration equation \cite{Crandall_Pierre82}, or the fast-diffusion equation \cite{
Vazquez06, Bonforte_Figalli24}. This family of inequalities have proven to be very relevant in the purely diffusive case since they provide strong properties for the solutions, like global well posedness, long-time asymptotic behaviors, regularity, etc. In the case of the porous medium equation, they also allow to better understand the hyperbolic behavior near the interfaces, \cite[Chapter 9]{Vazquez06}. Recovering this family of bounds for aggregation--diffusion equations has the potential to push the field towards a better understanding and the derivation of strong properties of the Keller--Segel system which are still not available in the literature. 

Therefore, in this work, we show that when adding a nonlocal drift modeling chemotaxis, similar inequalities still hold. More precisely, we will consider lower bounds on $\Delta v$ where $v=p-u$ for $p=\log\rho$ when $d=2$ or $p=\frac{m}{m-1}\,\rho^{m-1}$ when $d>2$. Let us point out that in both cases the Keller--Segel system can be written as $\partial_t\rho=\dive(\rho\nabla v )$. This implies robustness of the Li--Yau and Aronson--B\'enilan estimates even with a nonlocal term. In addition to this, we show an estimate of the form $\|\rho\|_{L^\infty}\leq C|\inf\Delta v|$. Thereby, from the Li--Yau and Aronson--B\'enilan type lower bounds on $\Delta v$ we deduce $L^\infty$ bounds on $\rho$. Finally, combining everything, we also obtain global existence results.

Furthermore, our result is consistent with the dichotomy of the mass and it allows us to recover and improve some of the classical results concerning global existence in two different ways: first, we prove an explicit $L^{\infty}$ estimate on the solutions, that it is also decreasing in time. Second, we are now able to consider initial data to be only measures, if the mass is small enough. If the mass is close to the critical mass, we only require the free energy of the initial condition to be bounded.

\begin{definition}[Critical mass]\label{def:critical_mass}
We define the critical mass as
\begin{equation}\label{def:Mc}
M_{c}(d)  = 
\begin{cases}
\displaystyle\,8\pi, & d=2,\\[3pt]
\displaystyle\,\left(\frac{2}{C_{*}(m-1)}\right)^{\frac{1}{2-m}}, & d>2.
\end{cases}
\end{equation}
Here,  $C_{*}$ is the best (smallest) constant in the inequality
\begin{equation}\label{eq:HLS}
\int_{\R^d}\rho(x) u(x)\diff x   \le   C_{*}\,\|\rho\|_{L^m}^{m} \,\|\rho\|_{L^1}^{2-m},
\end{equation}
where $m = 2 - \frac2d$ and $-\Delta u =\rho$. This is valid for all nonnegative $\rho\in L^{1}(\R^d)\cap L^{m}(\R^d)$. The existence of the constant $C_*$ is obtained from Hardy--Littlewood--Sobolev (HLS) inequalities combined with H\"older’s inequality, see for instance~\cite{Blanchet_Carrillo_Laurencot09, Calvez_Carrillo_Hoffmann17}. For convenience, we will loosely refer to \eqref{eq:HLS} as “the HLS inequality”, or the Lane--Emden inequality although it is slightly adapted to our Keller--Segel setting.

It is also well-known (see \cite[Proposition 3.5]{Blanchet_Carrillo_Laurencot09})  that the critical mass $M_c(d)$ for $d>2$ introduced above is equal to the mass of the density $\bar\rho \coloneqq \bar h^q$ where $q=d/(d-2)$ and $\bar h$ is a radially decreasing solution of 
$$
\begin{cases}
    \f{m}{m-1}\Delta \bar{h} + \bar{h}^q = 0 &\text{in } B_R(0),\\
    \bar h = 0 &\text{on } \partial B_R(0),
\end{cases}
$$
for some $R>0$. Given the value $\bar h(0)$ there exists unique a solution of this equation for a certain $R>0$, and the value of $R$ and/or of $\bar h(0)$ can be modified by scaling since $x\mapsto \lambda^{d-2}\bar h(\lambda x)$ is also a solution for every $\lambda>0$. 
\end{definition}

Equation \eqref{eq:KS} can be viewed as the gradient flow in the $2$-Wasserstein metric of the free-energy functional \cite{Otto01,Carrillo_McCann_Villani03,Ambrosio_Gigli_Savare08,Santambrogio15} defined by 
\begin{equation}\label{eq:free_energy}
\mathcal{F}[\rho]  =  
\begin{cases}
 \displaystyle \int_{\mathbb{R}^2}\rho(x)\log\rho(x) \,\diff x  - \f12 \int_{\mathbb{R}^2}\rho(x) u(x)\diff x, \quad & d=2,\\
 
  \displaystyle \int_{\mathbb{R}^d}\f{\rho(x)^m}{\,m-1\,}\,\mathrm{d}x  - \f12\,\int_{\mathbb{R}^d}\rho(x)\,u(x)\diff x, \quad & d\ge 3. 
\end{cases}
\end{equation}
The critical exponent $m=2-\f{2}{d}$ can be obtained by requiring that the two terms in the free energy compete on the same scaling. One takes $\rho_\lambda = \lambda ^d\rho(\lambda x)$ and each term should scale the same with respect to $\lambda$.

The critical mass dictates the behaviour of the free energy along the solutions of the Keller--Segel system. In dimension 2, the analog of the HLS inequality is the so-called {\it log HLS} inequality. This inequality is relevant to us since it provides a lower bound on $\mathcal{F}$ on the subcritical and critical mass regimes: $ \mathcal{F}[\rho]$ is bounded from below by a finite constant whenever $\rho$ has mass $M=8\pi$, while for $M<8\pi$ then we can bound the entropy $\int_{\R^d}\rho\log\rho$ in terms of $\mathcal{F}[\rho]$ and of a suitable moment of $\rho$. Let us also notice that the free energy is decreasing in time. Hence, in the subcritical mass regime, if the initial data has finite free energy,  we deduce a superlinear bound on $\rho$. This estimate can be used in order to obtain global bounded solutions of the Keller--Segel system. In the critical mass regime $M=8\pi$, we cannot control the entropy of $\rho$ via the log HLS inequality. Therefore, the argument is more delicate and requires more involved computations. A crucial tool which has been used to study the critical mass case is a new Lyapunov functional, introduced in~\cite{Blanchet_Carlen_Carrillo12}, given by.
 \[
 H_{\lambda}[\rho]  
   \coloneqq  
 \int_{\R^2}\! \left(\sqrt{\rho(x)} - \sqrt{u_\lambda(x)} \right)^{2}\,u_{\lambda}^{-1/2}(x)\,\diff x,
 \quad 
 u_{\lambda}(x)  
  =  
 \frac{8\,\lambda}{(\lambda + |x|^2)^2}.
 \]
%
%
In higher dimensions, when the mass is subcritical, we can bound the $L^{m}$ norm in terms of the free energy thanks to the HLS inequality. Once again at the critical mass case, we only have lower bounds on $\mathcal{F}[\rho]$ and not coercivity in $L^m$, and proving global existence of solutions becomes challenging. This was proved in~\cite{Blanchet_Carrillo_Laurencot09}.   

In our analysis, we find another notion of critical mass that is important and it is related to the minimal mass of subsolutions of a certain nonlinear elliptic equation. In order to introduce it, let us recall that we are interested in lower bounds on $\Delta v$ where $v=p-u$ (so that $\Delta v=\Delta p+\rho$) and $p=\log\rho$ when $d=2$ or $p=\frac{m}{m-1}\,\rho^{m-1}$ when $d>2$. Hence, for the analysis introduced in this manuscript, it is key to understand for which values of the mass $M$ it is possible to find subsolutions $\rho\geq 0$ satisfying $\|\rho\|_{L^1}=M$ and solving the problem
\begin{equation*}
    \Delta\log\rho + \rho\geq 0 \, \, (\text{if } d=2) \quad\text{or}\quad\Delta \left(\frac{m}{m-1}\,\rho^{m-1}\right) + \rho\geq 0 \, \, (\text{if } d>2).
\end{equation*}
In the case of dimension $d=2$, there exists an explicit radial solution of $\Delta\log\rho + \rho= 0$ with mass equal to $8\pi$. Furthermore, it is known (see \cite[Lemma 1.1]{Chen_Li91}) that any other solution has mass at least $8\pi$. Let us define  
\begin{equation*}
\begin{split}
M^{sub}_c(2)&\coloneqq\inf\left\{\int \rho(x)\diff x\,:\,\rho\geq 0, \, \Delta \left(\log\rho \right) + \rho\geq 0\text{ in $\{x: \rho(x)>t\}$ for all $t>0$}\right\}, \\
&=\inf\left\{\int e^h\diff x\,: \Delta h + e^h\geq 0\text{ in $\Omega_t=\{x: h(x)>t\}$ for all $t>-\infty$}\right\}.\notag
\end{split}
\end{equation*}
We can adapt the argument introduced in  \cite[Lemma 1.1]{Chen_Li91} in order to show that  $M^{sub}_c(2)=M_c(2)=8 \pi$ as well. However, in the case of dimension $d > 2$, the argument is less straight forward. Let us For $m=2-\frac2d$ and $q=\frac{d}{d-2}$, we define
\begin{equation}\label{def:Mc LaneEmden}
\begin{aligned}
M^{sub}_c(d)
&\coloneqq \inf\left\{
\int \rho(x)\diff x :
\rho\ge0,\ \rho\not\equiv0,\right. \left.
\Delta\left(\frac{m}{m-1}\rho^{m-1}\right)+\rho\ge0
\right\}\\
&=\inf\left\{
\int h^q\diff x :
h\ge0,\ h\not\equiv0,\right.\left.
\Delta h+\frac{1}{q+1}h^q\ge0
\right\}.
\end{aligned}
\end{equation}
It is known that there is no entire solution to the Lane-Emden equation $\Delta h + ch^q=0$ when $1< q< \frac{d+2}{d-2}$, see \cite{Bahri_Lions92}. However, subsolutions exist, and we prove below that their minimal mass is exactly $M_c(d)$.
Later in the paper, we will require $M < M^{sub}_c$ in order to prove some of the results, such as an $L^\infty$ bound on $\rho$, or $M = M^{sub}_c$ under some extra assumptions that we will explain later. 
Thus, we would like to show $M^{sub}_c(d)=M_c(d)$ so that the two notions of critical mass agree. To do so, we enlarge the minimization problem defining $M^{sub}_c(d)$ to a larger class stable by radial rearrangement, and then we prove the desired equality by studying the optimizers of this new problem.

\subsection{Main results}

We now present our main theorems. They provide pointwise bounds on the Laplacian of the pressure and yield uniform bounds on the density under some conditions on the total mass $M$. In dimension $d=2$, we recover a Li--Yau type differential inequality, whereas for higher dimension $d>2$ we obtain an Aronson--B\'enilan type estimate at the critical exponent $m=2-\f{2}{d}$. In order to introduce our results we need the notation 
$$
m_{V}[\rho]\coloneqq \int_{\R^d}\rho\,V \diff x.
$$
%
We say that a measure $\rho_0$ on $\R^2$ has finite log-moment if and only if there exists a radially increasing function $V$ such that $m_V[\rho_0]<+\infty$ and  $\int_{\R^2}e^{-V}\diff x<+\infty$, which is actually equivalent to $\int \log(1+|x|)d\rho_0<+\infty$. This condition is necessary in order to properly define $\Gamma*\rho_0$ in dimension 2.

\begin{remark}[Convention for the constants depending on the mass]
\label{rem:mass_constants}
When the argument is the mass $M$, we reserve the notation
$C_1(M)$, $C_2(M)$ and $C_3(M)$ for the three successive estimates
$$
\|\rho\|_{L^\infty} \le C_1(M)\alpha, \quad
\|\rho\|_{L^\infty} \le C_2(M)|\delta[\rho]|, \quad
Q(u)  \le C_3(M)|\delta[\rho]|.
$$
In the first estimate, $\alpha$ denotes a constant such that $\Delta\log\rho\ge-\alpha$ if $d=2$ or
$\Delta\rho^{1-2/d}\ge-\alpha$ if $d>2$. For the second estimate,
$\delta[\rho]$ is defined in \eqref{eq:def_delta}, and $Q(u)$ is
defined in Section~\ref{sec:minimum_laplacian} as the maximal difference of the eigenvalues of the Hessian $D^2u$, where $u=\Gamma*\rho$ is the function appearing in the nonlocal term of the Keller-Segel equation. The precise constants
are introduced in
\Cref{lem:Naive Linfty,prop:linf_very_small_mass,prop:linf_estimate,prop:split_cases_smallmass}.

Whenever
$$
\frac{d-1}{\sqrt{2d}}\,C_3(M)<1,
$$
we also define
\begin{equation}\label{eq:def_CAB}
C_{\mathrm{AB}}(M) \coloneqq \left( 1-\frac{(d-1)^2}{2d}C_3(M)^2 \right)^{-1}.
\end{equation}

Note that the optimal value of $C_1(M)$ is known and explicit, and that a non-optimal value for $C_2(M)$ can be computed explicitly from $C_1(M)$ when $C_1(M)$ (or, equivalently, $M$) is small enough. Yet, we will prove the existence of a finite value for $C_2(M)$ for any $M<M^{sub}_c$, without an explicit estimate on it. In what concerns $C_3(M)$, a non-optimal value can be computed from $C_2(M)$, but we will only give an explicit value for it in the case $d=2$, where the formula is quite easy.
\end{remark}

\begin{theorem}[Li--Yau estimate for Keller--Segel in dimension 2]\label{thm:Li--Yau}
Let $\rho_{0}\in \mathcal{M}(\R^2)$ be a nonnegative measure with finite log-moment, and call $M=\rho_{0}(\R^2)$ the mass of the initial condition. We distinguish three different cases.  
\begin{enumeratetheorem}

    \item\label{2d_small_mass}
Small mass. Assume $M<2\pi$. Then, there exists a global solution
$\rho$ of the Keller--Segel system~\eqref{eq:KS} with initial value
$\rho_0$ such that, for any $t\in(0,T]$, we have in the
distributional sense
$$
\Delta[\log\rho_t-u_t]
    \ge -\frac{C_{\mathrm{AB}}(M)}{t}, \quad C_{\mathrm{AB}}(M) =\frac{(4\pi-M)^2}{8\pi(2\pi-M)}.
$$
Moreover,
$$
\|\rho_t\|_{L^\infty}
    \le
    \frac{C_2(M)C_{\mathrm{AB}}(M)}{t}
    =
    \frac{M(4\pi-M)}{8\pi(2\pi-M)}\frac{1}{t},
$$
and
$$
\Delta[\log\rho_t]
    \ge
    -\frac{(1+C_2(M))C_{\mathrm{AB}}(M)}{t}
    =
    -\frac{4\pi-M}{2(2\pi-M)}\frac{1}{t}.
$$
    \item\label{2d_full_mass} Subcritical mass. If $M\in(0,8\pi)$ and the initial condition is such that $\mathcal{F}[\rho_0]<+\infty$, then there  exists a global solution of the Keller--Segel system~\eqref{eq:KS} $\rho$ with initial value $\rho_{0}$ and a constant $C(T)\ge 0$ such that for any $t\in (0,T]$ we have
    $$
\Delta [\log \rho_t - u_t]\ge -C(T)\left(1+\f{1}{t}\right),
$$
where the constant $C(T)$ satisfies $\lim_{T\to 0} C(T)=1$.
Moreover we have separate estimates
$$
\|\rho_t\|_{L^{\infty}}\le C(T)\left(1+\f{1}{t}\right),  \quad \Delta [\log \rho_t ]\ge- C(T)\left(1+\f{1}{t}\right),
$$
for some other finite constants $C(T)$.

\item\label{2d_exact_mass} Critical mass. If $M=8\pi$ and the initial condition is such that $\mathcal{F}[\rho_0]<+\infty$, the same result as in {\rm ii)} is true. Moreover, if $H_{\lambda}[\rho_0]  <+\infty$, then the constant $C(T)$ can be taken independent of $T$.

\end{enumeratetheorem}

\end{theorem}

\begin{remark} Note that we only assume that a \textit{log}-moment is finite instead of the usual second moment. This assumption is more natural and well-known in the Keller--Segel
system~\cite{Fernandez_Mischler16,Calvez_Corrias08}. 
\end{remark}

\begin{remark}
    [Sharpness on the mass for $d = 2$]
    The result presented in \Cref{thm:Li--Yau} is sharp in the sense of mass. If the mass is such that $M > 8 \pi$, i.e. supercritical, the solutions to the Keller--Segel system~\eqref{eq:KS} always present blow-up in finite time, see \cite{Blanchet_Dolbeault_Perthame06}. Therefore, a global-in-time Li--Yau estimate will not occur in this regime. 
\end{remark}
\begin{remark}[Non-sharpness of the small mass threshold]
Point~\ref{2d_small_mass} is presented for masses smaller than
$2\pi$, but this threshold is clearly not sharp. It is obtained by
combining several different arguments in order to finally reach an
inequality of the form $Q(u)\le C_3(M)|\delta|$ (see later for the
definitions of $\delta$ and $Q(u)$), together with the condition
$\frac{d-1}{\sqrt{2d}}C_3(M)<1$, which in dimension $2$ amounts to
selecting masses such that $C_3(M)<2$. To do so we first estimate
$Q(u)$ in terms of $\|\Delta\rho\|_{L^1}$, then
$\|\Delta\rho\|_{L^1}$ in terms of $M$, $|\delta|$ and
$\|\rho\|_{L^\infty}$ (this estimate being for sure not sharp, since
a term with the good sign is neglected), and finally
$\|\rho\|_{L^\infty}$ in terms of $M$ and $|\delta|$, in the form
$\|\rho\|_{L^\infty}\le C_2(M)|\delta|$. The explicit constant
$C_2(M)=\frac{M}{4\pi-M}$ that we use here blows up as $M\to4\pi$,
while we prove later in the paper that an estimate with a
nonexplicit constant also holds for $M<8\pi$. Moreover, the
combination of the different inequalities adds another layer of
non-sharpness, as their equality cases are different. All in all, we
can say that there exists a maximal mass $\varepsilon_2$ (where the
subscript $2$ stands for the dimension) such that bounds of the form
of those in~\ref{2d_small_mass} are satisfied, and we prove
$\varepsilon_2\ge2\pi$. Note that bounds such as
$\|\rho_t\|_{L^\infty}\le C/t$ were already proven for small mass
in~\cite{Bedrossian_Rodriguez_Bertozzi11}.
\end{remark}

\begin{theorem}[Aronson--B\'enilan estimate for Keller--Segel in dimension $d>2$]\label{thm: Aronson--B\'enilan}
Let $\rho_0 \in \mathcal{M}(\R^d)$ a nonnegative measure and $M=\rho_{0}(\R^d)$ the mass of the initial condition. We distinguish three different cases.
\begin{enumeratetheorem}

\item\label{2d_small_mass_multid}
Small mass. Assume that $M$ belongs to the explicit small-mass
range of \Cref{prop:linf_very_small_mass} and satisfies
$$
\frac{d-1}{\sqrt{2d}}C_3(M)<1.
$$
This defines an explicitly computable threshold $\varepsilon_d$.
Then, there exists a global solution $\rho$ of the Keller--Segel
system~\eqref{eq:KS} with initial value $\rho_0$ such that, for any
$t\in(0,T]$,
$$
\Delta\left[
    \frac{m}{m-1}\rho_t^{m-1}-u_t
\right]
    \ge -\frac{C_{\mathrm{AB}}(M)}{t}
$$
in the distributional sense. Moreover,
$$
\|\rho_t\|_{L^\infty}
    \le
    \frac{C_2(M)C_{\mathrm{AB}}(M)}{t},
\qquad
\Delta\left[
    \frac{m}{m-1}\rho_t^{m-1}
\right]
    \ge
    -\frac{(1+C_2(M))C_{\mathrm{AB}}(M)}{t}.
$$
All these constants are computable from $M$ and $d$.

    \item\label{2d_full_mass_multid} Subcritical mass. Assume $0<M<M_c(d)$ and that the initial condition has finite free energy, that is $\mathcal{F}[\rho_0]<+\infty$. Then, there exists a global solution of the Keller--Segel system~\eqref{eq:KS} $\rho$ with initial value $\rho_{0}$ and a constant $C(T)\ge 0$ such that for any $t\in (0,T]$ we have in distributional sense
    $$
\Delta \left[\frac{m}{m-1}\rho^{m-1}_t - u_t\right]\ge -C(T)\left(1+\frac{1}{t}\right),
$$
where the constant $C(T)$ is bounded for $T\to\infty$ and satisfies $\lim_{T\to 0}C(T)=1$.

Moreover we have separate estimates
$$
\|\rho_t\|_{L^{\infty}}\leq C\left(1+\f{1}{t}\right) , \quad \Delta \left[\f{m}{m-1}\rho^{m-1} \right]\ge -C\left(1+\f{1}{t}\right)
$$
for some other constants $C$.

\item\label{2d_exact_mass_multid} Critical mass. Assume $M=M_c(d)$ and $\rho_0\in L^m(\R^d)$. Then, there exists a global solution of the Keller--Segel system~\eqref{eq:KS} $\rho$ with initial value $\rho_{0}$ and a constant $C(T)\ge 0$ such that for any $t\in (0,T]$ we have in distributional sense
    $$
\Delta \left[ \f{m}{m-1}\rho^{m-1}_t - u_t \right]\ge -C(T)\left(1+\f{1}{t}\right).
$$
where the constant $C(T)$ satisfies $\lim_{T\to 0}C(T)=1$.

Moreover, we have separate estimates
$$
    \|\rho_t\|_{L^{\infty}}\leq C(T)\left(1+\f{1}{t}\right) , \quad \Delta \left[ \f{m}{m-1}\rho^{m-1} \right] \ge -C(T)\left(1+\f{1}{t}\right)
$$
for some other constants $C(T)$.
\end{enumeratetheorem}

\end{theorem}



Please note the dependence in $T$ of the constant $C(T)$ in \ref{2d_full_mass} and \ref{2d_full_mass_multid}. In the case $d>2$, the dependency in $T$ is used in order to describe its asymptotic behavior as $T\to 0$ (let us remark that $C(T)$ is bounded uniformly in $T$). 

\begin{remark}
[Sharpness on the mass for $d > 2$]
Let us recall that $M_c(d)$ is the critical mass of the Keller--Segel system~\eqref{eq:KS}. In particular, if the mass is such that $M > M_c(d)$ there exist initial data for which the solution blows-up in finite time, see~\cite{Blanchet_Carrillo_Laurencot09}. Therefore, a global-in-time Aronson-B\'enilan estimate with minimal assumptions on the initial datum will not hold whenever $M > M_c(d)$. 

In \Cref{sec:liouville}, we prove that
$M_c^{sub}(d)=M_c(d)$ for every $d\ge2$. Consequently, the mass
threshold in \Cref{thm: Aronson--B\'enilan} agrees with the critical mass of the Keller--Segel system in every dimension.
\end{remark}

\begin{remark}
    [The choice of the exponent $m$]
    \Cref{thm: Aronson--B\'enilan} concerns the choice of exponent $m = m_c = 2 - \frac{2}{d}$, which in the Keller--Segel literature is referred as the critical exponent or fair competition regime, see~\cite{Carrillo_Craig_Yao19}. Such choice of the exponent leads to self-similar solutions and a rich structure of the PDE, which shows a dichotomy between global-in-time solutions and finite time blow-up, depending only on the initial mass, see~\cite{Blanchet_Carrillo_Laurencot09}. Due to this interesting behaviour, the Keller--Segel system with critical exponent $m_c$ is the most popular within the mathematical literature. For general $m$ the available literature is more scarce, the equation looses its self-similarity and different techniques are needed in order to study it. If $m > 2 - \frac{2}{d}$ we are on the so-called diffusion-dominated regime, for which there exists global-in-time existence of solutions, \cite{Bedrossian_Rodriguez_Bertozzi11, Bedrossian_Rodriguez14} in dimension $d = 2$ and  \cite{Calvez_Carrillo06, Sugiyama07b} in dimension $d > 2$. If $0 < m < 2 - \frac{2}{d}$ we fall on the aggregation dominated regime. If $\frac{2d}{2+d} < m < 2 - \frac{2}{d}$ there exists a dichotomy between global-in-time existence and finite-time blow-up depending on the initial value of the $L^p$ norm for $p = \frac{d (2-m)}{2}$, \cite{Chen_Liu_Wang12, Bian_Liu13, Chen_Wang14}. For smaller values of $m$ there exist several partial result showing existence or blow-up for different various particular cases, see for example \cite{Perthame04, Bedrossian11, Carrillo_FJ_Gomez-Castro24, Glogic_Schorkhuber24}. In this manuscript, we will only focus on the case $m = 2 - \frac{2}{d}$ which is the most challenging one, and we leave the case $m \neq 2 - \frac{2}{d}$ as an open problem for future work, for which, in view of the existing literature, new techniques will have to be taken into account.
    \label{rem:Different m}

Furthermore, we remark that the proofs of \ref{2d_full_mass}, \ref{2d_exact_mass}, \ref{2d_full_mass_multid} and \ref{2d_exact_mass_multid} entail a new theoretical framework for Keller--Segel based on a careful analysis of the subsolutions of the Liouville and Lane--Emden equations and an optimal control argument. 
\end{remark}


\subsection{Comparison with the literature}

Our main results, the Li--Yau and Aronson--B\'enilan estimates and their proofs, imply a different series of results, which we highlight here. 

\paragraph{Existence of global solutions.}

Theorems~\ref{thm:Li--Yau} and~\ref{thm: Aronson--B\'enilan} provide very strong bounds and the existence of global-in-time solutions. Similar results of course already exist in the literature both for subcritical and critical masses (see \cite{Blanchet_Dolbeault_Perthame06, Blanchet_Carrillo_Masmoudi08, Blanchet_Dolbeault_Escobedo_Fernandez10, Campos_Dolbeault14}) but we believe that our statement is the most comprehensive one and the assumptions on the initial datum are very natural. Proving global existence under no assumption if the mass is very small, or assuming finite free energy for subcritical masses is very classical (in both cases, for $d=2$ we also need a logarithmic moment on $\rho_0$, which is easy to understand since otherwise one cannot define $u$ as a convolution with the logarithm). In what concerns the case of critical mass, our strategy is based on the use of a compact set $K$ (compact for the narrow convergence of measures) where the evolution takes place and which does not contain Dirac masses (measures concentrated on a single point). We prove (see below) that excluding Dirac masses is enough to obtain a control of the entropy or $L^m$ part of the free energy with the free energy itself (which, instead, also contains the negative interaction term). From this point of view, the classical result for the case $d=2$ with $H_\lambda(\rho_0)<+\infty$ is only a particular case of our general approach: indeed, it is proven in the literature that $H_\lambda$ is a Lyapunov function and that the entropy can be bounded in terms of $\mathcal F$ and of $H_\lambda$, which means that one can use the set $\{\rho\,:\, H_\lambda(\rho)\leq C\coloneq H_\lambda(\rho_0)\}$ as a compact set $K$. To the best of our knowledge, global existence in the 2D critical case without assumptions on the behavior at infinity of the initial datum (finite second moment, or closedness in some sense to some precise stationary profiles with infinite second moment) are difficult to find and we believe that our contribution could help in understand better the role of the key assumptions.


\paragraph{Smoothing of the initial condition.}

Our result also shows that the solution becomes instantaneously bounded and has second-order lower bounds independently of the initial datum, provided the mass is small enough. This is not surprising, since instantaneous smoothing for small mass can also be obtained with completely different methods. In particular, if one differentiates in time the entropy of the solution we observe that we have

\[
\frac{\mathrm{d}}{\mathrm{d}t}
  \int_{\R^d}\rho\,\log\rho
  =-m
\int_{\R^d}\frac{\lvert\nabla\rho\rvert^2}{\,\rho^{2/d}\,}
  +
\int_{\R^d}\rho^2 .
\]
We then use the Gagliardo--Nirenberg
inequality 
$$\|f\|_{L^{\frac{2d}{d-1}}}\leq C\|f \|_{L^{\frac{d}{d-1}}}^{1/d}\|\nabla f\|_{L^2}^{1-1/d}$$
applied to $f=\rho^{1-1/d}$, so that we obtain
$$\int_{\R^d}\rho^2\leq C\left(\int_{\R^d} \rho\right)^{2/d}\int_{\R^d}\frac{\lvert\nabla\rho\rvert^2}{\,\rho^{2/d}\,}.$$
This shows that if the mass $M=\int_{\R^d} \rho$ is small enough then we obtain for some $c_0>0$
\[
\frac{\mathrm{d}}{\mathrm{d}t}
  \int_{\R^d}\rho\,\log\rho
  \le 
- c_0\!\int_{\R^d}\rho^2.
\]
Next, we apply Jensen's inequality to $\rho^2$ in order to compare it to the entropy:
\[
\int_{\R^2}\rho^2
   \ge 
M\,\exp\! \left(\frac1M\!\int_{\R^2}\!\rho\,\log\rho\right).
\]
Hence, if we set 
$E(t) \coloneqq \frac1M\int_{\R^2}\rho(t,x)\,\log\rho(t,x)\,\mathrm{d}x$, we have
\[
E'(t)  \le  - c_0M\exp\! \left(\frac 1M E(t)\right).
\]
Thereby, this ODE-in-time yields the bound
\[
E(t)
  \le  -M
\log\! \left(c_0 t\right),
\]
 which shows a uniform bound on the entropy of the solution for $t>0$, thus preventing blow-up.  As usual, we begin by approximating the initial measure \(\rho_{0}\) by smooth data, but all the bounds above turn out to be independent of the chosen approximation.

\bigskip
Here, the novelty of our approach allows one to obtain a similar result, but with much stronger bounds than just finiteness of the entropy: we obtain $L^\infty$ and second-order bounds. The price to pay is a stronger requirement on the smallness of the mass, but as soon as the mass needs to be small compared to $M_c(d)$, we do not believe that this is crucial. In particular, for the method introduced above, the smallness of the mass depends on the corresponding Gagliardo-Nirenberg inequality. The explicit value of the constant is not known. However, we can compute it numerically. For instance, for dimension $d=2$ (see \cite{Weinstein82})  we observe that the entropy is non-increasing if $M \leq 1,862 \ldots \times (4 \pi)$.





\paragraph{Precise $L^{\infty}$ estimates and applications.} 

Theorems~\ref{thm:Li--Yau} and~\ref{thm: Aronson--B\'enilan} provide an explicit decay estimate of the $L^{\infty}$ norm. Such estimates can be particularly useful in applications. For instance, in~\cite{Bedrossian11}, a similar decay, although derived by a different method and under a different mass threshold, was used to prove the convergence to a self-similar profile. The author remarked that this kind of decay could be obtained with the Aronson–B\'enilan method, but suggested that such an estimate was not yet available for the Keller–Segel system. Here, we confirm that in fact this estimate is actually available, and we also find an explicit threshold in terms of the mass for this to happen. Similar results are also present in a recent preprint, \cite{LisMai}, with yet different techniques.

\paragraph{Minimal mass for the Liouville and the Lane--Emden equation.}

These two equations play a key role in the theory we develop in this manuscript. The Liouville equation
\begin{equation}\label{eq:Liouville}
    \Delta h + e^{h} = 0 \quad \text{in } \R^2
\end{equation}
was first studied by Liouville in \cite{Liouville1853}. In particular, he was the first one to give a representation formula for the equation. More recently, in \cite{Chen_Li91} Chen and Li used the \textit{moving planes} technique in order to show that if $\int_{\R^2} e^h < \infty$ then there exists a unique representation formula for the solutions of the equation \eqref{eq:Liouville}.
There has been further results concerning this problem. For instance, in \cite{Chou_Wan94}, the authors studied the problem when the assumption on the finiteness of the mass is dropped. For further literature review on this problem we recommend the interested reader to consult \cite{Chang04}. In particular, it is a well-known result that any nontrivial solution of \eqref{eq:Liouville} is such that $\int_{\R^2} e^h \geq 8\pi$, see \cite[Lemma 1.1]{Chen_Li91}. We adapt this result in order to show that any nonnegative, nontrivial subsolution of $\Delta h + e^h \geq 0$ also satisfies $\int_{\R^2} e^h \geq 8 \pi$. This result is key in our analysis since we use it in order to recover $L^\infty$ bounds of the solutions for \eqref{eq:KS} on dimension $d=2$ if and only if $\| \rho_0 \|_{L^1(\R^2)} < 8 \pi$. 

On dimension $d > 2$ in order to perform our analysis, we consider the Lane--Emden equation
\begin{equation}\label{eq:Lane Emden}
    \frac{m}{m-1}\Delta h + h^q = 0 \quad \text{in } \Rd
\end{equation}
in the case $q = \frac{d}{d-2}$. In the seminal paper \cite{Gidas_Spruck81}, Gidas and Spruck assert that \eqref{eq:Lane Emden} has no positive solutions. Afterwards, there has been further work focusing on understanding qualitative properties of the Lane--Emden equation \cite{Bahri_Lions92,  Berestycki_Caffarelli_Nirenberg97, Farina07}. In this manuscript we prove that any nonnegative, nontrivial subsolution of $\frac{m}{m-1}\Delta h + h^{\frac{d}{d-2}} \geq 0$ is such that $\int_{\Rd} h^{\frac{d}{d-2}} \geq M_c (d)$, where $M_c (d)$ refers to the critical mass from \Cref{def:critical_mass}. Up to our knowledge, this result is new. The proof relies on an optimal-control argument for a radially decreasing rearrangement of the subsolutions. Analogously to the $2$-dimensional case, this result helps us to understand better the dynamics of \eqref{eq:KS}. It allows us to recover $L^\infty$ estimates of the solutions for \eqref{eq:KS} on dimension $d > 2$ if and only if the mass is subcritical in the sense of \Cref{def:critical_mass}, see \ref{2d_full_mass_multid}.

\paragraph{Technical estimates.}

To prove our main results, we obtained several interesting functional inequalities or PDE estimates that we believe can be useful in other contexts. The first result concerns a general inequality involving the Hessian of the Newtonian potential in two dimensions.

\begin{lemma}[See Remark~\ref{rem:corollary}]
There exists $C>0$ such that for all $\rho\in L^{1}(\R^2)\cap L^{\infty}(\R^2)$ satisfying $\Delta\rho\in L^{1}(\R^2)$, setting $u=\Gamma*\rho$ we have
$$
\|D^2 u\|_{L^{\infty}}\le C \left(\|\rho\|_{L^{\infty}} +\|\Delta\rho\|_{L^{1}}\right).
$$

\end{lemma}

Additionally, in two dimensions, a key tool is the log- Hardy--Littlewood--Sobolev (HLS) inequality, that we reprove via convex analysis tools in order to adapt it to our scopes:

\begin{theorem}[See Theorem~\ref{thm:new_proof_hls}]
Assume that $\rho\ge 0$ is an absolutely continuous measure with finite logarithmic moment such that $\mathcal{F}[\rho]<+\infty$. Take a function $V$ such that $\int_{\R^2}\rho V(x)<+\infty$ and set $K_V \coloneqq \int_{\R^2}e^{-V(x)}$. 
\begin{enumerate}
    \item If the mass $M$ of $\rho$ satisfies $M<8\pi$, then there exists a constant $C=C(M,K_V)$ such that
    $$
        \int_{\R^2}\rho\log\rho\diff x\leq C \left(\mathcal{F}[\rho]+\int V\diff\rho \right).
    $$
    \item If $M=8\pi$, since $\rho$ is not a  Dirac mass,  there exist $\varepsilon>0$, and $R>0$ such that for all $x\in\R^2$ we have $\rho(
B_{R}(x))\leq 8\pi-\varepsilon$; take a constant $\kappa$ such that  $2 - \f{\varepsilon}{4\pi}<\kappa<2$. Then, there exists a constant $C=C(\kappa)$  such that
$$
\left(1-\frac{\kappa}{2}\right)\int_{\R^2}\rho\log\rho\diff x +  C(\kappa)(1+\log R)\le \mathcal{F}[\rho].
$$
\end{enumerate}
Moreover, the constants in the inequalities above are uniform on any compact set (for the narrow convergence) which does not contain Dirac masses.
\end{theorem}

While the first part of the above statement (i.e. the case of masses strictly below $8\pi$) is well-known, the second is diifferent from what is usually done in the literature. Indeed, it is usually said that $\mathcal F$ is bounded from below when $M=8\pi$, but the above result proves that it can also give a bound on the entropy, even at the critical mass $8\pi$, provided the distribution does not concentrate as a single Dirac delta (this is included in $\mathcal{F}[\rho]<+\infty$). This can be considered as a quantitative version of a similar result in \cite{KarWol}, which also proves that the entropy can be bounded by the free energy exploiting an assumption which corresponds to being far from single Dirac masses.

In dimensions larger than two, we also derive a related result at the critical mass analogous to the HLS (or Lane--Emden) inequality. In order to obtain this control in the critical regime we take advantage of a recent stability result by Carlen, Lewin, Lieb, and Seiringer \cite{Carlen_Lewin_Lieb_Seiringer24}.

\begin{theorem}[See Theorem~\ref{thm:new_proof_hls_multid}]
Let $d>2$ and set $m=2-\frac{2}{d}$. Consider a compact set $K\subset\mathcal{M}(\R^d)$ (compact with respect to narrow convergence) consisting of nonnegative measures with total mass $M_c$ that exclude Dirac masses. Then, given an arbitrary positive constant $C_1$, there exists a constant $C=C(K)>0$ such that
\[
\sup\left\{\|\rho\|_{L^m}: \rho\in K, \,  \mathcal{F}[\rho]\le C_1\right\}\le C.
\]
\end{theorem}

Finally, we insist on the intermediate results that we had to develop to justify the use of the comparison principle that we employ to obtain the Aronson--B\'enilan estimate, we relied on a comparison principle described in Section~\ref{sec:justification}.
To justify the existence of the minimum of the Laplacian of $v$ and the fact that it is attained on a compact set, we prove explicit polynomial decays of $\rho$ and of the spatial derivatives of the pressure $p$ (defined as $p=\log\rho$ if $d=2$ or $p=\frac{m}{m-1}\rho^{m-1}$ if $d>2$). More precisely, the statement that we prove can be taken from Propositions \ref{prop:upper_lower_bound_pressure} and \ref{prop:lipschitz-porous} and informally reads as follows:

{\it If $\rho$ is a bounded smooth solution on $[0,T]\times\R^d$ of the Keller-Segel system for some $T>0$ and $\rho_0(x)$ behaves at infinity as $ |x|^{-\beta}$ for an exponent $\beta>d$ (and $\nabla p_0$ and $D^2p_0$ behave as negative powers with exponents close to $-(d-1)$ and $-d$, respectively), then the same asymptotic behaviors, with possible very small changes in the exponents and changes in the multiplicative constants, stay true for $t>0$.}

\subsection{Notations and functional settings}

We denote by $L^{p}(\R^d)$, $W^{m,p}(\R^d)$ the usual Lebesgue and Sobolev spaces, and by $\|\cdot\|_{L^{p}}$, $\|\cdot\|_{W^{m,p}}$ their corresponding norms. As usual, $H^{s}(\R^d) = W^{s,2}(\R^d)$. The spaces $C^{0,\alpha}(\R^d)$ with $0<\alpha<1$ are the usual H\"older spaces, with corresponding norm $\|\cdot\|_{C^{0,\alpha}}$. When the norms concern a vector, or a matrix, we still write $\|\nabla u\|_{L^{p}}$ or $\|D^{2} u\|_{L^{p}}$ instead of $\|\nabla u\|_{L^{p}(\R^d; \R^d)}$ and $\|D^2 u\|_{L^{p}(\R^d;\R^{d^2})}$. We often write $C$ for a generic constant appearing in the different inequalities. Its value can change from one line to another, and its dependence to other constants can be specified by writing $C(a,...)$ if it depends on the parameter $a$ and other parameters. $A\lesssim B$ means that there exists a universal constant $C$ such that $A\le C B$. 

A key tool in our analysis will be the use of the function $v$ defined as follows: given a function $\rho\in L^1_+(\R^d)$, we define
\begin{equation}\label{defiv}
v(x) =
\begin{cases}
\log\rho(x) - u(x), & d=2,\\[3pt]
\displaystyle \f{m}{m-1}\,\rho^{m-1}(x) - u(x), \quad & d>2,\quad\left(m=2-\f{2}{d}\right),
\end{cases}
\end{equation}
where \(
u  = \Gamma\ast\rho,
\) is the Newtionian potential associated with $\rho$, solving $-\Delta u=\rho$.
We will often write $v[\rho]$ to denote the function $v$ defined in terms of $\rho$ as above, and $\delta[\rho]$ for its infimum.

\subsection{Structure of the paper}

In Section~\ref{sec:minimum_laplacian}, we derive a differential inequality on the Laplacian of the pressure. This inequality is at the basis of the Li--Yau and Aronson--B\'enilan estimates and must be treated carefully in two different regimes: $d=2$ and $d>2$. The reason for this split is that the diffusion is linear in the two-dimensional case, i.e. $\Delta\rho$, while it becomes a nonlinear diffusion term $\Delta \rho^m$ in higher dimensions. Compared to the classical heat and porous medium equations, an extra term is obtained in the differential inequality that depends on the $\|\rho\|_{L^\infty}$  norm (although not clearly at this stage). Therefore we show in Section~\ref{sec:Linfty} some estimates on $\|\rho\|_{L^{\infty}}$. Our argument depends both on the mass and on the dimension. The proof is based on properties of the Liouville and Lane--Emden equations (depending on whether $d=2$ or $d>2$ respectively), that we prove in  Section~\ref{sec:liouville}. Section~\ref{sec:d11-d22} shows precisely how the $\|\rho\|_{L^\infty}$ appears in the differential inequality from Section~\ref{sec:minimum_laplacian}. In Section~\ref{sec:endproof} we complete the proof of the Li--Yau and Aronson--B\'enilan estimates taking advantage of all the previous results.  Finally, Section~\ref{sec:justification} is dedicated to a rigorous justification of the formal steps used to prove the differential inequality (in particular the use of the maximum principle). 

\section{Differential inequality for the Laplacian of the pressure}\label{sec:minimum_laplacian}

In this section, we focus on the function $v=v[\rho_t]$. We compute the evolution in time of the minimal value of $\Delta v$  under the Keller--Segel flow at the critical exponent. Such an approach yields a Li--Yau or Aronson--B\'enilan type estimate, respectively for linear and nonlinear diffusion. We look at a point $x_0=x_{0}(t)$ that achieves the minimum of $\Delta v(t,\cdot)$ and we compute its time derivative $\f{\diff}{\diff t} \Delta v(t,x_0)$. 

For a given function $u\in C^2(\R^d)$ we denote by $Q(u)$ the following quantity
$$Q(u)\coloneqq\sup\{|D^2u\, e\cdot e-D^2 u \, e'\cdot e'|\;:\; e,e'\in\R^d,\, |e|=|e'|=1,e\cdot e'=0\}.$$
This quantity is equal to the maximal difference between two eigenvalues of the Hessian of $u$.

We prove the following proposition.
\begin{proposition}\label{prop:preliminary_aronson_benilan}
Set $v(t,x)\coloneqq v[\rho_t](x)$
and let $x_0 = x_0(t)$ be a point achieving the minimum of $\Delta v(t,\cdot)$. Then, if $\rho$ solves \eqref{eq:KS}, we have
%
%
\[
\f{\diff}{\diff t} \Delta v(t,x_0) 
  \ge \left(\Delta v(t,x_0)\right)^2 -\f{(d-1)^2}{2d} Q(u)^2.
\]
\end{proposition}

\begin{remark}
Strictly speaking, the infimum of $\Delta v$ may not be reached. Thus, part of the argument is formal, also because it is only presented for smooth solutions. However, it can be made rigorous. There are many difficulties to do this, and therefore we refer the reader to Section~\ref{sec:endproof} and Section~\ref{sec:justification}  for further details.
\end{remark}

 We break the proof of Proposition~\ref{prop:preliminary_aronson_benilan} into several steps:
\begin{itemize}
    \item In Proposition~\ref{lem:max_laplace_v}, we compute the evolution of $\Delta v$ at the minimum $x_0$. Here, a key observation is that some terms factorise conveniently.
    \item In Proposition~\ref{prop:estimate_scalar_product}, we show how these terms can be controlled by $Q(u)$.
    \item Finally, we combine the above estimates and apply a Young inequality to conclude the proof of Proposition~\ref{prop:preliminary_aronson_benilan}.
\end{itemize}
We present the first step.

\begin{proposition}\label{lem:max_laplace_v}
Considering $v(t,x) \coloneqq v[\rho_t](x)$ and assuming that $\rho$ solves \eqref{eq:KS}, we have 
\begin{multline*}
\partial_t \Delta v 
  =  
  m\rho^{\,m-1}\,\Delta(\Delta v) 
    +  \nabla (\Delta v)\, \cdot \, \left[2m\nabla v + (2m-1)\nabla u\right]
    +2\,\left|D^2 v\right|^2
  + (m-1)\,\left|\Delta v\right|^2\\ 
+2\,\left(D^2 v  -  \f1d\,\Delta v\,\mathrm{Id}\right)
  :\left(D^2 u -  \f1d\,\Delta u\,\mathrm{Id}\right).
\end{multline*}
If $x_0=x_{0}(t)$ is a minimum point of $\Delta v(t,\cdot)$, then at $(t,x_0)$, 
\begin{multline*}
\partial_t  \Delta v(t,x_0)
  \ge 
2\,\left|D^2 v(t,x_0)\right|^2
  + (m-1)\,\left|\Delta v(t,x_0)\right|^2\\
  + 
2\,\left(D^2 v(t,x_0)  -  \f1d\,\Delta v(t,x_0)\,\mathrm{Id}\right)
  :\left(D^2 u(t,x_0)  -  \f1d\,\Delta u(t,x_0)\,\mathrm{Id}\right).
\end{multline*}
\end{proposition}
Here $A:B$ stands for the scalar product between matrices $A:B:=\mathrm{Tr}(A^tB)$ and $|\cdot|$ is the Frobenius norm $|A|\coloneqq\sqrt{A:A}$.
\begin{proof} 
We divide the proof in several steps. 
\begin{enumeratesteps}
\step[Dimension $d=2$]
Using $-\,\Delta u = \rho$, we obtain 
\begin{align*}
\partial_t \Delta v 
&=  
\Delta\left(\partial_t \log\rho\right) 
 + \partial_t \rho = \Delta\left(\f{\partial_t \rho}{\,\rho\,}\right)  + \nabla \rho\cdot \nabla v  + \rho\,\Delta v\\
&= 
\Delta\left(\nabla\log\rho\,\cdot\,\nabla v + \Delta v\right)
 + \nabla\rho\cdot\nabla v 
 + \rho\,\Delta v\\
&= 
\nabla(\Delta\log\rho)\cdot\nabla v 
 + 2\,D^2\log\rho : D^2 v 
 + \nabla(\Delta v)\cdot\nabla\log\rho 
 + \Delta^2 v 
 + \nabla\rho\cdot\nabla v
 + \rho\,\Delta v.
\end{align*}
Using $\log\rho = u + v$, $-\Delta u =\rho$, reorganizing terms, and using the algebraic identity  
\[
2\,\left(D^2 v - \f1d\,\Delta v\,\mathrm{Id}\right) : 
\left(D^2 u - \f1d\,\Delta u\,\mathrm{Id}\right)
 = 
2\,\left(D^2 v : D^2 u\right) 
 + \f{2}{d}\rho\,\Delta v.
\]
we obtain the result in dimension $d=2$. 

\step[Dimension $d>2$]
For the case of dimension $d> 2$, we have to take into account nonlinear diffusion. We compute and obtain
\begin{align*}
\partial_t \Delta v 
=&  
\f{m}{m-1}\Delta\left(\partial_t \rho^{m-1}\right) 
 + \partial_t \rho = m\Delta\left(\rho^{m-2}\partial_t \rho\right)  + \nabla \rho\cdot \nabla v  + \rho\,\Delta v\\
=& 
m\Delta\left(\rho^{m-2}\nabla\rho\cdot\nabla v + \rho^{m-1}\Delta v\right)
 + \nabla\rho\cdot\nabla v 
 + \rho\,\Delta v\\
= &
\f{m}{m-1}\nabla\Delta\rho^{m-1}\cdot \nabla v
+ \f{2m}{m-1}D^2 \rho^{m-1} : D^2 v
+ \f{m}{m-1}\nabla\Delta v : \nabla \rho^{m-1} \\
&+ m \Delta \rho^{m-1}\Delta v
+ 2m\nabla\rho^{m-1}\cdot\nabla\Delta v + m\rho^{m-1}\Delta^2 v + \nabla\rho\cdot\nabla v
+\rho\Delta v.
\end{align*}
Analogously to the previous case we take into account  $\f{m}{m-1}\rho^{m-1}= u + v$ and  $-\Delta u =\rho$. Thus, reorganizing terms, we obtain the stated identity in the case of dimension $d>2$. 

\step[The differential inequality at $x_0$]
Let $x_0=x_{0}(t)$ be a minimum of $\Delta v(t,\cdot)$. Then, it follows  $\nabla (\Delta v)(x_0)=0$ and $D^2(\Delta v)(x_0)\ge 0$ in the sense of matrices (in particular the trace satisfies $\Delta(\Delta v)(x_0)\ge 0$). We obtain in either cases 
\begin{align*}
    \partial_t (\Delta v)(x_0)
 & \ge 
2\,\left|D^2 v(x_0)\right|^2
+ (m-1)\left|\Delta v(x_0)\right|^2 \\
& \quad +2\,\left(D^2 v(t,x_0)  -  \f1d\,\Delta v(t,x_0)\,\mathrm{Id}\right)
  :\left(D^2 u(t,x_0)  -  \f1d\,\Delta u(t,x_0)\,\mathrm{Id}\right). \qedhere
\end{align*}
\end{enumeratesteps}
\end{proof}

From the result obtained in \Cref{lem:max_laplace_v} we observe that we need to estimate  $ A_u (x_0) = \left(D^2 u(x_0) - \f{\Delta u(x_0)}{d} \,\mathrm{Id}\right)$. For instance, in dimension $2$, this matrix is written

\[
A_u  = 
\begin{pmatrix}
\f{\partial_{11} u - \partial_{22} u}{2}
&
\partial_{12} u 
\\[6pt]
\partial_{21} u 
&
\f{\partial_{22} u - \partial_{11} u}{2}
\end{pmatrix}.
\]

\begin{proposition}\label{prop:estimate_scalar_product}
Given $A_v \coloneqq D^2 v - \f1d\Delta v\,\mathrm{Id}$ and $A_u \coloneqq D^2 u - \f1d \Delta u\,\mathrm{Id}$, we have
\[
\left|\,A_v(x_0):A_u(x_0)\right|
  \le 
\frac{d-1}{\sqrt{d}}|A_v(x_0)|
Q(u)
\]
Moreover,
\[
\left|D^2 v(x_0)\right|^2 
  = 
\f1d\left(\Delta v(x_0)\right)^2 
  + 
\left|A_v(x_0)\right|^2.
\]

\end{proposition}



\begin{proof}
Up to choosing a suitable orthonormal basis, we can assume that $D^2 v(x_0)$ (which is symmetric) is diagonal, and $A_v(x_0)$ is diagonal as well. On the same basis, $A_u(x_0)$ is maybe not diagonal, but only its diagonal terms appear in the product $A_v(x_0):A_u(x_0)$. Moreover, each of the terms on the diagonal of $A_u(x_0)$ is bounded by $\frac{d-1}{d}Q(u)$. Indeed, if we consider $u_{ii}$, the diagonal terms of $D^2u(x_0)$, we have
$$
    (A_u(x_0))_{ii}=u_{ii}-\frac1d\sum_j u_{jj}=\frac{1}{d}\sum_{j\neq i} (u_{ii}-u_{jj}).
$$  
Hence, if we now consider $a_{ii}$, the terms on the diagonal of $A_u(x_0)$, we have
$$
    |A_v(x_0):A_u(x_0)|\leq \frac{d-1}{d}Q(u)\sum_i|a_{ii}|\leq \frac{d-1}{d}Q(u)\sqrt{d}\left(\sum_i|a_{ii}|^2\right)^{1/2}
$$
and the first part of the statement follows from  $A_v(x_0)=\left(\sum_i|a_{ii}|^2\right)^{1/2}$.
Lastly, the identity
\[
\left|D^2 v(x_0)\right|^2 
   =  
\f1d\,\left(\Delta v(x_0)\right)^2 
  + 
\left|\,D^2 v(x_0)  -  \f1d(\Delta v(x_0))\,\mathrm{Id}\right|^2
\]
comes from the orthogonal decomposition of any real matrix $M$ into a multiple of the identity and a trace-free part:
\[
M 
 = 
\f{\mathrm{tr}(M)}{d}\,\mathrm{Id} 
 + 
\left(M - \f{\mathrm{tr}(M)}{d}\,\mathrm{Id}\right)
\quad
\text{and thus}
\quad
|M|^2 
 = 
\f{(\mathrm{tr}(M))^2}{d} 
 + 
\left|M - \f{\mathrm{tr}(M)}{d}\,\mathrm{Id}\right|^2.\qedhere
\]
\end{proof}

We are now in a position to combine the results of Propositions~\ref{lem:max_laplace_v} and~\ref{prop:estimate_scalar_product}.

\begin{proof}[Proof of Proposition~\ref{prop:preliminary_aronson_benilan}]
The proof is an immediate consequence of Proposition~\ref{lem:max_laplace_v},    Proposition~\ref{prop:estimate_scalar_product} and the Young's inequality as
\begin{align*}
2\left|A_v(x_0):A_u(x_0)\right| &\le  2 \left|A_v(x_0)\right|^2 + \f{(d-1)^2}{2d}Q(u)^2.\qedhere
\end{align*}
\end{proof}

 With Proposition~\ref{prop:preliminary_aronson_benilan} in hand, the next step is to control $Q(u)$. Due to the definition of $Q(u)$, one way to do obtain this control is with $L^\infty$ estimates on some second derivatives of $u$. Since $u$ is related to $\rho$ through the Poisson equation, i.e. $-\Delta u=\rho$, we can obtain estimates on the $L^\infty$ norm of $\rho$ instead. Yet, this is not enough as the elliptic regularity fails in $L^\infty$ and $\rho\in L^\infty$ does not imply $D^2u\in L^\infty$. We devote the following two sections to this problem.

In \Cref{sec:Linfty}, we first derive an estimate on $\|\rho_t\|_{L^\infty}$ in terms of $\min(\Delta v(t,\cdot))$. We prove this bound in a \emph{static} setting, i.e.\ without the time dependence of $\rho$. Finally, in \Cref{sec:d11-d22} we study in detail the $L^\infty$ estimate on the differences of second derivatives of $u$ appearing in $Q(u)$.

\section{\texorpdfstring{$L^{\infty}$}{Linfty} estimates on the density}\label{sec:Linfty}

Our goal in this section is to derive bounds on $\|\rho\|_{L^\infty}$ in terms of the quantity
\begin{equation}\label{eq:def_delta}
\delta  \coloneqq  \inf_{x \in \mathbb{R}^d}\,\Delta v(x).
\end{equation}
We will often write $\delta[\rho]$ for $\inf \Delta v[\rho]$. 
We divide the analysis into three cases.
\begin{enumerate}
\item Small mass: There exists an explicitly computable threshold
such that, below this threshold, $\|\rho\|_{L^\infty}\le C_2(M)|\delta|$
with an explicit constant $C_2(M)$.

\item Subcritical mass: If $0<M<M^{sub}_c$, we still have
$\|\rho\|_{L^\infty}\le C_2(M)|\delta|,$ but the finite constant $C_2(M)$ is obtained by compactness and is
not explicit. We recall that $M^{sub}_c$ is defined
in~\eqref{def:Mc LaneEmden}.

\item Critical mass: When $M = M^{sub}_c$, we still manage to control $\rho$ in the $L^\infty$ norm. However, in this case, a bound of the form $\|\rho\|_{L^\infty}\le C|\delta|$ is not possible, as there exist non-zero functions for which $\delta=0$. Instead, we will prove a bound of the form  \text{$\|\rho\|_{L^\infty} \le C(1 + |\delta|)$}  whenever we impose $\rho\in K$, where the set $K$ is a suitable compact set of measures not containing any Dirac mass. In this case, $C$ cannot be computed explicitly either.
\end{enumerate}

Although this section is concerned with functional analysis results rather than the Keller--Segel system, the system’s critical mass appears in a surprisingly natural way. As we will see in the proof, this is because our analysis is based on particular properties of some elliptic PDEs: the Liouville equation in two dimensions and the Lane-Emden equation in higher dimensions. These equations happen to be associated with the minimizers of the log-HLS and HLS inequalities, respectively. We divide this section into three different subsections, one for each of the three cases discussed above.

\subsection{Small mass}

In this subsection, we assume the total mass $M = \|\rho\|_{L^1(\mathbb{R}^d)}$ is sufficiently small (smaller than some threshold explicitly computed). Under this condition, we can use properties of subharmonic functions to derive a direct inequality connecting $\|\rho\|_{L^\infty}$ and $\delta$. The main result of this subsection is as follows:

\begin{proposition}\label{prop:linf_very_small_mass}
Assume $\Delta v \ge \delta$ on $\mathbb{R}^d$. Let $M=\|\rho\|_{L^{1}(\mathbb{R}^d)}$. Hence, we have that:
\begin{enumerate}
\item If $d=2$ and $M < 4\pi$ then
$$
    \|\rho\|_{L^{\infty}}\le C_2(M)|\delta|,
    \qquad\text{where}\qquad
    C_2(M) = \f{M}{4\pi-M}.
    $$
\item If $d>2$ and $M < \left(\f{4d(d-1)}{d-2}\right)^{d/2}I(d)$,
where $I(d)$ is defined in \Cref{lem:Naive Linfty}, then
$$
\begin{aligned}
\|\rho\|_{L^{\infty}} &\le C_2(M)|\delta|,
\qquad
C_2(M) = \f{\theta_d C_1(M)}{1-\theta_d C_1(M)},\\
\theta_d C_1(M) &= \f{d-2}{4d(d-1)}
\left(\f{M}{I(d)}\right)^{2/d}.
\end{aligned}
$$
\end{enumerate}
In both cases the condition on $M$ is exactly $\theta_d C_1(M)<1$, with $\theta_d=1$ for $d=2$ and $\theta_d=\f{m-1}{m}$ for $d>2$.
\end{proposition}
Note that the necessity to use the constant $\theta_d$ comes from the choice of the normalization constants depending on $m$ in the definition of the pressure for $d>2$.

In order to prove this result, we need an intermediate step that we state in the  following lemma. This is a quite classical result, which can be rarely found stated as such. It is usually stated as an $L^\infty$ bound in terms of the $L^1$ norm and on lower bounds on the Laplacian of the logarithm in arbitrary dimension, but here we only use this formulation in the case $d=2$ and we adapt the proof to the critical power nonlinearity for $d>2$.
 
\begin{lemma}\label{lem:Naive Linfty} We distinguish two cases:
\begin{enumerate}
    \item If $d=2$ and $\Delta \log \rho \geq -\alpha$ on $\R^2$, then,
    $$
    \|\rho\|_{L^\infty} \leq C_1(M)\alpha,
    \qquad\text{where}\qquad
    C_1(M) = \f{M}{4\pi}.
    $$
    \item If $d>2$ and $\Delta \rho^{1-2/d} \geq -\alpha$ on $\R^d$, then,
    $$
    \|\rho\|_{L^\infty} \leq C_1(M)\alpha,
    \qquad\text{where}\qquad
    C_1(M) = \f{1}{2d}\left(\f{M}{I(d)}\right)^{2/d},
    $$
    and $I(d)=\int_{B_1} (1-|y|^2)^{d/(d-2)}dy$.
\end{enumerate}
\end{lemma}
\begin{remark}
    This estimate in the logarithmic case (generalized to arbitrary dimension) can be deduced from \cite[Theorem 1.16]{dePhilippis_Shenfeld25},  or found in the case $\alpha=1$ in \cite[Lemma 3.7]{GraKemLoeZak}, where the idea of integrating on spheres is presented, as well as in \cite{FanchPhD}. The extension to the critical exponent in higher dimension is an adaptation. We observe that the coefficient $I(d)$ can be computed in terms of the special function $\Gamma$, and in some cases admits more explicit expressions For instance, $I(3)=3\pi/14\approx \f{3}{14}3.14$ (note the double appearance of the numbers 3 and 14)

Let us also mention that these inequalities are sharp and they are attained by Gaussian profiles for $d=2$ ($\rho(x)=Ae^{-c|x|^2}$) and by Barenblatt profiles for $d>2$ ($\rho(x)=(A-c|x|^2)_+^{d/(d-2)}$).

Moreover, similar sharp estimates can be obtained if replacing $\rho^{1-2/d}$ ($\log\rho$ for $d=2$) with other powers of $\rho$, but only for this choice of the exponent the estimate is linear in $\alpha$, which is important in our analysis.
    
\end{remark}

\begin{proof}
We start with the two-dimensional case.
    From the definition of the Laplacian, we have
    \begin{equation*}
        \Delta \left( \log \rho + \f{\alpha}{4} |x|^2 \right) \geq 0.
    \end{equation*}
    Hence, for any $x_0$ and $r > 0$, using the properties of subharmonic functions and Jensen's inequality, we get

    \begin{equation*}
        \log \rho (x_0) \leq \dashint_{\partial B_r(x_0) } \log \rho + \f{\alpha}{4} r^2 \leq \log \left( \dashint_{\partial B_r (x_0)} \rho \right) + \f{\alpha}{4} r^2.
    \end{equation*}
    Hence, $\rho(x_0) e^{-\frac{\alpha}{4} r^2} \leq \dashint_{\partial B_r (x_0)} \rho$. We now multiply both sides times the $(d-1)-$dimensional measure of $\partial B_r$ and integrate in $r$, so that we obtain
   $$
         \rho (x_0) \int_\Rd e^{-\frac{\alpha}{4} |x-x_0|^2} \dx \leq \int_\Rd   \rho = M .
    $$
    To conclude we just need to compute the integral 
    $\int_\Rd e^{-\frac{\alpha}{4} |x-x_0|^2} \dx = \f{4\pi }{\alpha}$.

   We now consider the case of dimension $d > 2$.
We use 
\begin{equation*}
        \Delta \left(\rho^{1-2/d}+ \f{\alpha}{2d} |x|^2 \right) \geq 0.
\end{equation*}
We then obtain
    \begin{equation*}
        \rho (x_0)^{1-2/d} \leq \dashint_{\partial B_r(x_0) } \rho^{1-2/d} + \f{\alpha}{2d} r^2 \leq \left( \dashint_{\partial B_r (x_0)} \rho \right)^{1-2/d} +\f{\alpha}{2d} r^2.
    \end{equation*}
    We then deduce 
    $$\left(\rho(x_0)^{1-2/d}-\f{\alpha}{2d} r^2\right)_+^{\f{d}{d-2}}\leq \dashint_{\partial B_r (x_0)} \rho.$$ 
    By integrating as in the previous case we obtain
   $$
       \int_{\R^d}\left(\rho(x_0)^{1-2/d}-\f{\alpha}{2d}|x-x_0|^2\right)_+^{\f{d}{d-2}}dx\leq M .
    $$
We now do a change of variable in the integral, of the form $x=x_0+Ry$ with $R$ chosen such that $R^2\f{\alpha}{2d}=\rho(x_0)^{1-2/d}$. In this way we obtain
 $$\rho(x_0)R^d\int_{\R^d}(1-|y|^2)_+^{\f{d}{d-2}}dy\leq M.$$
 By using the choice of $R$, which gives
 $R^d\left(\f{\alpha}{2d}\right)^{d/2}=\rho(x_0)^{-1+d/2}$, we obtain
$$\left(\f{2d}{\alpha}\rho(x_0)\right)^{d/2}I(d)\le M,$$
 which gives the desired result.
\end{proof}

We now return to the main statement, which considers $\Delta v \ge \delta$.

\begin{proof}[Proof of Proposition~\ref{prop:linf_very_small_mass}]
Recall that $u$ satisfies $-\Delta u = \rho$. In dimension $d=2$ this gives $\Delta \log \rho \ge \delta - \rho$, whereas for $d>2$,
$$
\Delta \rho^{1-2/d} = \f{m-1}{m}\Delta\left(\f{m}{m-1}\rho^{m-1}\right) \ge \f{m-1}{m}(\delta-\rho).
$$
Thus we can apply \Cref{lem:Naive Linfty} with $\alpha = \theta_d(|\delta|+\|\rho\|_{L^\infty})$, where $\theta_d=1$ for $d=2$ and $\theta_d=\f{m-1}{m}$ for $d>2$. We obtain
$$
\|\rho\|_{L^\infty} \le \theta_d C_1(M)\left(|\delta|+\|\rho\|_{L^\infty}\right).
$$
If the mass is smaller than a certain explicit threshold, we have that $\theta_d C_1(M) < 1$. From here, we deduce
$$
\|\rho\|_{L^\infty} \le \f{\theta_d C_1(M)}{1-\theta_d C_1(M)}|\delta| = C_2(M)|\delta|,
$$
where the constant $C_1(M)$ is explicit as discussed in \Cref{lem:Naive Linfty}. This yields the result.
\end{proof}

\subsection{Subcritical mass}\label{subsec:subcritical_Linf}

We now consider the case $0<M<M^{sub}_c(d)$ with $M^{sub}_c(d)$ defined in \eqref{def:Mc LaneEmden}. In this regime, the explicit computations from the small mass case are not enough anymore. However, we can still prove a bound of the form
$$
\|\rho\|_{L^\infty}
    \le C_2(M)|\delta|
$$
with a finite, generally non-explicit, constant $C_2(M)$. The proof
is a contradiction argument based on subsolutions of the Liouville
and Lane--Emden equations.

\begin{proposition}[$L^\infty$-bound on $\rho$ for subcritical mass]\label{prop:linf_estimate}
Let $0<M<M^{sub}_{c}(d)$. Consider all nonnegative $\rho\in L^\infty(\mathbb{R}^d)\cap H^{1}_{loc}(\R^d)$ with
\[
\int_{\mathbb{R}^d}\rho = M,
\quad
\Delta v[\rho] \ge \delta .
\]
Then, there exists a finite constant $C_2(M)$, depending only on
$M$ and $d$, such that they all satisfy
\begin{equation}\label{eq:Subcritical mass Linfty bound}
\|\rho\|_{L^\infty}
    \le C_2(M)|\delta|.
\end{equation}
\end{proposition}

\begin{remark}
In dimension $2$, the explicit small-mass choice
$$
C_2(M)=\frac{M}{4\pi-M}
$$
diverges as $M\to4\pi$ (and a similar phenomenon occurs for $d>2$). This is only a limitation of the
explicit estimate: \Cref{prop:linf_estimate} provides a finite,
generally non-explicit, choice of $C_2(M)$ for every fixed
subcritical mass $M<8\pi$.
\end{remark}

The idea behind the proof of this result is based on a contradiction argument. 
If \Cref{prop:linf_estimate} does not hold then we are able to construct a nonnegative subsolution $\rho$ with mass $M$ such that $\Delta v \geq 0$ and $\| \rho \|_{L^\infty} \leq 1$ of the Liouville equation given by \eqref{eq:Liouville} if the dimension is $d=2$ or the Lane--Emden equation \eqref{eq:Lane Emden} if $d>2$. In the case of dimension $2$, we choose $h  = \log \rho$ in \eqref{eq:Liouville} and we construct a subsolution of the Liouville equation, which does not admit any solutions with mass strictly smaller than $8 \pi$. In the case of higher dimension, we take $h = \rho^{m-1}$ in \eqref{eq:Lane Emden} and we construct a subsolution of a Lane--Emden equation. By definition, this equation does not admit subsolutions with mass smaller than $M^{sub}_c(d)$, which produces the desired contradiction.

We divide the proof into the case of dimension $2$ and higher dimensions.

\begin{proof}[Proof of Proposition~\ref{prop:linf_estimate} for $d=2$] 
We divide the proof into several steps.

\begin{enumeratesteps}
\step[Setting the contradiction]
By the definition of $M^{sub}_c(2)$, it is not possible to have $M<M^{sub}_c(2)$ and $\delta \geq 0$. We can therefore assume that $\delta <  0$. Suppose, by contradiction, that for every integer $n$ there exists a nonnegative function $\rho_n$ with total mass $M$, i.e. $\int_{\mathbb{R}^2}\rho_n = M$, and a number $\delta_n<0$ such that 
\[
    \Delta(\log\rho_n) + \rho_n  \ge  \delta_n
    \quad\text{and}\quad
    \|\rho_n\|_{L^\infty}  \ge  n\,|\delta_n|.
\]
Define the rescaled sequence
\[
    \eta_n(x) = \lambda_n^2\,\rho_n\left(\lambda_n x\right),
    \quad
    \text{where}
    \quad
    \lambda_n^2 
    = 
    \f{1}{\|\rho_n\|_{L^\infty}}.
\]
Note that $\int_{\mathbb{R}^2}\eta_n = M$, and also $\|\eta_n\|_{L^\infty} =1$. Furthermore, since the sequence $\rho_n$ satisfies
\[
\Delta(\log\rho_n) + \rho_n  \ge  \delta_n,
\quad\text{we have that}\quad
\Delta(\log \eta_n) + \eta_n 
 = 
\lambda_n^2\,(\Delta\left(\log \rho_n(\lambda_n\cdot)\right) + \rho_n(\lambda_n\cdot) )
 \ge 
\delta_n\,\lambda_n^2.
\]
Thus, since $|\delta_n|\lambda_n^2 \leq \frac1n$, the right hand side converges to $0$ as $n \rightarrow \infty$.

\step[Compactness]
Since each $\eta_n$ has mass $M$ and $\|\eta_n\|_{L^\infty}= 1$, there is a subsequence (not relabeled) such that $\eta_n\rightharpoonup \eta$ weakly-* in $L^\infty(\mathbb{R}^2)$ and weakly in $L^p(\mathbb{R}^2)$ for $1<p<\infty$. Up to translating the function $\eta_n$, we can assume $\eta_n(0)>\frac12$ (since $\sup\eta_n=1$; we could prove that this sup is a max and hence assume $\eta_n(0)=1$, but we do not need it).

Let us remark 
\[
\eta_n\Delta(\log \eta_n)  =  \Delta \eta_n  - \f{|\nabla \eta_n|^2}{\,\eta_n\,}.
\]
Then, we get 
\[
\Delta \eta_n - \f{|\nabla \eta_n|^2}{\,\eta_n\,} + \eta_n^2  \ge  -|\delta_n|\,\lambda_n^2\,\eta_n.
\]
First, we notice that this implies $\Delta \eta_n\geq -2$ for large $n$ (using $\eta_n^2\leq 1$ and $|\delta_n|\lambda_n^2\to 0$). Hence, we have $\fint_{B(0,r)} \eta_n\geq\frac12-Cr^2$ and this inequality passes to the limit as $n\to\infty$. In particular, the limit function $\eta$ also satisfies $\fint_{B(0,r)}\eta>0 $ for small $r$, and it is not identically zero.

Thereby, if we multiply by $\varphi$ for some nonnegative $\varphi\in C_{c}^{\infty}(\R^d)$ we recover 
\[
\int_{\mathbb{R}^2}\f{|\nabla \eta_n|^2}{\,\eta_n\,}\varphi
 \le 
\int_{\mathbb{R}^2}\left(\eta_n^2\varphi + \eta_n\Delta\varphi + |\delta|\,\lambda_n^2\,\eta_n\varphi\right),
\]
which is uniformly bounded since $\int \eta_n=M$ and $\|\eta_n\|_{L^\infty}\le1$. Therefore, $\nabla \eta_n$ is uniformly bounded in $L^2_{loc}(\mathbb{R}^2)$ and thus $\eta_n$ is uniformly bounded in $H^{1}_{loc}(\mathbb{R}^2)$. Thus, up to a subsequence, $\eta_n\to \eta$ a.e.\ and strongly in $L^2_{\mathrm{loc}}(\mathbb{R}^2)$, and weakly in $H^1_{\mathrm{loc}}(\mathbb{R}^2)$. In particular, $\int_{\mathbb{R}^2}\eta\le M$ and $\eta\ge0$.

\step[Passing to the limit in the inequality]
For a nonnegative smooth compactly supported $\varphi$, we multiply the identity
\[
\Delta \eta_n  - \f{|\nabla \eta_n|^2}{\eta_n} + \eta_n^2 
 \ge 
-\delta\,\lambda_n^2\,\eta_n
\]
by $\varphi$, we integrate by parts and we recover 
\[
\int_{\mathbb{R}^2}\Delta \eta_n\,\varphi
 = 
\int_{\mathbb{R}^2}\eta_n^2\,\Delta\varphi
 \to 
\int_{\mathbb{R}^2}\eta^2\,\Delta\varphi,
\]
and
\[
\int_{\mathbb{R}^2}\eta_n\,\varphi
 \to 
\int_{\mathbb{R}^2}\eta\,\varphi,
\quad
\int_{\mathbb{R}^2}\delta\,\lambda_n^2\,\eta_n\,\varphi
 \to 
0.
\]
For the gradient term, the functional 
\[
h \mapsto \int_{\mathbb{R}^2}\f{|\nabla h|^2}{\,h\,}\,\varphi
\]
is lower-semicontinuous for the weak $H^1_{\mathrm{loc}}$ topology. This is the case because the integrand is convex in $\nabla h$ and continuous in $h$ and $x$ since $\varphi$ is smooth and nonnegative. Hence, we have 
\[
-\int_{\mathbb{R}^2}\f{|\nabla \eta|^2}{\,\eta\,}\,\varphi
 \ge 
\limsup_{n\to\infty}\left(-\int_{\mathbb{R}^2}\f{|\nabla \eta_n|^2}{\,\eta_n\,}\,\varphi\right).
\]
Letting $n\to\infty$, we conclude that $\eta$ satisfies
\[
\Delta \eta - \f{|\nabla \eta|^2}{\eta} + \eta^2  \ge  0,
\]
or equivalently $\Delta(\log \eta) + \eta \ge 0$ wherever $\eta>0$. Moreover $\eta\ge 0$, $\eta\in H^{1}_{loc}(\R^2)$, and $\|\eta\|_{L^{\infty}}\le 1$, $\int_{\R^2}\eta\le M.$ 

\step[Subsolution of the Liouville equation]
Let us define $h=\log \eta$, with $h=-\infty$ on $\{\eta=0\}$. Then $\eta=e^h$ and
\[
\eta\,(\Delta(\log \eta ) + \eta)
 = 
e^h\left(\Delta h + e^h\right).
\]
Hence, $h$ is such that 
\[
\Delta h + e^h  \ge  0
\quad
\text{in every level set }\Omega_t=\{\,x : h(x)>t\}.
\]
Moreover $0 <\int_{\mathbb{R}^2} e^h = \int_{\mathbb{R}^2}\eta \le M <  M_c^{sub} (2)$,  getting a contradiction which completes the proof. \qedhere
\end{enumeratesteps}
\end{proof}

Now, we proceed to prove the statement also for dimension $d > 2$ analogously.

\begin{proof}[Proof of Proposition~\ref{prop:linf_estimate} for $d>2$]
We adapt the strategy from the two-dimensional case. 

\begin{enumeratesteps}
\step[Setting the contradiction]
Suppose by contradiction that the statement of the proposition is false. Then, for every integer $n\in\mathbb{N}$, there exists a nonnegative function $\rho_n \in L^\infty(\mathbb{R}^d)$ with mass $\int_{\mathbb{R}^d}\rho_n = M$ and $\delta_n<0$ such that 
\[
\f{m}{m-1}\Delta \rho_n^{m-1} + \rho_n  \ge  \delta_n  \quad\text{and}\quad \|\rho_n\|_{L^\infty}  \ge  n \,|\delta_n|.
\]
Let us define the rescaled sequence
\[ 
\eta_n(x)  =  \lambda_n^d \,\rho_n\left(\lambda_n x\right), \quad \text{where}\quad \lambda_n^d  = \f{1}{\|\rho_n\|_{L^\infty}}.
\]
Then, we get that 
\[
\begin{aligned}
\f{m}{m-1}\Delta \eta_n^{m-1} + \eta_n
&= \lambda_n^d \left(
\f{m}{m-1}\Delta \rho_n(\lambda_n \cdot)^{m-1}
+ \rho_n(\lambda_n \cdot)\right) \ge -\frac1n,\\
\|\eta_n\|_{L^\infty} &= 1,
\qquad
\int_{\mathbb{R}^d} \eta_n = M.
\end{aligned}
\]
Analogously to the case $d=2$, we translate $\eta_n$ so that we assume $\eta_n(0)>\frac12$.

    \step[Compactness]
Since $0\le \eta_n \le 1$ and $\int_{\mathbb{R}^d} \eta_n = M$, we get uniform bounds in $L^p$ for all $1 \le p<\infty$. 
Once more, we work analogously to the two-dimensional case. We consider $\varphi$, a nonnegative smooth, compactly supported function. Thus, if we multiply  the equation by $\eta_n^{2-m}\varphi$ we obtain that
\[
\|\nabla \eta_n\|_{L^2_{loc}(\mathbb{R}^d)}
\quad\text{is uniformly bounded in }n.
\]
Hence, up to a subsequence, $\eta_n \to  \eta$ weakly in $H^1_{\mathrm{loc}}(\mathbb{R}^d)$, strongly in all $L^{p}(\R^d)$, and almost everywhere. The same holds for $\eta_n^{m-1}$, which converges to $\eta^{m-1}$. The identification of this limit is possible because of the a.e. convergence of $\eta_n$. We also get $\eta \in H^1_{\mathrm{loc}}(\mathbb{R}^d)$ with $\eta\ge 0$ and $\int_{\mathbb{R}^d} \eta \le M$. Moreover, from $\Delta\eta_n^{m-1}\geq -C$ we obtain  $\fint_{B(0,r)}\eta^{m-1}\geq C$ for small $r$. This inequality passes to the limit and proves that $\eta$ is not identically zero.

\step[Passing to the limit in the inequality]
Let $h_n=\eta_n^{m-1}$ and
$\varepsilon_n=|\delta_n|\lambda_n^d\le1/n$. Then $0\le h_n\le1$
and $\Delta h_n\ge-C$. Testing with $h_n\chi^2$, after
mollification if necessary, gives
$$
\int\chi^2|\nabla h_n|^2
\le2C\int\chi^2+4\int|\nabla\chi|^2.
$$
Thus $h_n\rightharpoonup h$ in $H^1_{\mathrm{loc}}$ and
$h_n\to h$ almost everywhere and strongly in every finite
$L^p_{\mathrm{loc}}$. Hence
$\eta_n\to\eta=h^{1/(m-1)}$ locally strongly, with
$\int\eta\le M$. The mean-value inequality shows that $h\not\equiv0$.
Finally,
$$
\frac{m}{m-1}\int h_n\Delta\varphi+\int\eta_n\varphi
\ge-\varepsilon_n\int\varphi
$$
passes to the limit and yields
$$
\frac{m}{m-1}\Delta h+h^{1/(m-1)}\ge0.
$$

Then $h\ge 0$, $h\in H^1_{\mathrm{loc}}(\mathbb{R}^d)$, and it satisfies
\[
\f{m}{m-1}\Delta h + h^{\f{d}{d-2}}  \ge  0.
\]
Yet, the mass is such that
\[
\int_{\mathbb{R}^d} h^{\f{d}{d-2}}  =  \int_{\mathbb{R}^d} \eta  \le  M .
\]
Then, the value of $M$ must be at least $M^{sub}_c(d)$, but we assumed $M<M^{sub}_c(d)$. This contradiction completes the proof.\qedhere
\end{enumeratesteps}
\end{proof}

\subsection{Critical mass}

Finally, we consider the borderline case $M = M^{sub}_c(d)$. Here, the approach used for subcritical mass fails to yield a bound of the form $C\,|\delta|$ alone, since there do exist nontrivial subsolutions at the critical mass for the Liouville and the Lane-Emden equation. Moreover, if we rescale such a subsolution as it is presented in the previous subsection, we obtain a sequence $\rho_n$ with $\delta[\rho_n]=0$ and $\rho_n\rightharpoonup\delta_0$. Thereby, our contradiction argument fails to work in this case, and we cannot obtain the same $L^\infty$ bound. 
Therefore, we will then need to prove the desired bound under the assumption $\rho\in K$ where $K$ is a compact set for the narrow convergence that does not contain any Dirac mass.

Here, we present the proof for a generic compact set $K$. Later, in \Cref{sec:endproof}, we choose the set $K$ that allows us to study the solutions of the Keller--Segel equation. To start this subsection, we first recall the definition of narrow convergence, tightness and Prokhorov's theorem. 

\begin{definition}[Narrow Convergence of Measures]
Let $\{\mu_n\}_{n=1}^\infty$ and $\mu$ be a (finite) Radon measures on $\mathbb{R}^d$. We say that $\mu_n$ \emph{converges narrowly} to $\mu$ if and only if
\[
\int_{\mathbb{R}^d} \varphi \,d\mu_n 
   \to
\int_{\mathbb{R}^d} \varphi \,d\mu
\quad\text{as }n\to\infty
\]
for every $\varphi\in C_{b}(\R^d)$ the space of bounded and continuous functions $\varphi: \mathbb{R}^d \to \mathbb{R}$. 
\end{definition}

\begin{definition}[Tightness for a family of measures with same mass]
Let $\{\mu_\alpha\}_{\alpha\in A}$ be a family of finite Radon measures on $\mathbb{R}^d$ with mass $M$. We say this family is \emph{tight} if, for every $\varepsilon>0$, there exists a compact set $K_\varepsilon \subset \mathbb{R}^d$ such that
\[
\mu_\alpha\left(K_\varepsilon\right)  \ge M - \varepsilon
\quad
\text{for all }\alpha\in A.
\]
\end{definition}
The intuition behind tightness is that our family of measures cannot ``escape to infinity". The reason is that most of their mass is contained in a compact set that can be chosen uniformly for all the measures of the family. Taking advantage of the notion of tightness one can also introduce Prokhorov's Theorem.

\begin{theorem}[Prokhorov's Theorem]
Let $\{\mu_n\}$ be a sequence of finite Radon measures on $\mathbb{R}^d$. Then, the following statements are equivalent:
\begin{enumerate}
\item The sequence $\{\mu_n\}$ is tight.
\item The set $\{\mu_n\}$ is relatively sequentially compact in the sense of narrow convergence, i.e.\ there exists a subsequence $\{\mu_{n_k}\}$ and a measure $\mu$ such that $\mu_{n_k}$ converges narrowly to $\mu$. 
\end{enumerate}
\end{theorem}

Once we have introduce the relevant notions from measure theory we proceed to introduce the main result of this subsection.

\begin{proposition}[$L^\infty$-bound on $\rho$ for critical mass]\label{prop:linf_estimate_critical}
Let $K$ be a set of positive measures with mass $M^{sub}_c(d)$ given by \eqref{def:Mc LaneEmden}. Let $K$ be compact for the narrow convergence and such that it does not contain any  Dirac mass, i.e. measures of the form $M^{sub}_c(d)\delta_a$ for $a\in \R^d$.  Then, there exists a constant $C=C(K)>0$ such that for all $\rho\in K\cap L^\infty(\mathbb{R}^d)\cap H^{1}_{loc}(\R^d)$ with $
\Delta v \ge \delta$ 
it follows that
\[
\|\rho\|_{L^\infty}  \le C( |\delta|+1).
\]
\end{proposition}

In order to prove this proposition, we need a further result on measure theory. 
Let us consider a measure sequence for which its rescalings are also tight.
Then, in the following lemma, we show that under the appropriate tightness assumptions, this measure sequence must converge to a single Dirac mass. 
Let us recall that $\nu=T\sharp\mu$ is called the pushforward of the measure $\mu$ under the map $T$.

\begin{lemma}\label{lem:contradiction_ Dirac_mass}
Let $\rho_n$ be a sequence of nonnegative Radon measures on $\mathbb{R}^d$ with the same total mass $M>0$. Assume that $\rho_n$ is tight and that there exist $\lambda_n\to+\infty$ and $x_n\in\R^d$ such that the sequence
\[
\tilde{\rho}_n  \coloneqq  T_n\sharp\rho_n
\quad \text{with} \quad
T_n(x)=\lambda_n\,(x - x_n)
\]
is also tight. Then, there exists a subsequence (not relabeled) and a point $x_\infty\in\R^d$ such that
\[
\rho_n  \rightharpoonup  M\,\delta_{x_\infty}
\quad
\text{narrowly in $\mathcal{M}(\R^d)$}.
\]
\end{lemma}

\begin{proof}[Proof of Lemma~\ref{lem:contradiction_ Dirac_mass}]
Since $\rho_n$ and $\tilde{\rho}_n$ are both tight, by Prokhorov’s theorem, there exist subsequences (not relabeled) such that
\[
\rho_n  \rightharpoonup \rho, 
\quad
\tilde{\rho}_n  \rightharpoonup \tilde{\rho}
\quad
\text{narrowly in } \mathcal{M}(\R^d),
\]
with $\rho,\tilde{\rho}$ having total mass $M$. We claim that $\rho$ must be a  Dirac mass. 
Let us consider a test function $\varphi\in C_b(\R^d)$. Then,  we note that
\[
\int \varphi\,\mathrm{d}\rho_n
  = 
\int \varphi\circ T_n^{-1}\,\mathrm{d}\tilde{\rho}_n
  = 
\int \varphi\!\left(x_n + \f{x}{\lambda_n}\right)\,\mathrm{d}\tilde{\rho}_n(x).
\]
Let us now define a smooth cutoff function $\chi$ such that $\chi=1$ on $B_R$ and $\chi=0$ on $B_{2R}^c$. Then, we have 
\begin{equation}\label{eq:conv_lemma_eq1}
\int \varphi\,\mathrm{d}\rho_n
  = 
\int \varphi\!\left(x_n + \f{x}{\lambda_n}\right)\,\chi(x)\,\mathrm{d}\tilde{\rho}_n
 + 
\int \varphi\!\left(x_n + \f{x}{\lambda_n}\right)\,\left(1-\chi(x)\right)\,\mathrm{d}\tilde{\rho}_n.
\end{equation}
We take into account the two possible cases
\begin{enumerate}
    \item $\{x_n\}_n$ is unbounded,
    \item $\{x_n\}_n$ is bounded. Therefore, up to a subsequence (not relabeled), we have $x_n\to x_{\infty}$.
\end{enumerate}
In the first case, since $\chi$ is compactly supported, $\lambda_n\to \infty$, and $\varphi$ is continuous, we deduce that
$$
\chi(x)\varphi\left(x_n + \f{x}{\lambda_n}\right) \to 0 \quad \text{uniformly on $\R^d$}.
$$
This yields that the first term on the right-hand side of~\eqref{eq:conv_lemma_eq1} vanishes as $n\to+\infty$. The second term of~\eqref{eq:conv_lemma_eq1} can be bounded by 
$$
\left|\int_{\R^d}\varphi\left(x_n + \f{x}{\lambda_n}\right) (1-\chi(x))\diff \tilde{\rho}_n\right| \le \|\varphi\|_{L^{\infty}} \sup_{n}\tilde{\rho}_n(B_{R}^c).
$$
Since $R$ is arbitrary and $\tilde{\rho}_n$ is tight, we can make this term as small as we want by taking the limit $R\to +\infty$.  Finally, we take the limits $n\to+\infty$ and $R\to +\infty$ in~\eqref{eq:conv_lemma_eq1} in order to obtain
$$
\int_{\R^d}\varphi\diff \rho = 0,
$$
for all $\varphi\in C_{b}(\R^d)$. This yields $\rho=0$, which is a contradiction as $\rho$ is supposed to be of mass $M>0$. 
Therefore, we can assume $x_n\to x_{\infty}$. Thereby, if we pass to the limit in~\eqref{eq:conv_lemma_eq1} we obtain
$$
\left|\int_{\R^d}\varphi\diff \rho - \varphi(x_\infty) \int_{\R^d}\chi \diff \tilde{\rho}\right|\le \|\varphi\|_{L^{\infty}}\sup_{n}\tilde{\rho}_{n}(B_{R}^c).
$$
Finally, we take $R\to+\infty$, in order to obtain  
$$
    \int_{\R^d}\varphi \diff \rho = M\varphi(x_\infty). 
$$
Since $\varphi$ is an arbitrary function of $C_{b}(\R^d)$, we conclude that $\rho = M \delta_{x_{\infty}}.$
\end{proof}

Now, we are ready to recover an $L^\infty$ bound on $\rho$ in terms of $\delta$.

\begin{proof}[Proof of Proposition~\ref{prop:linf_estimate_critical}]
The proof is based on a contradiction argument which applies to both cases, dimension $d=2$ and $d>2$. Assume the result is false. Then, we can construct a sequence $\rho_{n}\in K$ satisfying 
$$
\|\rho_{n}\|_{L^{\infty}}\ge n(|\delta_n| +1),
$$
which,  in particular, implies that $\|\rho_{n}\|_{L^{\infty}}\to\infty$.
We then define $\tilde{\rho}_{n}=T_{n}\sharp\rho_n$ the pushfoward of the measure $\rho_n$, where $T_n(x) = \lambda_n (x-x_n)$ with $\{x_n\}_n$ a sequence of $\R^d$  and $\{\lambda_n\}_n$ a sequence of real numbers such that $\lambda_n\to +\infty$. In particular, we have
$$
\tilde{\rho}_{n}(x) = \lambda_n^{-d} \rho_{n}(\lambda_n^{-1}(x-x_n)) .
$$
Choosing $\lambda_n^{d}= \|\rho_n\|_{L^{\infty}}\to\infty$ we obtain that
$$
\|\tilde{\rho}_n\|_{L^{\infty}}=1.
$$
We can also choose $x_n$ so that $0$ is a Lebesgue point of $\tilde{\rho}_n$ and $\tilde{\rho}_{n}(0)\ge \frac{1}{2}$.
Moreover, $\tilde{\rho}_{n}(\R^d) = \rho_n(\R^d)= M^{sub}_c(d)$, therefore, there exists a measure $\tilde{\rho}$ such that up to a subsequence (not relabeled), we have $
\tilde{\rho}_{n}\rightharpoonup\tilde{\rho}
$
where the weak convergence is understood against every continuous compactly supported function. Let us observe that some of the mass could escape at $\infty$ at this stage and that the convergence may not be narrow. In the following, we will show that this is not the case. First, let us notice that 
$\tilde{\rho}(\R^d)\le M^{sub}_c(d)$. Then, we use the hypothesis of the proposition, namely $\Delta v[\rho_n]\ge \delta_n$. Thus, we can repeat  the proof of Proposition~\ref{prop:linf_estimate} in order to conclude that $\tilde\rho$ is a nontrivial subsolution of 
\begin{align*}
&\Delta\log(\tilde{\rho}) + \tilde{\rho}\ge 0 \qquad \quad \, \, \text{in dimension $d=2$}, \\
&\f{m}{m-1}\Delta\tilde{\rho}^{m-1}+\tilde{\rho}\ge 0 \quad \, \text{in dimension $d>2$}.
\end{align*}
By definition this implies that $\tilde{\rho}(\R^d)\ge M^{sub}_c(d)$. Hence, we can conclude that $\tilde{\rho}(\R^d)=M^{sub}_c(d)$. From the weak convergence against $C_c(\R^d)$, functions and convergence of the mass, we can upgrade the convergence against every $C_{b}(\R^d)$ function (see \cite[Lemma 5.8]{Santambrogio15}). Using Lemma~\ref{lem:contradiction_ Dirac_mass}, we conclude that $\rho_{n}$ is tight and converges weakly against continuous bounded functions towards $M\delta_{x_{\infty}}$, that is, the convergence is narrow. This is in contradiction with the fact that $K$ does not contain Dirac masses. 
\end{proof}

\section{Estimating \texorpdfstring{$Q(u)$}{Q(u)}}\label{sec:d11-d22}

In this section, we focus on controlling the difference of the second derivatives of $u$, i.e. $Q(u)$. This is the next crucial step in order to prove the  Li--Yau and Aronson--B\'enilan estimates. We consider two cases. 
\begin{enumerate}
\item Small mass in dimension $2$: First, we provide an "easy" estimate in dimension $2$, via a rather intuitive argument that gives an explicit constant. We conjecture that at least a part of the estimate is optimal and cannot be improved. This leads to a bound of the form $Q(u)\leq C_3(M)|\delta|$, and it is only useful if the constant $C_3(M)$ is small enough. This is why we only apply this in the case where the mass $M$ itself is small.
\item Subcritical and critical mass and higher dimension: To treat all remaining cases ($d>2$ or large mass), we use a more involved argument. In particular, we want to underline the dependence of the bounds on $Q(u)$ in terms of other quantities than just $|\delta|$, such  $\|\rho\|_{L^2}$.
\end{enumerate}

The key point for these estimates is to use the identity $u=\Gamma*\rho$ where $\Gamma$ corresponds to the fundamental solution of the Laplacian in dimension $d$. In particular, $\Gamma$ is such that $-\Delta\Gamma=\delta_0$. For our proof, we will take advantage of other linear combinations of the second derivatives.

First of all, let us notice that we have that
\[
\partial_i \Gamma(x)
 = 
-\f{1}{d\omega_d}
\f{x_i}{|x|^d}.\]
This function belongs to $L^1_{loc}(\R^d)$ and it is a distribution. We now compute $\partial_{ii}\Gamma$ in the sense of distributions. It is such that
$$
    \langle \partial_{ii}\Gamma, \varphi \rangle \coloneqq -\langle \partial_{i}\Gamma,\partial_i\varphi \rangle = \f{1}{d\omega_d}\int_{\R^d}
    \f{x_i}{|x|^d}\partial_i\varphi(x)\diff x.
$$
We split the last integral into two parts:
$$
    \f{1}{d\omega_d}\int_{\R^d}
    \f{x_i}{|x|^d}\partial_i\varphi(x)\diff x=\f{1}{d\omega_d}\int_{B_r^c}
    \f{x_i}{|x|^d}\partial_i\varphi(x)\diff x+\f{1}{d\omega_d}\int_{B_r}
    \f{x_i}{|x|^d}\partial_i\varphi(x)\diff x.
$$
The second integral on the right hand side converges to  $0$ as $r\to 0$.
This follows since the integral is $L^1$ around the origin because $\partial_i\Gamma$ is integrable and $\partial_i \varphi$ is smooth and in particular bounded. 
Regarding the first integral on the right hand side we integrate by parts and we obtain  
$$
    \f{1}{d\omega_d} \int_{B_r^c} 
    \f{x_i}{|x|^d}\partial_i \varphi(x) \diff x = -\f{1}{d\omega_d} \int_{B_r^c} \left(\f{1}{|x|^d} - d\f{(x_i)^2 }{|x|^{d+2}}\right) \varphi(x) \diff x + \int_{\partial B_r}\f{x_in_i}{|x|^d} \varphi(x) \diff x.
$$
For the boundary integral, we use $\varphi(x)\to \varphi(0)$ as $r\to 0$ and $\int_{\partial B_r}\f{x_in_i}{|x|^d}\diff x=\omega_d$. Therefore, combining everything  we deduce 
$$
    \partial_{ii}\Gamma = \frac{1}{\sigma_d} \left(\f{1}{|x|^d}- d\f{(x_i)^2 }{|x|^{d+2}}\right) - \frac1d \delta_0.
$$
This distributional derivative has two parts. The first one is not an $L^1_{loc}$ function and should be considered in the principal value sense, i.e. when tested against a smooth function $\varphi$ the limit of the integral on $B_r^c$ exists and is finite. This is because $\f{1}{|x|^d}-
d\f{(x_i)^2 }{|x|^{d+2}}$ is of zero average on each sphere, so  we can replace $\varphi(x)$ with $
\varphi(x)-\varphi(0)$ and obtain integrability around the origin.
With regard to the second part, when computing a difference 
$\partial_{ii}\Gamma-\partial_{jj}\Gamma$, the Dirac mass disappears and we only have to estimate the integral part.

\subsection{Small mass case in dimension 2}\label{subsec:small_mass_case_dimension_2}
First, we provide an estimate which is useful in the small mass case in dimension 2 (case~\ref{2d_small_mass}). 

\begin{proposition}
\label{prop:estimate_Linf_difference}
Let $u$ be the solution of $-\Delta u=\rho$ given by $u=\Gamma*\rho$. Then
\[
Q(u)
 \le 
\f{1}{\,4\pi\,}\,\|\Delta \rho\|_{\,L^1}.
\]
 Moreover, if we take into account $\delta =  \inf_{x \in \mathbb{R}^d}\,\Delta v(x)$, then we also have
\[
\|\Delta \rho\|_{L^1}
 \le 
2M|\delta| + 2M\|\rho\|_{L^\infty}.
\]
Finally, using
\[
\|\rho\|_{L^\infty}
    \le C_2(M)|\delta|,
\]
we obtain
\[
Q(u)
\le
\frac{M}{2\pi}
\bigl(1+C_2(M)\bigr)|\delta|
\eqqcolon
C_3(M)|\delta|.
\]
For the explicit two-dimensional small-mass choice
$C_2(M)=M/(4\pi-M)$, this gives
\[
C_3(M)
=
\frac{M}{2\pi}
\bigl(1+C_2(M)\bigr)
=
\frac{2M}{4\pi-M}.
\]
\end{proposition}


\begin{remark}\label{rem:corollary}
A noteworthy side observation is the following: if $\rho \in L^\infty(\mathbb{R}^2)$ and $\Delta\rho \in L^{1}(\mathbb{R}^2)$, then the Hessian $D^2 u$ is automatically in $L^\infty(\mathbb{R}^2)$. Indeed, a bound on $\Delta u$ controls the sum of the eigenvalues of $D^2 u$, while a bound on $Q(u)$ controls their difference, giving a uniform bound on the full Hessian. Note that the easiest way to bound $D^2u$ in $L^\infty$ would require to asusme $\rho\in C^{0,\alpha}$ (so that $u\in C^{2,\alpha}$).
Nevertheless, we can observe that the assumption $\rho \in L^\infty$ and $\Delta\rho \in L^{1}$ does not imply any modulus of continuity for $\rho$. For a counterexample, let us consider a regularisation of the family of functions $x\mapsto \min\{\max\{-\log|x|,k\},k+1\}-k$. This family of functions is bounded by $1$ and the mass of the Laplacian is $2$. However, they take the value $1$ at the origin and $0$ outside of the ball $B(0,e^{-k})$.
\end{remark}

\begin{proof}[Proof of \Cref{prop:estimate_Linf_difference}]
Without loss of generality, we consider the point $x=0$ and the vectors $e=(1,0)$ and $e'=(0,1)$  (a change of coordinates reduces the general case to this one), 
denoting $\partial_1$ and $\partial_2$ the two partial derivatives in the direction of the vectors of the canonical basis, and using the computation for the distributional second derivatives of $\Gamma$, we have  
\[
(\partial_{11}u-\partial_{22}u)(0)
 = 
-\f{1}{\pi}
\int_{\mathbb{R}^2}\rho(z)\,\f{z_2^2 - z_1^2}{|z|^4}\,\mathrm{d}z,
\]
where the integral has to be taken in the principal value sense.

We observe that in polar coordinates $z=(r\cos\theta,\,r\sin\theta)$ the function  $\f{z_2^2 - z_1^2}{|z|^4}$ equals $\f{\cos(2\theta) }{r^2}$. 
%
We will show that this new expression can be rewritten as the Laplacian of a bounded function. We recall that in dimension 2 for $f(r,\theta)$ the Laplacian in radial coordinates is
$$
\Delta f= \f{1}{r}\partial_{r}(r\partial_{r}f) + \f{1}
{r^2}\partial_{\theta\theta}f.$$
Then, we point out that the function $f(r,\theta)= \f{1}{4}\cos(2\theta)$ satisfies $\Delta f = -\f{\cos(2\theta)}{r^2}$. With similar computations to the ones for the second derivative of $\Gamma$, we can show that this is the true distributional derivative and no other term at the origin is present. This function can be written in the Euclidian basis $f(y_1,y_2) = -\f{1}{4}\f{y_1^2-y_2^2}{y_1^2+y_2^2}$ and it is such that $\|f\|_{L^{\infty}}\le 1/4.$
Thereby, if we come back to Euclidean coordinates, we obtain 
$$
    (\partial_{11}u-\partial_{22}u)(0,0) = -\f{1}{\pi}\int_{\R^2}\Delta f\diff \rho.
$$
Performing two integrations by parts, we can estimate 
$$
|(\partial_{11}u-\partial_{22}u)(0,0)|\le \f{1}{4\pi}\int_{\R^2}|\Delta\rho|.
$$
This proves the first claim. 

Let us now recover the $L^1$-bound on $\Delta\rho$.
We recall $\Delta v = \Delta[\log\rho-u]= \Delta\log\rho + \rho$ since $-\Delta u =\rho$. Therefore, we have
$$
    \Delta\rho = \rho\Delta v -\rho^2 +\f{|\nabla\rho|^2}{\rho}\ge \rho\Delta v - \rho^2.
$$
This means $(\Delta\rho)_{-}\le \rho|\delta|+\rho^2$ where $s_{-}$ denotes the negative part of a number $s$. Since $\int \Delta\rho=0$ then $\int (\Delta\rho)_{-} = \int (\Delta\rho)_{+} $ and
$
\int_{\R^2} |\Delta\rho |= 2\int_{\R^2} (\Delta\rho)_{-}$. Therefore 
$$
    \int_{\R^2} |\Delta\rho| \le 2M |\delta| + 2\int_{\R^2}\rho^2.
$$
Then, we conclude the proof by observing
\begin{equation*}
    \int_{\R^2}\rho^2\le \|\rho\|_{L^{\infty}}\int_{\R^2}\rho. \qedhere
\end{equation*}
\end{proof}

\subsection{General estimate}

The previous proof from Subsection~\ref{subsec:small_mass_case_dimension_2} does not extend to the higher-dimensional small mass case~\ref{2d_small_mass_multid}. 
Thus, in this subsection, we develop a more general theory that covers all the remaining cases.
\begin{proposition}
\label{prop:estimate_Linf_difference_full_mass_multid}
Let $u$ be such that $-\Delta u=\rho$. Let us also consider $\delta = \min(\Delta v)$ and $\bar{\delta}=\max(|\delta|,\|\rho\|_{L^{\infty}})$. Then, we have
$$
Q(u) \le A_d\bar\delta^{d/(d+2)}\|\rho\|_{L^2}^{4/(d+2)}.
$$
Furthermore, the constant $A_d$ in the previous inequality can be explicitly computed, and it only depends on $d$.
\end{proposition}

\begin{remark}
At first glance, the $\|\rho\|_{L^2}$ might seem unnatural. However, it includes a result similar to the one of the previous section if one just notices that $\|\rho\|_{L^2}\leq \|\rho\|_{L^\infty}^{1/2}\|\rho\|_{L^1}^{1/2}\leq M^{1/2}\bar\delta^{1/2}$. In particular, this observation together with \Cref{prop:estimate_Linf_difference_full_mass_multid} implies $Q(u)\leq A_d M^{2/(d+2)}\bar\delta$. Whenever the estimate $\|\rho\|_{L^\infty}\leq C_2(M)|\delta|$ of the previous section is available, we have $\bar\delta\leq(1+C_2(M))|\delta|$, so that for $d>2$ one may take
$$
C_3(M) = A_d M^{2/(d+2)}C_2(M)^{2/(d+2)}(1+C_2(M))^{d/(d+2)},
$$
which is less explicit than in dimension $2$. Nevertheless, we state the result using the $L^2$ norm because it will be relevant later in the paper. Combined with an additional control on the free energy, we will use the above estimate in order to obtain that $Q(u)= o(|\delta|)$.
\end{remark}
Before proving this result, we state two useful inequalities:
\begin{lemma}\label{lemma estimate laplacian}
 We have  $\Delta\rho\geq \f{2}{d}\f{|\nabla\rho|^2}{\rho}-(\rho-\delta)\f{\rho^{2/d}}{m}$. In particular, 
$\Delta\rho\geq -\f{2}{m}\bar\delta \rho^{2/d}$.
\end{lemma}
\begin{proof}
We start from $\Delta v\ge \delta$ and expand $\Delta\log\rho$ (if $d=2$) or $\Delta\rho^{1-2/d}$ (if $d>2$).
This gives
$$m\rho^{-2/d}\Delta\rho-m\frac{2}{d}\rho^{-1-2/d}|\nabla\rho|^2+\rho\geq\delta.$$
We then multiply times $\rho^{2/d}/m$ in order to obtain the first estimate. Then, we use $|\delta|,\rho\leq\bar\delta$ and the positivity of the term $|\nabla\rho|^2/\rho$.
\end{proof}
\begin{lemma}\label{lemma local fischer}
 We have  $
\int_{B_{2R}}\f{|\nabla \rho|^2}{\rho}\lesssim \|\rho\|_{L^{2}} R^{d/2-2} + \|\rho\|_{L^{2}}^{2/d} R^{d-1}\bar\delta.
$
\end{lemma}
\begin{proof}
    
We start from the estimate on $\f{|\nabla\rho|^2}{\rho}$ contained in the previous lemma. Then, we choose $\chi:\R^{d}\to [0,1]$ a smooth function such that
\[
\chi(x)  = 
\begin{cases}
1, & \text{if } |x|\le 2R,\\
0, & \text{if }|x|\ge 3R,
\end{cases}
\quad
|\nabla\chi| \lesssim \f{1}{\,R\,},
\quad
|\Delta\chi| \lesssim \f{1}{\,R^2\,}.
\]
Hence, from the properties of $\chi$ and the Young's inequality it follows that
\begin{align*}
\int_{B_{2R}}\f{|\nabla \rho|^2}{\rho}&\lesssim \int_{\R^d}\f{|\nabla\rho|^2}{\rho}\, \chi\, \diff x \lesssim \int_{\R^d} \rho\, \Delta\chi\diff x + \bar\delta \int_{\R^d} \rho^{2/d}\, \chi\diff x\\
& \lesssim \f{1}{R^2} \|\rho\|_{L^2} R^{d/2} + \bar\delta\|\rho\|_{L^2}^{d/2} R^{d-1}. \qedhere
\end{align*}
\end{proof}

We are now ready to show \Cref{prop:estimate_Linf_difference_full_mass_multid}. We consider the cases of dimension 2 and dimension $d>2$ individually. 

\begin{proof}[Proof of Proposition~\ref{prop:estimate_Linf_difference_full_mass_multid} in dimension 2] 
We divide the proof into several steps.
\begin{enumeratesteps}
\step[A new expression for $Q(u)$]
We recall
\[
(\partial_{11}u-\partial_{22}u)(0)
 = 
-\f{1}{\pi}
\int_{\mathbb{R}^2}\rho(z)\,\f{z_2^2 - z_1^2}{|z|^4}\,\mathrm{d}z.
\]
%
%
%
We define $R>0$ and we take $\chi:\mathbb{R}^2\to[0,1]$  a smooth, monotone decreasing function such that
\[
\chi(x)  = 
\begin{cases}
1, & \text{if }|x|\le R,\\
0, & \text{if }|x|\ge 2R,
\end{cases}
\quad
|\nabla\chi| \lesssim \f{1}{\,R\,},
\quad
|\Delta\chi| \lesssim \f{1}{\,R^2\,}.
\]
We split the integral into two terms:
\[
(\partial_{11}u - \partial_{22}u)(0)
 = 
 - 
\f{1}{\pi}
\int_{\R^2}
\rho\,\f{\cos(2\theta)}{r^2}\,\chi  \diff x
-\f{1}{\pi}
\int_{\R^2}
\rho\,\f{\cos(2\theta)}{r^2}\,(1-\chi)   \diff x = A+B.
\]
In the following two steps we bound both integral terms, $A$ and $B$, respectively.

\step[An estimate for $A$]
We have already proven that the function
\[
F(r,\theta) 
= \f{\cos(2\theta)}{r^2}
\]
is the Laplacian in radial coordinates of 
$\displaystyle g(\theta) = -\f{1}{4}\cos(2\theta).$
We perform an integration by parts on $A$ to obtain that 
$$
    A = -\f{1}{\pi}\int_{\R^2} g(\theta)\, \Delta (\rho\chi)\, \diff x .
$$
Since $|g|$ is bounded uniformly, we obtain 
\[
    |A| \lesssim \|\Delta (\rho\chi)\|_{L^1(\R^2)}=2\|(\Delta (\rho\chi))_-\|_{L^1(\R^2)}.
\]
Lemma \ref{lemma estimate laplacian} in the case $d=2$ provides $
(\Delta\rho)_{-} \lesssim \bar\delta\rho
$. Thus, by properties of $\chi$ we deduce
$$
    \left[\Delta(\rho\chi)\right]_{-} = \left[(\Delta \rho)\chi + 2\nabla\rho\cdot\nabla \chi + \rho\Delta\chi \right]_{-}\lesssim \bar\delta\rho \mathbbm{1}_{B_{2R}} + \left(\f{|\nabla \rho|}{R} + \f{\rho}{R^2}\right) \mathbbm{1}_{B_{2R}\setminus B_{R}}.
$$
With this inequality, we obtain 
$$
    |A|\lesssim \bar\delta\int_{B_{2R}} \rho + \f{1}{R} \int_{B_{2R}}|\nabla\rho|+ \f{1}{R^{2}} \int_{B_{2R}}\rho.
$$
The first and the last term on the right-hand side above can be estimated using H\"older's inequality and we obtain $\int_{B_{2R}} \rho\lesssim \|\rho\|_{L^2} R$. The estimate of the second one requires the use of Lemma \ref{lemma local fischer} is more delicate.
%
%
%
%
We compute and we get that 
\begin{align*}
\f{1}{R} \int_{B_{2R}}|\nabla\rho| &\lesssim \f{1}{R}\left(\int_{B_{2R}}\rho\right)^{1/2}\left(\int_{B_{2R}}\f{|\nabla \rho|^2}{\rho}\right)^{1/2}\\
&\lesssim \f{1}{R}\left(\|\rho\|_{L^2} R\right)^{1/2}\left(\|\rho\|_{L^{2}} R^{-1} + \|\rho\|_{L^{2}} R\bar\delta\right)^{1/2}\lesssim \|\rho\|_{L^2}R^{-1} + \sqrt{\bar\delta}\|\rho\|_{L^2}.
\end{align*}
Therefore, if we combine everything and we use Young's inequality, we deduce
\begin{align*}
|A| \lesssim \|\rho\|_{L^2} (\bar\delta R+R^{-1}+ \sqrt{\bar\delta}) \lesssim \|\rho\|_{L^2} (\bar\delta R+R^{-1}).
\end{align*}

\step[An estimate for $B$]
Since $|\cos(2\theta)(1-\chi)| \le 1$ and $1-\chi$ vanishes on $B_R$, we obtain 
\[
    |B| \lesssim \int_{\R^2\setminus B_R} \f{\rho}{r^2} 
        \lesssim \|\tilde{\rho}\|_{L^2}\,\sqrt{ \int_{\R^2\setminus B_R} \f{1}{r^4} }
        \lesssim \f{\|\rho\|_{L^2}}{\,R\,}.
\]
Thereby, combining the results from the three steps, we get 
$$
|(\partial_{11}u - \partial_{22}u)(0)|\lesssim \|\rho\|_{L^2} (\bar\delta R+R^{-1}).
$$
Finally, we optimise with respect to $R$ and we choose  $R= \bar\delta^{-1/2}$. 
This yields the result.   \qedhere 
\end{enumeratesteps}
\end{proof}

We now prove the statement in the case of higher dimensions.

\begin{proof}[Proof of Proposition~\ref{prop:estimate_Linf_difference_full_mass_multid} in dimension $d>2$]
We work analogously to the previous case of dimension $2$.
\begin{enumeratesteps}
\step[A new expression for $Q(u)$]    
As before, without loss of generality, we consider $(\partial_{11}{u} - \partial_{22}u)(0)$ and
we compute
\[
    (\partial_{11} u - \partial_{22} u)(0) = \f{d}{\sigma_d} \int_{\mathbb{R}^d} \rho(y)\,\f{y_1^2 - y_2^2}{|y|^{d+2}} \,\diff y.
\]
We define $R>0$ and we take $\chi:\mathbb{R}^d\to[0,1]$  a smooth, monotone decreasing function such that
\[
\chi(x)  = 
\begin{cases}
1, & \text{if }|x|\le R,\\
0, & \text{if }|x|\ge 2R,
\end{cases}
\quad
|\nabla\chi| \lesssim \f{1}{\,R\,},
\quad
|\Delta\chi| \lesssim \f{1}{\,R^2\,}.
\]
We split the integral into two terms:
\[
(\partial_{11}u - \partial_{22}u)(0)
=
\int_{\R^{d}}
\rho\,\chi\,
\f{y_1^2 - y_2^2}{|y|^{d+2}}\,
\diff x+\int_{\R^{d}}
\rho\,(1-\chi)\,
\f{y_1^2 - y_2^2}{|y|^{d+2}}\,
\diff x = A + B.
\]
We devote the next two steps to bounding both of them, $A$ and $B$ respectively. 

\step[An estimate for $A$]
We claim that the function
\[
F(y) 
= \f{y_1^2 - y_2^2}{|y|^{d+2}}
\]
is the Laplacian of a function which writes, in polar coordinates $\frac{g(w)}{r^{d-2}}$. More precisely, if we consider a function of the form $f(y)=|y|^{-d} a(y)$, we obtain, after some computations,
$$\Delta f(y)=\frac{\Delta a(y)}{|y|^d}+\frac{2d}{|y|^{d+2}}(a(y)-\nabla a(y)\cdot y).$$
If we choose $a(y)=y_1^2-y_2^2$ we have $\Delta a(y)= 0$ and $\nabla a(y)\cdot y=2a(y)$, so that we obtain $\Delta f(y)=-2d F(y)$. This shows that it is enough to choose $g(\omega)=-\frac{1}{2d}(\omega_1^2-\omega_2^2)$.


Afterwards, we integrate by parts the term $A$ and we get  
$$
    A = \int_{\R^{d}} \f{g(\omega)}{r^{d-2}}\, \Delta (\rho\chi)\, \diff x.
$$
Since $|g|$ is bounded uniformly  we can use that  $|s| = 2 s_{-} + s$ in order to obtain that  
\begin{align*}
|A| &\lesssim \int_{\R^{d}}
\f{1}{r^{d-2}}\, |\Delta (\rho\chi)|\,
\diff x\\
&= \int_{\R^{d}}
\f{2}{r^{d-2}}\, [\Delta (\rho\chi)]_-\,
\diff x+\int_{\R^{d}}
\Delta\left(\f{1}{r^{d-2}}\right)\rho\chi\,
\diff x\\
&\lesssim \int_{\R^{d}}
\f{1}{r^{d-2}}\, [\Delta (\rho\chi)]_-\,
\diff x.
\end{align*}
Here, the last inequality follows from
$
\Delta\left(\f{1}{r^{d-2}}\right) = -C\delta_{0}\leq 0.
$
Therefore the second integral is nonpositive. From \Cref{lemma estimate laplacian} we know  $
(\Delta\rho)_{-} \lesssim \bar\delta\rho^{2/d}
$. From the properties of $\chi$ we deduce
\begin{align*}
\left[\Delta(\rho\chi)\right]_{-} &= \left[(\Delta \rho)\chi + 2\nabla\rho\cdot\nabla \chi + \rho\Delta\chi \right]_{-} \lesssim \bar\delta\rho^{2/d}\mathbbm{1}_{B_{2R}} + \left(\f{|\nabla \rho|}{R} + \f{\rho}{R^2}\right)\mathbbm{1}_{B_{2R}\setminus B_{R}}.
\end{align*}
Thanks to this inequality, we obtain 
\begin{align*}
|A|&\lesssim \bar\delta\int_{B_{2R}} \f{\rho^{2/d}}{r^{d-2}} + \f{1}{R^{d-1}} \int_{B_{2R}} |\nabla\rho|\,  + \f{1}{R^{d}} \int_{B_{2R}} \rho .
\end{align*}
Let us now study each one of the three integral terms on the right-hand side. 
Concerning the first one, H\"older's inequality allows us to bound it by
$$
\bar\delta\int_{B_{2R}} \f{\rho^{2/d}}{r^{d-2}} \leq \bar\delta\left(\int_{B_{2R}} \rho^{2}\right)^{1/d}\left(\int_{B_{2R}}r^{-\frac{d}{d-1}(d-2)} \right)^{(d-1)/d} =C \bar\delta\|\rho\|_{L^2}^{2/d} R.
$$
The third one can  be estimated by $ \|\rho\|_{L^2}R^{-d/2}$. It also follows from H\"older's inequality.

For the second term, we use  \Cref{lemma local fischer} and we obtain 
\begin{align*}
\f{1}{R^{d-1}} \int_{B_{2R}} |\nabla\rho| &\lesssim \f{1}{R^{d-1}}\left(\int_{B_{2R}}\rho\right)^{1/2}\left(\int_{B_{2R}}\f{|\nabla \rho|^2}{\rho}\right)^{1/2}\\
&\lesssim \f{1}{R^{d-1}}\left(\|\rho\|_{L^2} R^{d/2}\right)^{1/2}\left(\|\rho\|_{L^{2}} R^{d/2-2} + \|\rho\|_{L^{2}}^{2/d} R^{d-1}\bar\delta\right)^{1/2}\\
&\lesssim \|\rho\|_{L^2}R^{-d/2} + \sqrt{\bar\delta}\|\rho\|_{L^2}^{1/2+1/d}R^{-d/4 +1/2}.
\end{align*}
Finally, if we combine all the estimates and we use Young's inequality, we deduce 
\begin{align*}
|A| \lesssim \bar\delta R\|\rho\|_{L^2}^{2/d} + \|\rho\|_{L^2} R^{-d/2} + \sqrt{\bar\delta}\|\rho\|_{L^2}^{1/2+1/d}R^{-d/4 +1/2} \lesssim \bar\delta R\|\rho\|_{L^2}^{2/d} + \|\rho\|_{L^2} R^{-d/2}.
\end{align*}

\step[An estimate for $B$]
In order to bound $B$, we follow a similar strategy. Combining H\"older's inequality, the definition of $\chi$ and the bound $|\omega_{1}^{2}-\omega_{2}^2|\leq 1$ we deduce that
$$
    |B|\lesssim \int_{B_R^c} \frac{\rho(r\omega)}{r^d} \lesssim \|\rho\|_{L^{2}} \left( \int_{B_R^c} \f{1}{r^{2d}} \right)^{1/2} = C \|\rho\|_{L^{2}} R^{-d/2}.
$$
Thus, from the combination of the three steps, we obtain
$$
    |(\partial_{11}u - \partial_{22}u)(0)|\lesssim \bar\delta R\|\rho\|_{L^2}^{2/d} + \|\rho\|_{L^2} R^{-d/2}.
$$
Analogously to the $2$-dimensional case, we can optimize with respect to $R$. Therefore, we can choose $R= \|\rho\|_{L^2}^{2(d-2)/d(d+2)}\bar\delta^{-2/(d+2)}$ in order to yield the result. \qedhere
\end{enumeratesteps}
\end{proof}

\section{Subsolutions of Liouville and Lane--Emden equations}\label{sec:liouville}

In this section, we study the subsolutions of the Lioville and the Lane--Emden equations. In particular, we obtain the minimal mass of their nonnegative subsolutions. As explained before in \Cref{sec:Linfty}, obtaining these values is key for our analysis of the Keller--Segel equation. 
Let us start with a short summary of our analysis for both problems:

\begin{itemize}
\item The Liouville equation in dimension $d=2$
\[
    \Delta h  +  e^{h}  \ge  0.
\]
We show that any nontrivial subsolution satisfies $\int_{\R^2} e^h\ge 8\pi$. This result is well known for solutions of the equation \cite[Lemma 1.1]{Chen_Li91}, and the adaptation of this result to subsolutions is straightforward. Therefore $M_c^{sub}(2)=M_c(2)=8\pi$.
\item The Lane--Emden equation in dimension $d>2$
\[
\f{m}{m-1}\Delta h  +  h^{\f{d}{\,d-2\,}}  \ge  0.
\]
For any nonnegative, nontrivial subsolution we define its ``mass'' $\int h^{\f{d}{\,d-2\,}}$ and call $M_c^{sub}(d)$ the minimal possible mass of nonnegative, nontrivial subsolution, as it is given by \eqref{def:Mc LaneEmden}. We prove that this critical mass coincides with the critical mass from Definition~\ref{def:critical_mass}, i.e.  $M_c^{sub}(d)=M_{c}(d)$ for every $d\geq 3$.
This result appears to be new and relies on an optimal-control argument for a radially decreasing rearrangement of the subsolution.  
\end{itemize}

\begin{remark}
It is straightforward to see that the subsolutions of the Lane-Emden equations must have a mass strictly greater than some $\varepsilon>0$. Indeed, we multiply the equation
\[
\frac{m}{m-1}\,\Delta h + h^q \ge 0
\]
(for $q=\f{d}{d-2}$)
by $h$ and we integrate over $\R^d$ (assuming decay properties on $h$). Hence, we obtain 
$$
\|\nabla h\|_{L^{2}}^2\le c\|h\|_{L^{q+1}}^{q+1}.
$$
However the Gagliardo-Nirenberg inequality states, for this choice of the exponents (remember $q+1=2\f{d-1}{d-2}$) for some $C_{GN}>0$
$$
\|h\|_{L^{q+1}}\le C_{GN} \|\nabla h\|_{L^{2}}^{\f{d-2}{d-1}}\|h\|_{L^q}^{\f{1}{d-1}}.
$$
Combining these two inequalities we obtain
$$\|h\|_{L^{q+1}}\le C\|h\|_{L^{q+1}}\|h\|_{L^q}^{\f{1}{d-1}},
$$
which proves a strictly positive lower bound on $\|h\|_{L^q}$.
Our goal is to prove that this bound can be extended up to  $M_c$.   
\end{remark}

In a more precise way, in this section we prove the following two main results for the Liouville and the Lane--Emden equations. 

\begin{lemma}[Subsolutions of the Liouville equation in $\mathbb{R}^2$]
\label{lem:subsolution_liouville}
Suppose $h\colon \mathbb{R}^2 \to [-\infty,\infty)$ is a locally bounded function belonging to $H^{1}_{loc}(\R^2)$ that satisfies:
\[
\begin{cases}
\displaystyle
\Delta h + e^h  \ge  0
&\text{in the set }\{\,x\in \R^2 : h(x) > t\,\}\text{ for all }t\in\R,\\[5pt]
\displaystyle
\int_{\mathbb{R}^2} e^h < +\infty.
\end{cases}
\]
Then $\int_{\mathbb{R}^2} e^h \,\mathrm{d}x \ge 8\pi.$
\end{lemma}
 We now state the analogous result for the Lane--Emden equation.

\begin{theorem}[Subsolutions of the Lame-Emden equation in $\mathbb{R}^d$]
\label{lem:subsolution_lane-emden}
Suppose $h\colon \mathbb{R}^d \to [0,+\infty)$ is a function in $L^{\infty}(\R^d)\cap L^{\f{d}{d-2}}(\R^d)\cap H^{1}_{loc}(\R^d)$, $h\not\equiv 0$ and that:
\[
\begin{cases}
\displaystyle
\f{m}{m-1}\Delta h + h^{\f{d}{d-2}}  \ge  0
&\text{ on $\R^d$ }\\[5pt]
\displaystyle
\int_{\mathbb{R}^2} h^{\f{d}{d-2}} < +\infty,
\end{cases}
\]
with $m=2-\f{2}{d}$. Then $\int_{\mathbb{R}^2} h^{\f{d}{d-2}} \,\mathrm{d}x \ge M_c(d)$, where $M_c (d)$ is given by \eqref{def:Mc}.
\end{theorem}

\subsection{Subsolutions of the Liouville equation in dimension 2}\label{sec:Liouville d 2}

We begin with the proof of the minimal mass condition for subsolutions to the Liouville equation. The statement was given previously in Lemma~\ref{lem:subsolution_liouville}. 

\begin{proof}[Proof of \Cref{lem:subsolution_liouville}]
The proof is just an adaptation of the proof of~\cite[Lemma 1.1]{Chen_Li91} in the case of subsolutions on level sets, that we rewrite here for the sake of completeness.
Let $\Omega_t = \{x\in \mathbb{R}^2:\,h(x)>t\}$.  If $h$ is a subsolution in the sense above, then on each $\Omega_t$ we have $\Delta h + e^h  \ge  0.$
On the one hand, since $\nabla h$ is orthogonal to $\partial\Omega_t$, from the divergence theorem it follows that
\[
    \int_{\Omega_t} e^h \geq -\int_{\Omega_t} \Delta h = -\int_{\partial\Omega_t} \f{\partial h}{\partial \nu} = \int_{\partial\Omega_t} \left|\nabla h\right|.
\]
On the other hand, by the coarea formula, we have
\[
-\f{\mathrm{d}}{\mathrm{d}t}\,\left|\Omega_t\right|
 = 
\int_{\partial\Omega_t} \f{1}{\,|\nabla h|\,}.
\]
Hence, by the Cauchy--Schwarz and the isoperimetric inequality in 2D, we recover that
\[
4\pi\,|\Omega_t|
 \le 
\left|\partial\Omega_t\right|^2
 \le 
\int_{\partial\Omega_t}\!|\nabla h| \int_{\partial\Omega_t} \!\f{1}{\,|\nabla h|\,}
 = 
-\f{\mathrm{d}}{\mathrm{d}t}\,\left|\Omega_t\right| \int_{\Omega_t} e^h.
\]
Thus
\[
\f{\mathrm{d}}{\mathrm{d}t}\,\left(\int_{\Omega_t} e^h\right)^2
 = 
2\,e^t\,\f{\mathrm{d}}{\mathrm{d}t}\,\left|\Omega_t\right| \int_{\Omega_t} e^h
 \le 
-8\pi\,e^t\,|\Omega_t|.
\]
Next, we integrate in $t$ from $-\infty$ to $+\infty$. Since $\Omega_t$ converges to $\emptyset$ as $t\to+\infty$ and to $\mathbb{R}^2$ as $t\to -\infty$, it follows that we have
\[
-\left(\int_{\mathbb{R}^2}e^h\right)^2
 \le 
-\,8\pi\int_{\mathbb{R}^2} e^h,
\]
and, in particular, $\int_{\mathbb{R}^2} e^h \ge 8\pi$.  This completes the argument.
\end{proof}

\subsection{Subsolutions of the Lane--Emden equation in dimension \texorpdfstring{$d>2$}{d>2}}\label{sec:Lane-Emden}

We now turn to the proof of \Cref{lem:subsolution_lane-emden}. The very same arguments of the Liouville equation do not apply in this case. The problem comes from the exponent of the isoperimetric inequality which is not 2 in higher dimensions. 
Recall that $m=2-2/d$, that $M_c$ is defined in \eqref{def:Mc},
and that $M_c^{sub}$ is defined in \eqref{def:Mc LaneEmden}.
Our goal is to prove $M_c^{sub}=M_c$. The known reference profile
$\bar h$ has mass $\int_{\R^d}\bar h^q=M_c$ and satisfies the
following properties.
\begin{itemize}
\item $\bar h$ is radially decreasing, compactly supported and continuous on $\R^d$;
\item on its support $B_R$ the function $\bar h$ satisfies the equality $\f{m}{m-1}\Delta h + h^{\frac{d}{\,d-2\,}}  =  0$;
\item $\bar h$ is unique up to scaling, i.e. replacing $\bar h$ with $x\mapsto \lambda^{d-2}\bar h (\lambda x)$.
\end{itemize}
Since the coefficient $m/(m-1)$ does not play a relevant role in the proof, for the sake of simplicity, we will remove it from the equation during this section. 
From a normalisation argument, the above statement is equivalent to showing that any subsolution of 
$$
    \Delta h + h^{\frac{d}{\,d-2\,}}  \ge  0
$$
has at least the same mass as the radial subsolution $\bar{h}$ that it is defined such that  $\Delta \bar{h} + \bar{h}^{\frac{d}{\,d-2\,}}  = 0$ on its support.

\paragraph{First intuition of the result.} Assume  we have
\begin{equation}\label{intuition alpha}
    \Delta h + h^{\frac{d}{\,d-2\,}}  = \alpha  \ge  0
\end{equation}
for some continuous function $\alpha\ge 0$. If there is a point $x_0$ where $\alpha(x_0) > 0$, then by continuity we have $\alpha\geq c_0>0$ on an open neighborhood $U$ around $x_0$ with $\alpha(x)>0$ for all $x \in U$. Taking $\varphi\ge 0$ smooth compactly supported in $U$ and setting $\psi \coloneqq h - \varepsilon \varphi$ for $\varepsilon>0$ small enough, we have decreased the $L^\f{d}{d-2}(\R^d)$ norm but
\[
    \Delta \psi + \psi^{\frac{d}{\,d-2\,}}   >  0.
\]
 This argument shows that $\Delta h + h^{\f{d}{\,d-2\,}}$ cannot be strictly greater than 0 on an open set. However, and this is where the difficulty lies, the Laplacian $\Delta h$ could be singular and the distribution  $\Delta h + h^{\f{d}{\,d-2\,}}$ could be a nonnegative measure whose support has empty interior. For instance, this is what happens for the function $\bar h$, extended to $0$ outside the ball $B_R$, and we expect this function to be the minimiser.

\paragraph{Minimisation in a larger class and reduction to an optimal control problem.}

We already explained that the two-dimensional strategy used for the Liouville equation cannot be reproduced here. Nevertheless, we are inspired by the strategy we follow in dimension 2 in our previous approach in subsection \ref{sec:Liouville d 2}. Let us take a more careful look to the strategy that we developed. We observe two facts.
\begin{itemize}
\item First, we only exploited $\int_{\{h>t\}}\Delta h+e^h\geq 0$ for every level set $t$, and not the validity of the inequality $\Delta h+e^h\geq 0$ everywhere. 
\item Second, we exploited some inequalities that are equalities when the functions are radial. The isoperimetric inequality is optimal on level sets being balls, and the Jensen inequality on the norm of the gradient on each level set is optimal when this norm is constant.
\end{itemize}
A priori, we should consider the minimization problem
\begin{equation*}\tag{\normalfont{P}}
    \min\left\{\,M(h)  =  \int_{\mathbb{R}^d} h^{\f{d}{\,d-2\,}},\quad h \in \mathcal{S}(\R^d)\right\}, 
\end{equation*}
where $\mathcal{S}(\R^d)$ is the set of nontrivial nonnegative subsolutions given by
\[
    \mathcal{S}(\R^d)  = \left\{\,h \in L^\infty(\R^d)\cap L^{\f{d}{\,d-2\,}}(\R^d)\cap H^1_{\mathrm{loc}}(\R^d), \,  h\ge0, \, \Delta h + h^{\f{d}{\,d-2\,}}\ge 0 \right\} \setminus \{0\}.
\]
We can prove that this problem has a minimizer. Yet, finding this solution does not appear to be simple, so, instead, we consider an extended problem:
\begin{equation}\tag{\normalfont{P}}
    \min\left\{\,M(h)  =  \int_{\mathbb{R}^d} h^{\f{d}{\,d-2\,}},\quad h \in \mathcal{S'}(\R^d)\right\}, \label{eq:Min problem}
\end{equation}
where $\mathcal{S'}(\R^d)$ is now given by
\[
\mathcal{S'}(\R^d)=
\left\{
\begin{aligned}
h\in{}& L^\infty(\R^d)\cap L^{\f{d}{d-2}}(\R^d)
\cap H^1_{\mathrm{loc}}(\R^d),\\
&h\ge0,\qquad
\int_{\{h>t\}}\!\!\bigl(\Delta h+h^{\f{d}{d-2}}\bigr)\ge0
\quad\text{for every }t
\end{aligned}
\right\}\setminus\{0\}.
\]
Let us point out that $\mathcal{S'}(\R^d)$ can  also be expressed as the set of functions satisfying
$$
    - \int \nabla\phi\cdot\nabla h+\int \phi h^{\f{d}{\,d-2\,}}\geq 0
$$
for every test function $\phi$ of the form $\phi=g(h)$ where $g$ is nondecreasing. This means 
$$
    \int g(h)\left( \Delta h + h^{\f{d}{\,d-2\,}}\right)\ge 0
$$ 
and the limit case $g=\ind_{[t,\infty)}$ provides the characterization on level sets. Let us also notice that we have that $\mathcal{S}(\R^d)\subset \mathcal{S'}(\R^d)$. 

Moreover, the set $\mathcal{S'}(\R^d)$ is invariant by symmetric-decreasing rearrangement. This comes from standard properties of the rearrangement, whose proofs can be found, for instance, in \cite{Elliot_Loss01}.
Indeed, let us take $h\in \mathcal{S'}(\R^d)$ and let $h^*$ be its rearrangement. We would like to prove 
$$ 
    \int \nabla(g(h^*))\cdot\nabla h^*\leq\int g(h^*) (h^*)^q 
$$
for $q = \frac{d}{d-2}$.
From well-known properties of the rearrangement, the right-hand side is equal to   $\int g(h) h^q$. With regard to   the left-hand side, it corresponds to the $L^2$ norm of the gradient of $\tilde g'(h^*)$ where $\tilde g'=\sqrt{g'}$, and  $\tilde g(h^*)$ is the rearrangement of  $\tilde g(h)$. Hence, from the Polya-Szego inequality, we have
\begin{align*}
    \int \nabla(g(h^*))\cdot\nabla h^* & = \int |\nabla \tilde g(h^*)|^2\leq  \int |\nabla \tilde g(h)|^2\\
    &= \int \nabla(g(h))\cdot\nabla h\leq \int g(h) h^q=\int g(h^*) (h^*)^q,
\end{align*}
which shows $h^*\in \mathcal{S'}(\R^d)$.
This allows us to simplify our optimization problem. Instead of considering a minimisation over the class $\mathcal S'(\Rd)$ we move to the one-dimensional calculus of variations problem: 
$$
\min\left\{\int_0^{R_0} r^{d-1}f^q(r)\,:\, f(0)=1,\, f \text{ is nonincreasing and }  x\mapsto f(|x|) \text{ belongs to } \mathcal{S'}(\R^d)\right\}.
$$
We point out that, ideally, we would like to consider $R_0=+\infty$. We will see later why it is more convenient to choose a very large value of $R_0$ and consider it to be finite.

It is interesting to re-write the condition  ``$x\mapsto f(|x|)$ belongs to $\mathcal{S'}(\R^d)$'' in terms of $f$. For radially decreasing functions $f$, the level sets are balls centered at the origin. Therefore, the condition reads 
$$
    \int_{B_R}\left( \Delta f(|x|)+f(|x|)^q\right)\geq 0, \quad \mbox{ i.e. } \quad R^{d-1}f'(R)+\int_0^R r^{d-1}f(r)^q\diff r\geq 0
$$
for every $R\in [0,R_0]$. We need to be careful because the condition on the level sets implies this new simpler condition only for those values of $R$ for which the ball $B_R$ is a level set of $f(|x|)$. If $f$ is not strictly decreasing, there exist $R_1, R_2$ for which $f(|x|)$ is constant over $[R_1, R_2]$ and we should not impose the new condition for $R \in (R_1, R_2)$. However, if this condition is satisfied at $R=R_1$ and $\Delta f(|x|) = 0$ for $x \in B_{R_2} \setminus B_{R_1}$ due to the positivity of $f$ we can obtain it also for $R\in (R_1,R_2)$.

Thus, we  characterise the condition  ``$x\mapsto f(|x|)$ belongs to $\mathcal{S'}(\R^d)$'' as ``$R^{d-1}f'(R)\geq -M(R)$, where $M(R)\coloneqq\int_0^R r^{d-1}f(r)^q\diff r$''. We can exploit this language in order to impose that $f$ is nonincreasing and also to compute its mass. Thus, combining everything, we can set the following optimal control problem that will allow us to find the value of $M_c^{sub} (d)$.
\begin{problem}[Optimal control problem]
We consider trajectories $(f,M):[0,R_0]\to \R^2$ satisfying
\begin{equation}\label{eq:Pontryagin}
\begin{cases} r^{d-1}f'(r)=-M(r)+\alpha(r),\\
			M'(r)=r^{d-1}f_+^q(r),\\
			f(0)=1,\\
			M(0)=0,\end{cases}
\end{equation}
where $\alpha:[0,R_0]\to \R_+$ is a nonnegative measurable control, subject to the constraint $\alpha(r)\in [0,M(r)]$ (a state-dependent constraint on the control). The quantity that we want to minimize is $M(R_0)$, where $R_0$ is the final horizon.
\end{problem}

Note that we wrote $f_+$ instead of $f$ in the equation defining $M'$ so that we allow $f$ to become negative (we do not impose state constraints). Moreover, it is important to understand the role of the control $\alpha$. This $\alpha$ is not the same as in \eqref{intuition alpha}, i.e. it does not represent the positivity of the equation, but only of its cumulated mass from $0$ to $r$. In some sense, this $\alpha$ is the antiderivative of the one in \eqref{intuition alpha} and we are only requiring that it is nonnegative instead of requiring that it is the antiderivative of a nonnegative function. Finally, the condition $\alpha\leq M$ has been imposed to guarantee that $f$ is nonincreasing. 
\label{prob:OC}

\paragraph{The solutions of \Cref{prob:OC}.}

Since our goal is to use necessary optimality conditions, in the form of Pontryagin's principle, on the optimal control problem that we just introduced, we first need to prove that it admits a solution.

\begin{lemma}\label{lem:minimizer_P_tilde}
The optimal control problem \ref{prob:OC} admits a minimizer that we call $f$.    
\end{lemma}

\begin{proof}[Proof of Lemma~\ref{lem:minimizer_P_tilde}]
Let $(f_n,M_n,\alpha_n)$ be a minimizing sequence. We can assume that $f_n$ is nonnegative, up to considering $(f_n)_{+}$, which is also a competitor and does not increase the mass. Moreover, $f_n(0)=1$ and $f_n$ is nonincreasing, so that $f_n$ is bounded in $BV(0,R_0)$ and therefore we can extract a subsequence which converges strongly to some $f$. On the same subsequence, we also have convergence of $M_n$ to a certain $M$. The controls $\alpha_{n}$ are bounded uniformly in $L^{\infty}$ by the constraint and the condition $M_n(R_0)\leq C$, which comes from the minimization. Thus, the weak convergence of $\alpha_n$ is enough to pass to the limit in the first equation. Thereby, we are left to prove that we have $f(0)=1$. This fact is not obvious since we only have a BV bound on $f_n$. However, the uniform $L^\infty$ bound on $f_n$ provides $M_n(r)\leq Cr^d$ and it also implies $\alpha_n(r)\leq Cr^d$ due to the constraint. This provides $0\geq f_n'(r)\geq -Cr$ and shows that the functions $f_n$ are equi-Lipschitz continuous. Hence, the initial value passes to the limit and it concludes the proof.
\end{proof}

We recall here the main tool that we intend to use, namely Pontryagin's maximum principle. We provide a very informal introduction to it, borrowing notations from the standard language used in optimal control theory. We refer for instance to \cite{Pontryagin87, Macki_Strauss82}.  Given a contol set $V$, a function $F \in C^0([0,T] \times \mathbb{R}^n \times V, \mathbb{R}^n)$, and $x \in \mathbb{R}^n$, we consider the Cauchy problem:
\begin{equation}
\dot{y}(t) = F(t, y(t), u(t)), \quad t \in [0,T], \quad y(0) = x .
\end{equation}
Here, $u$ is the control variable. Under suitable assumptions, the Cauchy problem admits a unique solution denoted $y_u$.  

Given a Lagrangian $L \in C^0([0, T] \times \mathbb{R}^n \times V, \mathbb{R})$ and a terminal gain function $g \in C^0(\mathbb{R}^n, \mathbb{R})$, we consider the optimal control problem:
\begin{equation}
\sup_{u\in V} J(u) \coloneqq \int_0^T L(t, y_u(t), u(t)) \, dt + g(y_u(T)).
\end{equation}
For $p \in \mathbb{R}^n$, we define the pre-Hamiltonian $\underline{H} : [0,T] \times \mathbb{R}^n \times V \times \mathbb{R}^n \to \mathbb{R}$ by:
\begin{equation}
\underline{H}(t, x, u, p) \coloneqq L(t, x, u) + p \cdot F(t, x, u).
\end{equation}
We then define the Hamiltonian $H : [0,T] \times \mathbb{R}^n \times \mathbb{R}^n \to \mathbb{R}$ by:
\begin{equation}
H(t, x, p) \coloneqq \sup_{u \in V} \{ L(t, x, u) + p \cdot F(t, x, u) \} = \sup_{u \in V} \underline{H}(t, x, u, p).
\end{equation}
We assume that $H$ is continuous, differentiable with respect to $x$ and $p$, and we denote the corresponding partial gradients by $\nabla_x H$ and $\nabla_p H$. Hence, the statement of Pontryagin's principle reads as follows.

\begin{theorem}[Pontryagin's principle]\label{lem:Pontryagin} If $u$ is an optimal control and $y \coloneqq y_u$ denotes the associated trajectory, then there exists $p$ such that:
\begin{equation}
H(t, y(t), p(t)) = \underline{H}(t, y(t), u(t), p(t)).
\end{equation} 
Furthermore, the pair $(y(\cdot), p(\cdot))$ is a solution to the Hamiltonian system
\begin{equation*}
\begin{cases}
\dot{p}(t) &= -\nabla_x H(t, y(t), p(t)), \\
\dot{y}(t) &= \nabla_p H(t, y(t), p(t)), 
\end{cases}
\end{equation*}
with the boundary conditions
\begin{align*}
y(0) &= x, \\
p(T) &=  g'(y(T)) \quad \text{(transversality condition)}.
\end{align*}
\end{theorem}

In our case we have $x=(f,M)$, so that we need to define a pair $(p_1,p_2)$ of dual variables for the Pontryagin's principle. The pre-Hamiltonian is given by $\underline{H}(r,x,\alpha,p)= p\cdot F(r,x,\alpha)$ where 
$$
F(r,x,\alpha)=\left( \f{-x_2+\alpha}{r^{d-1}},\,  (x_{1})_{+}^q r^{d-1}\right).
$$
The terminal gain function is $g(x_1,x_2)=-x_2$. To avoid the degeneracy at $r=0$ we apply Pontryagin's principle for $r\in[r_0 ,R_0]$, with any initial value at $r=r_0$ for arbitrary positive $r_0$.  The Hamiltonian can be computed and gives
\begin{align*}
    H(r,(f,M),p) &= \sup_{0\le \alpha\le M}p_{1}\left( \f{-M+\alpha}{r^{d-1}}\right) + p_2 r^{d-1}f_+^q\\
    &=  p_2 r^{d-1}f_+^q + (p_{1})_-\f{M}{r^{d-1}}. 
\end{align*}
Furthermore, the optimality implies the following conditions:
\begin{equation}\label{lien p alpha}
p_1>0\impl \alpha=M ; \quad p_1<0 \impl \alpha=0.
\end{equation}
Pontryagin's principle  yields 
\begin{equation}\label{eq:p1 and p2}
    p_{1}'(r) = - qp_2(r) r^{d-1} f_{+}(r)^{q-1}, \quad  p_{2}'(r) = - \f{(p_{1}(r))_{-}}{r^{d-1}},
\end{equation}
for $r\in[r_0, R_0]$, $p_{1}(R_0)=0$, $p_{2}(R_0)=-1$. Moreover, the ODE system~\eqref{eq:Pontryagin} is satisfied.  Our goal is to prove $f_{+}(r)\alpha(r)=0$ for every $r$. If such a statement is true we can conclude that the optimizer is indeed a {\it solution} of the Lane-Emden equation till the first moment where it touches $0$, i.e. that the minimizer is a rescaling of $\bar h$.

 We recall $f_{+}(0)>0$. If $f_{+}$ vanishes somewhere in $(r_0,R_0)$ we call $R_1$ the first time such that $f_{+}$ vanishes. Otherwise, we set $R_1=R_0$. If $R_1<R_0$, then $f_{+}$ must identically vanish after $R_1$, otherwise replacing it with $0$ would improve the mass $M(R_0)$. 

On $(R_1,R_0)$ we have $f_{+}(r)=0$ therefore $p_{1}(r)=0$ and $p_{2}(r)=-1$. By continuity, close to $r=R_1^-$ we have $p_2<0$ and hence $p_1'>0$, which implies $p_1<0$. If $p_1$ is strictly negative everywhere we obtain $\alpha(r)=0$ by \eqref{lien p alpha}. If not, let us consider $r_1=\sup\{r\in[r_0,R_1]\,:\, p_1(r)\geq 0\}$. Since $p_2$ is a nonincreasing function, we necessarily have $p_2(r_1)>0$, otherwie $p_2\leq 0$ on $[r_1,R_1$ and $p_1(r_1)<0$. But this implies that $p_2$ is strictly positive on $[r_0,r_1]$ and hence $p_1(r)>p_1(r_1)=0$ for all $r<r_1$. On this interval we then deduce $\alpha=M$. In particular, it is important to observe that we do not have non-trivial intervals on which $p_1=0$, where it would be impossible to deduce the value of $\alpha$.

\begin{proposition}[One-parameter family]\label{prop:family-gamma}
If $(f,M,\alpha)$ is a solution of the optimal control problem \ref{prob:OC}, with $f\geq 0$, then there exists a unique $\gamma\in[0,R_1]$ such that:
\begin{enumerate}
\item\label{item:flat-zone}
  $f(r)=1$ and $\alpha(r)=M(r)$ for $r\in[0,\gamma]$.
  \item There  exists a number $R_1\le R_0$ such that $f=0$ on $[R_1,R_0]$.
\item\label{item:ode-region}
  $\displaystyle \,f'(r)\,r^{d-1}+M(r)=0$
  and $\alpha(r)=0$ for $r\in(\gamma,R_1)$; in particular, the second-order ODE
  \begin{equation}
     \label{2ndOgamma}
   \,\left(f'(r)\,r^{d-1}\right)'+r^{d-1}f^q(r)=0   
  \end{equation}
  is satisfied on $(\gamma,R_1)$.
  \item  $f'(\gamma)=-\frac{\gamma}{d}$ where $f'$ denotes here the derivative of $f$ restricted to $[\gamma,R_1]$ (i.e. the right-derivative of $f$ at $\gamma$).
\end{enumerate}
The parameter $\gamma$ measures the length of the possible flat plateau where $f\equiv1$.
\end{proposition}

\begin{proof}
The Pontryagin's principle described above tells us that, for any $r_0<R_0$, the dual function $p_1$ is strictly positive on an interval $[r_0,\gamma)$ and strictly negative on $(\gamma,R_0]$. This implies $\alpha=M$ on $[r_0,\gamma)$ and hence on $(0,\gamma)$, and $\alpha=0$ after $r=\gamma$. In particular, we have that $f'(r)=0$ on  $(0,\gamma)$ and therefore $f=1$ on the same interval. Moreover, $M(\gamma)=\int_0^\gamma r^{d-1}f_+^q(r)\diff r= r^d/d$. Then, we also have  $\lim_{r\to\gamma^+} r^{d-1}f'(r)=-M(\gamma)=-\gamma^d/d$, which gives the condition on $f'(\gamma)$.  

Afterwards, we define $R_1$ as the first point where $f$ vanishes.
Since we choose $R_0$ large enough, from Proposition~\ref{prop:def_R0} we know that we have  $R_1<R_0$. 
Necessarily, $f$ is identically $0$ after $R_1$ since $f$ is nonincreasing and nonnegative.
\end{proof}

It is important to observe that for every $\gamma\geq 0$ there exists a unique function $f_\gamma$ satisfying the ODE \eqref{2ndOgamma} together with the two initial conditions $f(\gamma)=1$ and $f'(\gamma)=-\gamma/d$. This function is defined on $(\gamma,\infty)$. If it vanishes somewhere, we call the first vanishing point  $R(\gamma)$. Otherwise, if it does not vanish, we set $R(\gamma)=+\infty$.

\paragraph{Choice of the horizon $R_0$.}

In this step we justify the fact that we can choose a universal value of $R_0$ such that, for sure, the vanishing point $R_1$ satisfies $R_1<+\infty$. We call $\bar M$ the mass associated with an arbitrary compactly supported competitor, for instance, the one corresponding to the function $\bar h$.
\begin{proposition}\label{prop:def_R0}
 We consider the functions $\{f_{\gamma}\}_{\gamma}$ that we just defined.
Then, there exists $R_0>0$ such that $R(\gamma)\leq R_0$ for all $\gamma$ such that
$$
\int_{0}^{R(\gamma)} f_{\gamma}(r)^{q}r^{d-1}\diff r \le \bar{M}. 
$$
\end{proposition}

\begin{proof}
First, we prove that, for all $\gamma$, $R(\gamma)<+\infty$. We make the change of variables $f_{\gamma}(r) = \frac{u_{\gamma}(r^{d-2})}{r^{d-2}}$ for some nonnegative $u_{\gamma}$. Since $f_{\gamma}(\gamma)=1$ and $f_{\gamma}'(\gamma)=-\frac{\gamma}{d}$ we obtain $u_{\gamma}\left(\gamma^{d-2}\right)=\gamma^{d-2}$ and $u_{\gamma}'(\gamma^{d-2})=1-\gamma^2\frac{1}{d(d-2)}$. We extend $u_\gamma$ on $[0,\gamma^{d-2}]$ through $u_\gamma(s)=s$. The equation for $u_\gamma$ on $(\gamma^{d-2},R(\gamma)^{d-2})$ reads
\begin{equation}\label{eq:eq_on_u}
    (d-2)^2 u_{\gamma}''(s) + \frac{u_{\gamma}(s)^{q}}{s^{2}}=0.
\end{equation}
In particular, $u_{\gamma}$ is concave. We compute the mass in terms of $u_{\gamma}$ and we obtain
$$
    \int_{0}^{R(\gamma)} f_{\gamma}(r)^{q}r^{d-1}\diff r= \int_{0}^{R(\gamma)^{d-2}}\frac{u_{\gamma}(s)^{q}}{s}\diff s.
$$
Therefore, 
from the initial conditions on $u_{\gamma}$ and the fact that $u_{\gamma}$ is concave we have clearly that either $u_{\gamma}$ vanishes in finite time -- and therefore $R(\gamma)<+\infty$ --  or its derivative remains nonnegative and in that case $u_{\gamma}\ge \gamma^{d-2}$ and $R(\gamma)=+\infty$. The second option is in contradiction with the finiteness of the integral. Moreover, in order to estimate $R(\gamma)$ uniformly it is enough to find $s_0>0$ and $\varepsilon_0$ such that $u_\gamma(s_0)>\varepsilon_0$ independently of $\gamma$. 
From the concavity, we have that $u_\gamma(s)\leq s$. If we insert this observation into \eqref{eq:eq_on_u} we obtain $|u_\gamma''(s)|\leq Cs^{q-2}$, so that $u_\gamma''$ is uniformly integrable independently of $\gamma$. Hence, if $\gamma$ is small such that $u_\gamma'(\gamma^{d-2})>1/2$ we obtain $u_\gamma'>1/4$ on $[0,s_0]$. Thereby, $u_\gamma(s_0)>s_0/4$ for a value of $s_0$ which can be explicitly computed. If $\gamma$ is large, then $u_\gamma(s_0)=s_0$.
\end{proof}

\paragraph{The one-parameter family of competitors.}

In the previous steps we have proven that the only reasonable competitors are those of a one-parameter family constructed in the following way: Choose $\gamma>0$, take $f=1$ on $[0,\gamma]$, then solve the ODE with initial condition $f(\gamma)=1$ and $f'(\gamma)=-\frac\gamma d$, and stop at a point $R(\gamma)$ where $f(R(\gamma))=0$, then extend by $0$.
 
In order to conclude with our argument, it only remains to show that the best $\gamma$ is $\gamma=0$. Hence, we study the dependence in $\gamma$. Note that we present our pathway as a reduction to a one-parameter problem since we believe that it is easier to understand, but, as we will see later, the optimality conditions in terms of $\gamma$ could also be obtained by only looking at the Pontryagin principle. Finally, note that in case (but we will prove it is not the case) the optimal value of $\gamma$ would not be $0$, then we can simply conclude that we have enlarged too much the original problem by allowing the functions to be constant equal to $1$ on some interval near the origin. Indeed, these functions fail to satisfy the original problem $ \Delta f + f^\frac{d}{d-2}\ge 0$: close to the point where $f$ becomes less than one, $\Delta f$ is singular and negative. The only function satisfying the original constraint is the one with $\gamma=0$. 


In the following, we write $f_{\gamma}(r) \coloneqq f(\gamma,r)$ in order to clarify the dependence in both parameters and in order to avoid confusion. 
%
We define $\mathrm{M}(\gamma)$ the mass associated with the function $f (\gamma, \cdot)$. Using the equation \eqref{eq:Pontryagin} satisfied by $f (\gamma, \cdot)$ we have
$$
    \mathrm{M}(\gamma) =-   R_{\gamma}^{d-1} \partial_r f (\gamma,R_{\gamma} ).
$$
This leads to the following characterisation of the problem.

\begin{lemma}
For any fixed $\gamma$, the function $w (r) \coloneqq \partial_\gamma f (\gamma , r)$ satisfies the ODE system 
\begin{equation}
\label{system for v}
\begin{cases}
    (r^{d-1}w'(r))'+qf^{q-1}(\gamma , r)r^{d-1}w(r)=0 & \mbox{ in }(\gamma,R(\gamma)),\\
    w(\gamma)=\frac\gamma d , &\\
    w'(\gamma)=0. &
\end{cases}
\end{equation}
Moreover, if $\gamma\geq 0$ is such that $\mathrm{M}'(\gamma) =0$ then, $w$ also satisfies $w'(R(\gamma)) = 0$.
\end{lemma}

Before proving this result, let us observe that $w(r)$ is only defined for $r \in (\gamma , R(\gamma))$. Moreover, on this interval, we note that $w$ coincides with the function $p_2$ obtained as a dual variable in the Pontryagin principle.

\begin{proof}
First of all, we can obtain the ODE satisfied by $w$ if we differentiate the ODE  \eqref{2ndOgamma}.  We then use
$$
    f(\gamma,\gamma)=1,\quad \partial_r f (\gamma,\gamma)=-\frac\gamma d.
$$
If we differentiate them in terms of $\gamma$ we obtain that  
$$
    w(\gamma)+ \partial_r f (\gamma,\gamma)=0,\quad w'(\gamma) + \partial_{rr} f (\gamma,\gamma) = -\frac 1d.
$$
From $\partial_{rr} f + \frac{d-1}{r} \partial_r f + f^q = 0$ we have $\partial_{rr} f (\gamma,\gamma) =  -\frac{1}{d}$ and in particular $w'(\gamma)=0$.

For the second part of the statement, 
%
%
%
we  compute  $\frac{\diff}{\diff \gamma} \mathrm{M}(\gamma)$ and we obtain 
\begin{align*}
    \frac{\diff}{\diff \gamma} \mathrm{M}(\gamma)& =- \left( \frac{\diff R(\gamma)}{\diff \gamma} \right) (d-1)R(\gamma)^{d-2} \partial_r f (\gamma,R(\gamma)) \\
    & \quad \, -  R(\gamma)^{d-1} \left[ \partial_{r \gamma} f (\gamma,R(\gamma)) + \partial_{rr} f (\gamma, R(\gamma)) \left( \frac{\diff R(\gamma)}{\diff \gamma} \right) \right]
    \\ &= \left( \frac{\diff R(\gamma)}{\diff \gamma} \right) R(\gamma)^{d-1} \left[ \frac{d-1}{R(\gamma)} \partial_r f (\gamma,R(\gamma)) + \partial_{rr} f (\gamma,R(\gamma)) \right] - R(\gamma)^{d-1} \partial_{r \gamma} f (\gamma,R(\gamma)).
\end{align*}
We use the equation on $f$ and we observe that
$$
    \frac{d-1}{R(\gamma)} \partial_r f (\gamma,R(\gamma)) + \partial_{rr} f (\gamma,R(\gamma))= -f(\gamma,R(\gamma))^q = 0.
$$
Therefore, we deduce
$$
    \frac{\diff}{\diff \gamma} \mathrm{M}(\gamma)= -  R(\gamma)^{d-1} \partial_{r \gamma} f (\gamma,R(\gamma))=-  R(\gamma)^{d-1} w'(R(\gamma)).
$$
With this we conclude the second part of the statement.
\end{proof}

\begin{proposition}\label{ChatGPT}
If $\gamma>0$ there is no solution to \eqref{system for v} with $w'(R(\gamma))=0$.
\end{proposition}

\begin{proof}
The key point in the proof is the fact that the function $z$ defined via
$$z(r)=rf'(r)+(d-2)f(r)$$
is a solution of the same second-order linear ODE solved by $w$. This - which can of course be checked by explicit computations - comes from the fact that $\lambda^{d-2} f(\lambda\cdot)$ solves the same nonlinear equation as $f$, and $z$ woincides with the derivative w.r.t. $\lambda$ of this expression at $\lambda=1$.

We can also compute
\begin{equation}\label{tildew'}z'(r)=f'(r)+rf''(r)+(d-2)f'(r)=rf''(r)+(d-1)f'(r)=-rf^q(r)
\end{equation}
and observe that we have $$
z'(R(\gamma))=-R(\gamma)f^q(R(\gamma))=0. 
$$
On the other and, we have $z(R(\gamma)=R(\gamma)f'(R(\gamma)\neq 0$. The function $w$ as well does not vanish at $r=R(\gamma)$ otherwise it would be a solution of a linear second-order ODE satisfying homogeneous Cauchy conditions for $w$ and $w'$ at $r=R(\gamma)$ and, by uniqueness, we would have $w=0$ (but we have $w(\gamma)\neq 0$). Thus, setting $c:=w(R(\gamma))/z(R(\gamma))$, the functions $cz$ and $w$ solve the same second-order linear equation and their values and the values of their derivatives agree at $r=R(\gamma)$. Hence they coincide. Yet, we have $w'(\gamma)=0$ while (using \ref{tildew'}) $cz'(\gamma)=-c\gamma f^q(\gamma)=-c\gamma\neq 0$. This contradiction concludes the proof. 
\end{proof}

\section{Proof of the Li--Yau and Aronson--B\'enilan estimates}\label{sec:endproof}

In this section, we complete the proof of the main  Li--Yau and Aronson--B\'enilan estimates stated in Theorems~\ref{thm:Li--Yau} and~\ref{thm: Aronson--B\'enilan}. We recall that these estimates depend on three cases:
\begin{enumerate}
\item Small mass case (\ref{2d_small_mass} and \ref{2d_small_mass_multid}),
\item Subcritical mass case (\ref{2d_full_mass} and \ref{2d_full_mass_multid}),
\item Critical mass case (\ref{2d_exact_mass} and \ref{2d_exact_mass_multid}).
\end{enumerate}
We begin here with the small mass case, where the proofs are simpler and yield  explicit constants.

\subsection{Small mass case}
In the small mass case, we obtain a threshold mass $\varepsilon_d$ depending on the dimension for which the Li--Yau and Aronson--B\'enilan type estimates hold. In dimension $2$ we provide an estimate on $\varepsilon_2$ ($\varepsilon_2\ge 2\pi$). For dimension $d>2$ we can also compute a threshold via a complicated formula. However, we want to stretch that the value for $\varepsilon_d$ may not be optimal.

We start our analysis showing that $Q(u) \leq C_3(M)|\delta(t)|$ for dimension $d>2$. The case $d=2$ was already introduced in Proposition~\ref{prop:estimate_Linf_difference}:

\begin{proposition}\label{prop:split_cases_smallmass}
Let $\delta(t)= \min\Delta v(t,\cdot)$. Assume we are under the small mass assumptions of \ref{2d_small_mass} (dimension 2) or \ref{2d_small_mass_multid} (dimension $d>2$). Then for each $t\in(0,T)$, the following estimates hold:
$$Q(u)\leq C_{3}(M)\left|\delta(t)\right|.$$
Here $\lim_{M\to 0}C_3(M)=0$ and, if $d=2$, we have
$$C_3(M)=\f{M}{2\pi}(1+C_2(M))=\frac{2M}{4\pi-M},$$ whereas for $d>2$ one may take the constant given in the remark following \Cref{prop:estimate_Linf_difference_full_mass_multid}.
\end{proposition}
\begin{proof}
The estimate in the case $d=2$ was proven in Proposition~\ref{prop:estimate_Linf_difference}. The general case is a consequence of Proposition~\ref{prop:linf_very_small_mass} and Proposition~\ref{prop:estimate_Linf_difference_full_mass_multid} combined with the inequalities $\bar{\delta}\le(1+C_2(M))|\delta|$ and $\|\rho\|_{L^{2}}\le M^{1/2}\|\rho\|_{L^{\infty}}^{1/2}$.
\end{proof}

\begin{proof}[Proof of \ref{2d_small_mass} and \ref{2d_small_mass_multid}]
Let $\delta(t)=\min\Delta v(t,\cdot)$. Combining Proposition~\ref{prop:preliminary_aronson_benilan} and Proposition~\ref{prop:split_cases_smallmass}, we obtain
$$
\delta'(t)\ge\left(1-\f{(d-1)^2}{2d}C_3(M)^2\right)\delta(t)^2.
$$
Assume now that $\f{d-1}{\sqrt{2d}}C_3(M)<1$, so that the factor above is positive. Then,
$$
\delta(t)\ge-\f{C_{\mathrm{AB}}(M)}{t},
\qquad\text{where}\qquad
C_{\mathrm{AB}}(M)=\left(1-\f{(d-1)^2}{2d}C_3(M)^2\right)^{-1}.
$$
Since $M<M^{sub}_c(d)$, one has $\delta(t)<0$, and Proposition~\ref{prop:linf_very_small_mass} gives
$$
\|\rho_t\|_{L^\infty}\le C_2(M)|\delta(t)|\le\f{C_2(M)C_{\mathrm{AB}}(M)}{t}.
$$
Finally, writing $p_t=\log\rho_t$ if $d=2$ and $p_t=\f{m}{m-1}\rho_t^{m-1}$ if $d>2$, we have $\Delta p_t=\Delta v_t-\rho_t$, and therefore
$$
\Delta p_t\ge-\f{(1+C_2(M))C_{\mathrm{AB}}(M)}{t}.
$$
\end{proof}

Thus, we have established the Li--Yau and Aronson--B\'enilan estimates in the case of small mass $M < \varepsilon_d$. Before handling the case of subcritical and critical mass cases, we prove here that the solutions can be considered smooth up to an approximating scheme, even when the initial condition is merely a measure, which is necessary in order to justify all the above computations in the small mass case. 
We prove that if a sequence of initial conditions $\rho_{n,0}$ has mass below the small-mass threshold, then we can pass to the limit in the corresponding sequence of smooth solutions $\rho_n$ to the Keller--Segel system. The limit function $\rho$ is a global solution in a weak sense of the same system with the initial data $\rho_0$.

We use the following weak formulation: for every
$\varphi\in C_c^\infty([0,\infty)\times\R^d)$,
$$
\int_0^\infty\int_{\R^d}
\left(\rho\,\partial_t\varphi+\rho^m\Delta\varphi
+\rho\nabla u\cdot\nabla\varphi\right)\diff x\diff t
+\int_{\R^d}\varphi(0,x)\diff\rho_0(x)=0,
$$
where $u=\Gamma*\rho$ and the terms are locally
integrable up to $t=0$. The solutions constructed below preserve
their mass and have a narrowly continuous representative.

\begin{proposition}\label{prop:existence-small-mass}
Let $\rho_{n}$ be a sequence of solutions of the Keller--Segel system 
\[
\begin{cases}
\partial_t\rho_n = \Delta\rho_n^m - \mathrm{div}\left(\rho_n \nabla u_n\right),\\
-\Delta u_n = \rho_n,
\end{cases}
\]
on $(0,T)\times \R^d$
with initial conditions $\rho_{n,0}$ satisfying:
\begin{enumerate}
\item $\rho_{n,0}$ is smooth,
\item $\rho_{n,0}\to\rho_0$ narrowly and 
$\rho_{n,0}(\R^d)=\rho_0(\R^d)=M<\varepsilon_d$ where $\varepsilon_d$ is the small-mass thresholds from Theorems~\ref{thm:Li--Yau}  and~\ref{thm: Aronson--B\'enilan},.
In dimension two, assume also a uniform logarithmic-moment bound.
\item $\rho_n\rightharpoonup \rho$ weak-* in $\mathcal{M}([0,T]\times \mathbb{R}^d)$.
\end{enumerate}
Then, $\rho$ is a weak solution  of the same Keller--Segel system such that:
\[
\lim_{t\to 0}\rho(t,\cdot)   =   \rho_0
\quad \text{(in the sense of measures, $\mathcal{M}(\mathbb{R}^d))$}.
\]
Moreover, $\rho$ satisfies the conclusions of \ref{2d_small_mass} if $d=2$ or \ref{2d_small_mass_multid} if $d>2$. 
\end{proposition} 


\begin{proof}
We divide the proof in several steps.

\begin{enumeratesteps}
    
\step[Uniform existence and first estimates]
By the small-mass estimates from \ref{2d_small_mass}    (for $d=2$) or \ref{2d_small_mass_multid} (for $d>2$), we observe that
\[
    \|\rho_n(t,\cdot)\|_{L^\infty(\mathbb{R}^d)}   \lesssim   \f{C}{\,t\,}
\]
for $t>0$, with a uniform constant $C$ independent of $n$. Furthermore, the total mass satisfies $\rho_n(\mathbb{R}^d)=M$ (constant in time). Therefore, it follows that
\[
    \rho_n(t,\cdot)  \text{ is uniformly bounded in }L^\infty\left((\varepsilon,T); L^\infty(\mathbb{R}^d)\right)\cap L^\infty\left((0,T); L^1(\mathbb{R}^d)\right)
\]
for any $0<\varepsilon<T$.  

\step[Compactness and passage to the limit]
Beyond these uniform bounds, the Keller--Segel equation provides additional regularity. One can use the structure $\partial_t\rho_n - \Delta\rho_n^m + \mathrm{div}(\rho_n \nabla u_n)=0$ (with $-\Delta u_n = \rho_n$) in order to obtain further compactness (for instance by computing the dissipation of the free energy). Therefore one shows that $\{\rho_n\}$ is relatively compact in a topology strong enough to pass to the limit in the nonlinearities (e.g.\ strong convergence in a suitable Sobolev space on compact time intervals away from $t=0$).

Hence, up to extracting a subsequence (not relabeled), we deduce
\[
\rho_n   \longrightarrow   \rho
\quad
\text{strongly in an appropriate sense on }(0,T)\times \mathbb{R}^d.
\]
It follows that $\rho$ is a (weak) solution to
\[
\partial_t\rho
   =  
\Delta\rho^m
 - \mathrm{div}\left(\rho\,\nabla u\right),
\quad
-\Delta u = \rho,
\]
satisfied almost everywhere in $(0,T)\times \mathbb{R}^d$.

\step[Recovering the initial condition $\rho(0,\cdot)=\rho_0$]
We are now left with the last subtle part, which is that the estimate $\|\rho_n\|_{L^\infty(\varepsilon,T; L^\infty(\mathbb{R}^d))}\lesssim 1$ blows up as $\varepsilon\to0$. In order to identify the limit at $t=0$, we need time equicontinuity in a negative Sobolev norm, which allows us to apply the Arzel\`a--Ascoli theorem in time.

To establish time equicontinuity in some $H^{-s}$ norm, we want to obtain
\[
\|\partial_t\rho_n\|_{L^p(0,T;H^{-s}(\mathbb{R}^2))}   \le   C,
\]
independently of $n$. We split the argument into estimating $\Delta\rho_n^m$ and $\mathrm{div}(\rho_n \nabla u_n)$ in some dual space.

The case $d=2$ can be done via an easier argument since we have that $m=1$, i.e. linear diffusion. In particular, $\rho_n$ is uniformly bounded in $L^{\infty}(0,T; L^1(\R^2))$ due to the mass conservation. This implies that for each $t$ the distribution $\Delta\rho_n$ is uniformly bounded in $H^{-s_1}(\R^2)$ for some $s_1>0$. For the case  of dimension $d>2$, where $m>1$, we require a more elaborate argument. We observe that 
$$
    \rho_n^{m}(t,\cdot)=\rho_n^{m-1}(t,\cdot)\rho_n (t,\cdot)\lesssim \f{1}{t^{m-1}}\rho_{n}(t,\cdot).
$$
Therefore, we have that
$$
    \|\rho_n^m\|_{L^{p_1}(0,T; L^{1}(\R^d))}\lesssim \left\|\f{1}{t^{m-1}}\right\|_{L^{p_1}(0,T)}\|\rho_{n}\|_{L^{\infty}(0,T; L^{1}(\R^d))}.
$$
Since $0<m-1<1$, it is always possible to choose $p_1$ close enough to $1$ such that the previous quantity is bounded. We deduce that there exists $p_1>1$ and $s_1>0$ large enough such that
$$
\|\Delta \rho_n^m\|_{L^{p_1}(0,T; H^{-s_1}(\R^d))}\lesssim C.
$$

For the convective term given by $\mathrm{div}(\rho_n\,\nabla u_n)$ it is enough to estimate $\rho_n \nabla u_n$ in a suitable $L^p(0,T; L^1(\R^d))$ norm.  We have
\[
\|\rho_n\nabla u_n\|_{L^p(0,T;L^1(\mathbb{R}^d))}
   \le  
\|\rho_n\|_{L^\infty(0,T;L^1(\mathbb{R}^d))}  \|\nabla u_n\|_{L^p(0,T;L^\infty(\mathbb{R}^d))}.
\]
The first factor is uniformly bounded due to the mass conservation. In order to control the term $\|\nabla u_n\|_{L^p(0,T;L^\infty(\R^d))}$, we take advantage of the form of the Newtonian potential. It  provides that 
\begin{align*}
|\nabla u_n(t,x)|
 &  \lesssim   
\int_{|x-y|\le 1}\f{\rho_n(t,y)}{|x-y|^{d-1}}\,\mathrm{d}y
    +  
\int_{|x-y|\ge 1}\f{\rho_n(t,y)}{|x-y|^{d-1}}\,\mathrm{d}y\\
 &  \lesssim  
\int_{|x-y|^{d-1}\le 1}\!\f{\rho_n(t,y)}{|x-y|}\,\mathrm{d}y
    +  
\|\rho_n(t,\cdot)\|_{L^1}
   \lesssim  
\|\rho_n(t)\|_{L^{d+\varepsilon}} + M.
\end{align*}
Since $\rho_n(t,\cdot)\leq c/t$ and also $\int\rho_n=M$, we get that
\[
\|\nabla u_n(t)\|_{L^\infty}
   \lesssim  
\f{1}{\,t^{(d-1+\varepsilon)/(d+\varepsilon)}\,}   +  1.
\]
Therefore, we can choose $1<p<(d+\varepsilon)/(d-1+\varepsilon)$ so that
\[
\|\nabla u_n\|_{L^p(0,T;L^\infty(\R^d))}
   \lesssim   
1.
\]
Hence, we conclude that $\mathrm{div}(\rho_n \nabla u_n)$ is bounded in $L^p(0,T;H^{-s_2}(\mathbb{R}^d))$ for some $s_2>0$. Summarizing, $\partial_t\rho_n$ is uniformly bounded in $L^p(0,T;H^{-s}(\mathbb{R}^d))$, implying equicontinuity in time in the sense of $H^{-s}$. This allows to obtain uniform convergence as $n\to\infty$, which implies the convergence of the initial datum, at least in the distributional sense. But since the limit has the same mass, this can be upgraded to a convergence in $\mathcal{M}(\R^d)$.\qedhere
\end{enumeratesteps}
\end{proof}

The next two subsections handle the case of subcritical and critical mass respectively. 

\subsection{Subcritical mass case}

When the total mass $M$ is strictly smaller than the critical mass $M_c$, but not necessarily small in the sense of the previous subsection, we cannot directly conclude a Li--Yau or Aronson--B\'enilan type estimate from the  bound
\[
    \f{\diff}{\diff t}\,\delta   \ge   \delta^2 - C\,\delta^2.
\]
In this case, we do not have an explicit value for the constant $C > 0$ and if it is too large,  the coefficient in front of $\delta^2$ might be negative.  
To obtain the result, we assume the free energy of the initial condition to be finite.  Under this assumption, we can derive additional integrability (e.g.\ via HLS or Log-HLS inequalities) of $\rho$.  
Therefore, we can refine the differential inequality for $\delta$, in order to obtain that 
\[
\f{\diff}{\diff t}\,\delta   \ge   \delta^2 \left[1 - o(1)\right].
\]
From this inequality on $\delta (t)$, we can  recover the Li--Yau or Aronson--B\'enilan estimates. The price to pay is that this differential inequality yields a weaker control on $\delta$.

\begin{proposition}\label{prop:split_cases_fullmass}
Let $\delta(t)= \min\Delta v(t,\cdot)$. Assume we are under the small mass assumptions of~\ref{2d_full_mass} (dimension $2$) or~\ref{2d_full_mass_multid} (dimension $d>2$). Then for each $t\in(0,T)$, the following estimates hold:
\begin{enumerate}
\item Dimension $d=2$:
\[
Q(u) \lesssim C_1(T) \f{|\delta(t)|}{\sqrt{\log(|\delta(t)|+e)}}.
\]
where $C_1(T)$ depends on $T$.   

\item Dimension $d>2$:
\[
Q(u)  \le  
C_{2}|\delta(t)|^{\f{d^2+4}{d(d+2)}},
\]
 where $C_{2}$ is a constant independent of $T$.
\end{enumerate}
\end{proposition}

\begin{proof}
We divide the proof into the case of dimension $d=2$ and $d > 2$.

\begin{enumeratesteps}
    \step[Dimension $d=2$]
We start with the case $d=2$. Under \ref{2d_full_mass}, we assume that the free energy and a logarithmic moment of the initial condition are bounded.  Therefore, the solution $\rho(t)$ maintains bounded free energy (since the free energy is nonincreasing in time) and bounded logarithmic-moment (however the logarithmic-moment bound depends on $T$). The latter bound can be obtain similarly than in~\cite{Fernandez_Mischler16, Calvez_Corrias08}, see also Theorem~\ref{thm:new_proof_hls}.

From the sharp logarithmic Hardy--Littlewood--Sobolev (Log-HLS) inequality (see for instance \cite{Carlen_Loss92,Beckner93} or \cite[Lemma~0.3]{Dolbeault_Perthame04}, or Theorem \ref{thm:new_proof_hls}), we deduce
\[
\int_{\R^2} \rho(t,x)\log\rho(t,x) \mathrm{d}x   \lesssim C(T)
\quad\text{for all }t\in(0,T).
\]
Using the bound on the logarithmic moment of $\rho_t$
and partitioning the domain in the regions where $\rho\le e$ or $\rho\ge e$, we also obtain
\[
\int_{\R^2}\rho  \log(\rho + e)
   \lesssim  
C(T)
\]
as well.
On the other hand, we already know from Proposition~\ref{prop:linf_estimate} that we have
$
\|\rho_t\|_{L^\infty}
   \lesssim  
C\,\lvert\delta(t)\rvert.
$
Hence, combining both results, we recover
\begin{align*}
\|\rho_t\|_{L^2} & = \sqrt{\int_{\R^2}\rho^2} =  
\sqrt{\int_{\R^2}\rho\,\log(\rho+e)  \f{\rho}{\log(\rho+e)}}
   \lesssim  \sqrt{C(T)}  \sqrt{\left\|\f{\rho}{\log(\rho+e)}\right\|_{L^\infty}} \\
   & \lesssim  
\sqrt{C(T)}  \sqrt{\f{\lvert\delta(t)\rvert}{\log(\lvert\delta(t)\rvert+e) }}.
\end{align*}
Applying Proposition~\ref{prop:estimate_Linf_difference_full_mass_multid}, we obtain the stated inequality
\[
Q(u)
   \lesssim  
\sqrt{1+T}  \f{\lvert\delta(t)\rvert}{\sqrt{\log(\lvert\delta(t)\rvert + e)}}.
\]

\step[Dimension $d > 2$]
Concerning the higher dimensional case $d>2$, we do not need to assume that a logarithmic moment is bounded. The Hardy-Littlewood-Sobolev inequality (see for instance \cite{Blanchet_Carrillo_Laurencot09, Calvez_Carrillo_Hoffmann17}) allows us to prove a uniform bound on the $L^m$ norm of $\rho_t$, independently of $T$. Since $\|\rho_t\|_{L^\infty}\le C|\delta(t)|$ from Proposition~\ref{prop:linf_estimate} we deduce that 
$$
    \|\rho\|_{L^2}=\sqrt{\int_{\R^d}\rho^2\diff x}\lesssim \sqrt{\|\rho\|_{L^{\infty}}^{2-m}}\lesssim |\delta|^{(2-m)/2}= |\delta|^{1/d} .
$$
Once more we use Proposition~\ref{prop:estimate_Linf_difference_full_mass_multid} and we deduce
\[
Q(u)  \lesssim |\delta|^{\f{d^2+4}{d(d+2)}} .
\]
In particular, we observe that the exponent is strictly smaller than $1$ for $d>2$. \qedhere
\end{enumeratesteps}
\end{proof}

Thanks to the auxiliary result given by \Cref{prop:split_cases_fullmass} we can now proceed to prove the main result.

\begin{proof}[Proof of \ref{2d_full_mass} and  \ref{2d_full_mass_multid}]
Let $\delta(t,\cdot)= \min \Delta v(t,\cdot)$.
By Proposition~\ref{prop:split_cases_fullmass}, and Proposition~\ref{prop:preliminary_aronson_benilan} we deduce a differential inequality of $\delta(t)$. Let us cover the cases $d=2$ and $d > 2$.

\begin{enumeratesteps}
\step[Dimension $d=2$]  
In the case $d=2$, this differential inequality reads for a constant $C(T)$ depending on $T$:
\[
  \f{\diff}{\diff t}\,\delta(t)
   \ge \delta^2(t)\,\left(1 - \f{C(T)}{\log\left(|\delta(t)| + e\right)}\right).
\]
We split into two subcases:

\begin{enumerate}
\item \emph{If $\delta(t)$ becomes large:} When $\left|\delta(t)\right|\gg e^{C(T)}$, then 
  $\log\left(|\delta(t)|+e\right)\approx \log|\delta(t)| \gg C(T)$.  Hence
  \[
    1 - \f{1+T}{\log\left(|\delta(t)| + e\right)}
     \ge  \widetilde{C}
  \]
  for some $\widetilde{C} > 0$.  Therefore, we have
  \[
    \f{\diff}{\diff t}\,\delta(t)
     \ge  \widetilde{C}\,\delta^2(t),
  \]
  which again forces $\delta(t)$ to be larger than $- \widetilde{C}/t$.

\item \emph{If $\delta(t)$ stays bounded:} Suppose $|\delta(t)| \le e^{T}$. 
  Then automatically $\delta(t) \ge -e^{T}$.  
\end{enumerate}

Thus, in either subcase, we have
$\delta(t) \ge -C(T)\left(1+\f{1}{t}\right)$for some  constant $C(T)$ depending on the final time $T$. Thereby, with the same argument already presented in  Proposition~\ref{prop:linf_estimate}, we conclude that $\|\rho\|_{L^\infty}$ remains bounded and it is such that
\[
  \|\rho_t\|_{L^\infty}
   \lesssim  C(T)\left(1+\f{1}{t}\right),
\]
for $t\in(0,T)$. This completes the proof of \ref{2d_full_mass}.

\step[Dimension $d>2$]
In dimension $d>2$, the differential inequality reads for a constant $C$:
$$
\f{\diff}{\diff t}\ge \delta^2(t)\left(1-C\delta^{\f{2(d^2+4)}{d(d+2)}-2}\right).
$$
Since $\f{2(d^2+4)}{d(d+2)}-2<0$ for $d>2$ we conclude in a similar way as in the two dimensional case. This completes the proof of \ref{2d_full_mass_multid}. \qedhere
\end{enumeratesteps}
\end{proof}

\subsection{Critical mass case}

When the total mass $M$ equals the critical mass $M_c$, it is no longer possible to rely solely on the logarithmic or HLS inequalities to obtain uniform integrability bounds on $\rho$. In fact, if we assume only bounded free energy, we cannot automatically derive an $L^m$ (or $L \log L$ in the case of dimension~2) control that remains uniform in time. Instead, our idea is to identify a compact set $K$ such that any Keller--Segel solution with bounded free energy is forced to remain in $K$. On this compact set, the stronger bounds needed to conclude the proof of the Li--Yau /  Aronson--B\'enilan estimate (e.g.\ an $L^\infty\lesssim |\delta|+1$ type estimate) become available.

We explain here the way we prove \ref{2d_exact_mass} and \ref{2d_exact_mass_multid}  and we will show later the technical proofs which we need to carry out this strategy.

\begin{enumeratesteps}
    \step Our initial datum $\rho_0$ is such that it is not a Dirac mass. We approximate it with smooth initial data from which we can start a smooth solution. Then, we will find a compact set $K$ in the space of finite measures such that $K$ does not contain any Dirac mass and such that the evolution of the Keller--Segel stays inside $K$ at least for some time $T_0>0$. The set $K$ and the time $T_0$ are allowed to depend on the initial datum $\rho_0$, but not on its approximation (for instance, they can depend on the moments of $\rho_0$, or similar quantities). 
    \step On $K$ we know (cf. Proposition~\ref{prop:linf_estimate_critical}) that we can estimate $\|\rho\|_{L^\infty}$ in terms of $|\delta[\rho]|$. In order to prove existence we only need to handle the case where $|\delta[\rho]|$ is very large, so the fact that the estimate is $\|\rho\|_{L^\infty}\leq C(|\delta[\rho]|+1)$ instead of $\|\rho\|_{L^\infty}\leq C|\delta[\rho]|$ is not an issue. Along the evolution we also know that $\mathcal F$ decreases, thus it stays bounded. We need to estimate the entropy (for $d=2$) or the $L^m$ norm (for $d>2$) of $\rho$ in terms of $\mathcal F$, which is in general not possible when $M=M^c$. Yet, it is possible to do so when restricting to $\rho\in K$ and this is the object of Subsections~\ref{sec:Estimates entropy} and~\ref{sec:Estimates Lm}. With these ingredients it is possible to prove also in this setting the differential inequality $\delta'\geq \delta^2(1+o(1))$, which provides existence of the solution.
    \step If $T_0<+\infty$ we need to iterate the same construction after $t=T_0$. This can be done starting from the measure that we reached at such a time, but we need to prove that the existence time has not worsened. If $T_0=+\infty$ we obtain that the estimates are global in time and in particular the same $\delta'\geq \delta^2(1+o(1))$ holds on the whole real line $t\in\R_+$. This implies $\delta(t)\geq -C(1+1/t)$.
\end{enumeratesteps}
We start by explaining how to choose the set $K$.

\subsubsection{\texorpdfstring{$H_\lambda[\rho_{0}]<+\infty$}{Hlambda[rho0]<+infty} in dimension 2}

The reader has maybe noticed that in Theorem~\ref{thm:Li--Yau} for dimension $2$, one of the statements requires the finiteness of an auxiliary functional $H_{\lambda}[\rho_0]$ for some $\lambda>0$. We recall the definition
\[
H_{\lambda}[\rho]
   =   
\int_{\mathbb{R}^2}\left(\sqrt{\rho} - \sqrt{u_\lambda}\right)^2\,u_{\lambda}^{-1/2},
\]
where $u_{\lambda}(x) = \f{8\lambda}{\,(\lambda + |x|^2)^2\,}$. This functional $H_{\lambda}$ has been a key tool of the study of the critical-mass 2D Keller--Segel system, since it has been proven in~\cite{Blanchet_Carlen_Carrillo12} that this quantity is also nonincreasing in time under the flow. Therefore, $H_\lambda$ is such that
\[
    H_{\lambda}[\rho(t)]   \le   H_{\lambda}[\rho_0] \quad \text{for all }t\ge0.
\]
The first global existence result in this setting was produced under the assumption $H_{\lambda}[\rho_0]<+\infty$.

The goal of this subsection is to convince the reader that this is just a particular case of our general approach (which, of course, also provides stronger estimates on the solution).

In order to carry out our strategy, we first define a suitable set $K$. 

\begin{proposition}\label{prop:compact_set_Hlambda}
Define 
\[
    K_{\lambda,C}  \coloneqq  \left\{\, \rho\in \mathcal{M}(\R^2)  \,:\, \rho(\R^2)= M_c, \, H_{\lambda}[\rho]\le C \right\}.
\]
Then:
\begin{enumerate}
    \item $K_{\lambda , C}$ contains no  Dirac masses.
    \item $K_{\lambda , C}$ is sequentially compact for the narrow convergence.
\end{enumerate}
\end{proposition}

\begin{proof}
We will prove at the same time that any $\rho\in K_{\lambda , C}$ has infinite second moment (which implies that $K_{\lambda , C}$ does not contain any Dirac masses) and that there is a uniform bound on the moments of order $r<2$, i.e. $\int |x|^r\diff\rho\leq C$ for every $\rho\in K_{\lambda , C}$ (which implies that $K_{\lambda , C}$ is tight).

In order to show this let us notice that we have that 
$$
    \rho\leq 2u_\lambda+2(\sqrt{\rho}-\sqrt{u_\lambda})^2 \; \mbox{ and } \; u_\lambda\leq 2\rho+2(\sqrt{\rho}-\sqrt{u_\lambda})^2.
$$ 
This implies that for any $r>0$ such that $\int (\sqrt\rho-\sqrt{u_\lambda})^2|x|^r<+\infty$, then the two conditions $\int |x|^r\diff\rho<+\infty$ and $\int |x|^r\diff u_\lambda<+\infty$ are equivalent, and that we can bound one quantity with the other up to constants.

Now, from the fact $\rho\in K_{\lambda , C}$ and the fact that $u_\lambda^{-1/2}$ behaves as $1+|x|^2$, we deduce that $\int (\sqrt\rho-\sqrt{u_\lambda})^2|x|^r$ can be bounded for any $r\in [1,2)$. Yet, using the fact that
%
%
%
$u_\lambda$ behaves as $|x|^{-4}$, we see that we have $\int |x|^r\diff u_\lambda<+\infty$ if and only if $r\in [0,2)$. This shows the desired claim.
 


In order to prove compactness of $K_{\lambda , C}$ we only need to prove that it is closed for the narrow convergence, i.e. it remains to prove that  $H_\lambda$ is lower-semicontinuous w.r.t. to the same convergence. In particular, this fact follows from~\cite{Blanchet_Carlen_Carrillo12}, observing that the integral expression defining $H_\lambda$ is convex in the density $\rho$ (of course $H_\lambda$ has to be properly defined for measures $\rho$ which are not absolutely continuous, using the recession function).
This concludes the proof. 
\end{proof}

Since the sets $K_{\lambda,C}$ are preserved by the evolution of the Keller--Segel equation, thanks to \Cref{prop:compact_set_Hlambda}  we can conclude that the result presented in Proposition~\ref{prop:linf_estimate_critical} is valid and we can repeat the proof of the previous subsection for the Li--Yau estimate in dimension 2, that is~\ref{2d_exact_mass}. Regarding the bound on the entropy, we will see in subsection \ref{sec:Estimates entropy} that this holds whenever we restrict to compact sets with no Dirac masses. Nevertheless, it has also been explicitly proven in~\cite{Blanchet_Carlen_Carrillo12} that a bound on $H_\lambda$ and on $\mathcal F$ imply bounds on the entropy. 

\subsubsection{Assuming only that \texorpdfstring{$\rho_{0}$}{rho0} is not a  Dirac in any dimension}

In fact, it is possible to avoid the assumption that $H_{\lambda}[\rho_{0}]<+\infty$ and to treat the problem in a unified way the multi-dimensional case.  

We focus here on the choice of the compact set $K$. In the rest of the section, we will often use the Wasserstein distance $W_2$. However, we do not need to impose the finiteness of the second moment $m_2[\rho]$. Only, in dimension $2$ we need to impose finiteness of the logarithmic moment, which is always required at least to define $u=\Gamma\ast\rho$. This means that the quantity $W_2(\mu,\nu)$ is defined via
$$
       W_2 (\mu , \nu) \coloneqq \inf_{\gamma \in \Gamma (\mu , \nu)} \left\{ \int_{\Rd \times \Rd} | x - y |^2 \diff \gamma (x, y) \right\}^{\frac{1}{2}},
$$
where $\Gamma (\mu , \nu)$ is the set of transport plans between $\mu$ and $\nu$, assumed to be of same mass,
    \begin{equation*}
        \Gamma (\mu , \nu ) = \left\{ \gamma \in \mathcal{M} (\Rd \times \Rd) : (\pi_x)_\# \gamma = \mu, (\pi_y)_\# \gamma = \nu \right\}.
    \end{equation*}
This infimum could a priory be infinite. Also, we use the Wasserstein distance on measures which have finite mass without being probability measures, which requires very small adaptations of the language.

We need to consider the following quantity: given a nonnegative finite measure $\rho$ we consider 
$$
D[\rho]=\inf_{x_0}W_2(\rho,M\delta_{x_0})
=\left(\inf_{x_0}\int|x-x_0|^2\diff\rho(x)\right)^{1/2},
\qquad M=\rho(\R^d).
$$
Let us point out that the above infimum takes the value $+\infty$ if $m_2[\rho]=+\infty$. If, instead, $\rho$ has finite second moment, then the infimum is achieved and we have a  minimum, realized by taking $x_0$ to be the barycenter of $\rho$. Moreover, we have $D[\rho]>0$ unless $\rho$ is a Dirac mass.

We now consider a new compact set, depending on $\rho_0$. We introduce this set and present some of its key properties.

\begin{proposition}\label{prop:compact_set_no_ Dirac}
Given a measure $\rho_0$ which is not a Dirac mass, define
\begin{equation}\label{def:set K}
K[\rho_0]
  \coloneqq 
 \left\{
\rho \in \mathcal{M}_{+}(\R^d)
   \Bigm| 
\rho(\R^d)=\rho_0(\R^d), 
W_{2} \left(\rho,\rho_0 \right) \le  \min\left\{\frac{D[\rho_0]}{2},1\right\}
 \right\}.
\end{equation}
Then, the set $K[\rho_0]$ is such that: 
\begin{enumerate}
\item $K[\rho_0]$ does \emph{not} contain  Dirac masses,
\item $K[\rho_0]$ is sequentially compact for the narrow convergence.
\end{enumerate}
\end{proposition}

\begin{remark}
We allow $K[\rho_0]$ to depend on the initial measure $\rho_0$  since the estimate $W_2^2(\rho,\rho_0)\le C$ typically only holds for relatively small times. Afterwards,  we construct a new set $K[\rho(t_1, \cdot)]$ after a short time $t_1$ to iterate the argument further. Note that the minimum is only added to take into account the fact that $D[\rho_0]$ could be infinite and the constant $1$ in this minimum is arbitrary. This addition is not necessary if the initial condition is assumed to have a bounded second moment. 
\end{remark}

\begin{proof}
The number $D[\rho_0]$ (strictly positive since $\rho_0$ is not a Dirac mass) is, by definition, the minimal $W_2$ distance from $\rho_0$ to Dirac masses, so of course no Dirac mass can be at distance smaller that $\frac{D[\rho_0]}{2}<D[\rho_0]$ to $\rho_0$.

It remains to prove that $K[\rho_0]$ is sequentially compact. In particular, we need to prove that it is a tight set. This would be easy if $m_2[\rho_0]<+\infty$, since in this case $\rho\in K$ would imply a bound on $m_2[\rho]$. Dealing with the general case is slightly harder. We let $\eta\coloneqq \min\left\{\frac{D[\rho_0]}{2},1\right\}$.

 Let $\rho\in K[\rho_0]$. We set $\rho_0(\mathbb{R}^d) \coloneqq M$ and $\rho_0$ is a Radon measure. For any $\varepsilon>0$, there exists a compact set $K \subset \mathbb{R}^d$ such that
$$
  \rho_0(K)  >  M -  \frac{\varepsilon}{2}.
$$
We aim to show that, for a suitable $\delta>0$, the measure $\rho$ places at least $(M-\varepsilon)$ of its mass in the $\delta$-enlargement
$$
  K^\delta  \coloneqq  \left\{\,y \in \mathbb{R}^d : \mathrm{dist}(y,K)\le \delta \right\}.
$$
Since $K$ is compact, $K^\delta$ is also compact.
By the definition of the Wasserstein distance, there exists a coupling $\gamma$ with marginals $\rho_0$ and $\rho$ such that
$$
    \int_{\mathbb{R}^d\times \mathbb{R}^d} |x-y|^2\, d\gamma(x,y) \le \eta^2 .
$$
For $K\subset \mathbb{R}^d$ compact and $\delta>0$, every $(x,y)\in K\times (K^\delta)^c$ satisfies $|x-y|\ge \delta$, hence
$$
\delta^2\,\gamma(K\times (K^\delta)^c)
   \le 
   \int_{K\times (K^\delta)^c} |x-y|^2\, d\gamma
   \le 
   \eta^2.
$$
Therefore, it follows that
$$
\gamma(K\times (K^\delta)^c)\le \frac{\eta^2}{\delta^2}.
$$
If we use the marginal conditions we also get
$$
\gamma(K\times K^\delta)
  =\rho_0(K)- \gamma(K\times (K^\delta)^c)
  \ge \rho_0(K)-\frac{\eta^2}{\delta^2}.
$$
Let us choose $\delta$ so that $\eta^2/\delta^2\le \varepsilon/2$. Then, we obtain that
$$
\gamma(K\times K^\delta)\ge M-\varepsilon.
$$
Since $\gamma(K\times K^\delta)\le \gamma(\mathbb{R}^d\times K^\delta)=\rho(K^\delta)$, $\rho(K^\delta)\ge M-\varepsilon.
$
Thereby, we obtain the tightness of the set. In order to conclude, we need to show that $K$ is closed for the narrow convergence. This follows as a consequence of the lower semicontinuity of the transport costs, see \cite[Proposition 7.4]{Santambrogio15}.
\end{proof}

With $K[\rho_0]$ defined in \eqref{def:set K}, we still have to show that solutions $\rho_t$ of the Keller--Segel system remain in $K[\rho_0]$ for a short interval $[0,T_0]$. Once that is established, we can define a second set $K[\rho_{T_0}]$ around the measure $\rho_{T_0}$ at time $T_0$, and to keep iterating. This procedure is repeated over successive intervals $[kT_0,(k+1)T_0]$ to obtain global-in-time control.

\begin{proposition}\label{prop:short_time_stay_in_K}
Let $\rho_t$ be a solution of the Keller--Segel system with critical mass $M_c$ and initial condition $\rho_0$ with $\mathcal F[\rho_0]<+\infty$ (and finite logarithmic moment if $d=2$).  Then:
\begin{enumerate}
\item There is a time $T_0>0$, depending on $\eta$, such that for all $t\in[0,T_0]$ we have that $\rho_t \in K[\rho_0]$.
\item At $t=T_0$, the measure $\rho_{T_0}$ satisfies $D[\rho_{T_0}]\geq D[\rho_0]$. 
Hence we can reset and iterate on $[k T_0,(k+1) T_0]$ successively for every $k \in \mathbb{Z}_+$.
\end{enumerate}
\end{proposition}

\begin{proof}
The Keller--Segel system is a gradient flow with respect to the $W_2$ distance. As a consequence, and up to considering solutions constructed by the JKO scheme (see e.g. \cite{Blanchet_Calvez_Carrillo08}), the solutions of the Keller--Segel system satisfies 
\begin{equation}\label{eq:estim_w2}
W_{2}^2(\rho_t, \rho_{0})\le T\mathcal{F}[\rho_0].
\end{equation}
This inequality is obtained by looking at the dissipation of $\mathcal F$ along solutions (or considering the EDE formulation of gradient flows \cite[Chapter 14]{Ambrosio_Brue_Semola21}). This proves the first item as soon as one takes $T_0\mathcal F(\rho_0)\leq \frac{1}{2}\min\left\{\frac{D[\rho_0]}{2},1\right\} $.

It remains to prove the second item. First let us assume that the second moment of the initial condition is infinite. Then the second moment of the solution stays infinite, this can be seen as a consequence of~\eqref{eq:estim_w2} for instance. Moreover we have $\mathcal{F}[\rho_{T_0}]\le \mathcal{F}[\rho_0]$. Therefore the proof is concluded in a similar way in this case. It remains to treat the case where the second moment of the initial condition is assumed to be finite. Let us assume this and let $\rho_t$ be a solution of the Keller--Segel system. We recall that for $d=2$ we have (see, for instance, \cite{Blanchet_Dolbeault_Perthame06}) that 
$$
    \f{\diff}{\diff t}\int_{\R^2}|x|^2\rho_t =4M \left(1-\f{M}{8\pi} \right) 
$$
and, in the case of the critical mass $M=M_c=8\pi$, we obtain that $\int_{\R^2}|x|^2\rho_t$ is constant in time. Of course, the same holds, just by translating the reference system, for $\int_{\R^2}|x-x_0|^2\rho_t$. This shows that, for $d=2$, the quantity $D[\rho_t]$ is constant in time. 

In higher dimension the computation is different and we have (see~\cite{Blanchet_Carrillo_Laurencot09}): 
$$
\f{\diff}{\diff t}\int_{\R^d}|x-x_0|^2 \rho_t = 2(d-2)\mathcal{F}[\rho]\ge 0.
$$
This proves, as soon as $M\leq M_c$, that $t\mapsto D[\rho_t]$ is actually non-decreasing and concludes the proof. \qedhere
\end{proof}

We have proved that the solutions of the Keller-Segel system stay inside a series of compact $K[\rho_{kT_0}]$ for $t\in[kT_0, (k+1)T_0)$. To conclude the argument as in the subscritical case, we need to show that on these compacts, we can bound the entropy ($d=2$) and the $L^m$ norm ($d\ge 2$) of the solutions.

\subsubsection{Estimates on the entropy in dimension 2}\label{sec:Estimates entropy}

All the analysis of this section will be performed on a compact set of measures not including Dirac masses. We first need a fundamental property of those sets.

\begin{proposition}
    Let $K\subset\mathcal M_+(\R^d)$ be a set of measures, all with the same mass $M$, which is compact for the narrow convergence and does not contain any Dirac mass. Then, there exist $\ve,R>0$ such that for any $\rho\in K$ and any $x\in\R^d$ we have $\rho(B_R(x))<M-\ve$.
\end{proposition}
\begin{proof}
    The proof is done by contradiction. If the statement is false there exists a sequence $\rho_n\in K$ and a sequence of points $x_n$ such that  $\rho_n(B_{1/n}(x_n))\ge M - 1/n$. By tightness of $K$ the sequence $x_n$ should be bounded, since the measures in $K$ give small mass out of a common big ball $B_R$, so that we have $x_n\in B_{R+1/n}$. We can then suppose, up to taking a subsequence, that we have $x_n\to x^*$. Moreover, by compactness of $K$, we can also assume, up to taking a subsequence, that we have $\rho_{n}\rightharpoonup\rho^*$ narrowly, for some $\rho^*\in K$. Then, for any $r>0$, we have $\rho^*(B(x^*,r)\geq \limsup_n \rho_n(B(x^*,r))\geq \limsup_n \rho_n(B(x_n,1/n))=M$, since we have $B(x_n,1/n)\subset B(x^*,r)$ for large $n$. This proves $\rho^*(\{x^*\})\geq M$, and hence $\rho^*=M\delta_{x^*}$, which is a contradiction because $K$ does not contain Dirac masses.
\end{proof}
We then prove a version of the log-HLS inequality. The same statement contains the cases of subcritical and critical masses.

\begin{theorem}\label{thm:new_proof_hls}
Assume that $\rho\ge 0$ is an absolutely continuous measure with finite logarithmic moment such that $\mathcal{F}[\rho]<+\infty$. Take $V\ge0$ such that $K_V=\int_{\R^2}e^{-V}\diff x<+\infty$ and
$\int_{\R^2}\rho V\diff x<+\infty$.
\begin{enumerate}
    \item If the mass $M$ of $\rho$ satisfies $M<8\pi$, then there exists a constant $C=C(M,K_V)$ such that
    $$
        \int_{\R^2}\rho\log\rho\diff x\leq C \left(\mathcal{F}[\rho]+\int V\diff\rho \right).
    $$
    \item If $M=8\pi$, since $\rho$ is not a  Dirac mass,  there exist $\varepsilon>0$, and $R>0$ such that for all $x\in\R^2$ we have $\rho(
B_{R}(x))\leq 8\pi-\varepsilon$; take a constant $\kappa$ such that  $2 - \f{\varepsilon}{4\pi}<\kappa<2$. Then, there exists a constant $C=C(\kappa)$  such that
$$
\left(1-\frac{\kappa}{2}\right)\int_{\R^2}\rho\log\rho\diff x +  C(\kappa)(1+\log R)\le \mathcal{F}[\rho].
$$
\end{enumerate}
\end{theorem}


\begin{proof}
We divide the proof in several steps.
\begin{enumeratesteps}
\step[Setting the inequality]
We first exploit the inequality 
$$
ab\le a\log a + e^{b-1}, 
$$
which holds for all $a,b \in \R$ with $a\ge 0$ and it follows directly from the Legendre's duality. Therefore, if we choose $a=\rho$ and $b=\f{1}{\kappa}u-V(x)$ we obtain
\begin{equation}\label{eq:Inequality for entropy}
    \int_{\R^2}\rho u\diff x \le \kappa\int_{\R^2}\rho\log\rho +\kappa e^{-1}\int_{\R^2}e^{\f{u}{\kappa}-V}\diff x+\kappa\int V\diff\rho. 
\end{equation}
In order to estimate the second term of the right-hand side, we observe that, in both the cases of our statement, for every $x$ there exists a value $R(x)\in [R,+\infty]$ such that $\rho(
B_{R}(x))= 8\pi-\varepsilon$ and that it satisfies $R(x)\ge R$. In the first case, we choose $\ve=8\pi-M>0$ and we have $R(x)=+\infty$ for every $x$, while in the second case we have $R(x)<+\infty$ for every $x$.

Let us take into account the assumptions on $R$ from the statement of the proposition. Let us also take advantage of the Jensen's inequality. Then, we can compute in order to obtain that  
\begin{align*}
e^{\f{u(x)}{\kappa}}
&=\exp\left(
\int_{\R^2}\f{1}{2\pi\kappa}
\log\left(\f{1}{|x-y|}\right)\rho(y)\diff y
\right)\\
&=\exp\left(
\int_{\R^2}
\log\left(\f{1}{|x-y|^{\f{1}{2\pi\kappa}}}\right)
\rho(y)\diff y
\right)\\
&=\exp\left(
\begin{aligned}
&\int_{B_{R(x)}(x)}
\log\left(\f{1}{|x-y|^{\f{1}{2\pi\kappa}}}\right)
\rho(y)\diff y\\
&+\int_{(B_{R(x)}(x))^c}
\log\left(\f{1}{|x-y|^{\f{1}{2\pi\kappa}}}\right)
\rho(y)\diff y
\end{aligned}
\right)\\
&=\underbrace{
\exp\left(
\int_{B_{R(x)}(x)}
\log\left(\f{1}{|x-y|^{\f{1}{2\pi\kappa}}}\right)
\rho(y)\diff y
\right)}_{\eqqcolon A}\\
&\quad\times
\underbrace{
\exp\left(
\int_{(B_{R(x)}(x))^c}
\log\left(\f{1}{|x-y|^{\f{1}{2\pi\kappa}}}\right)
\rho(y)\diff y
\right)}_{\eqqcolon B}.
\end{align*}
In the following we bound the exponential factors $A$ and $B$. We observe that in the first case, i.e. $R(x) = + \infty$, we have $B_{R(x)}^c=\emptyset$ and $B=1$.

\step[A bound on $A$ and the case of subcritical mass]
Let us start with the bound on $A$. Let us recall that $\rho(B_{R(x)}(x))=8\pi-\varepsilon$. Therefore, we can use  Jensen's inequality in order to recover that 
\begin{align*}
A&= \exp\left(\int_{B_{R(x)}(x)}\log\left(\f{1}{|x-y|^{\f{1}{2\pi\kappa}}}\right)\rho(y)\diff y\right)\\
&= \exp\left(\int_{B_{R(x)}(x)}\log\left(\f{1}{|x-y|^{\f{\rho \left( B_{R(x)}(x) \right)}{2\pi\kappa}}}\right)\f{\rho(y)}{\rho(B_{R(x)}(x))}\diff y\right)\\
&\le \f{1}{8\pi-\varepsilon}\int_{B_{R(x)}(x)}\f{1}{|x-y|^{\f{8\pi-\varepsilon}{2\pi\kappa}}}\rho(y)\diff y.
\end{align*}
We choose $\kappa$ such that $2 - \f{\varepsilon}{4\pi}<\kappa<2$ so that we have $s \coloneqq \f{8\pi-\varepsilon}{2\pi\kappa} < 2$.

Let us now concentrate on the first case, i.e. subcritical mass $M<8\pi$. We then have to estimate
$$
    \int_{\R^2}e^{\f{u}{\kappa}-V}\diff x\leq \f{1}{M}\int_{\R^2}e^{-V(x)}\diff x\int_{\R^2}\f{1}{|x-y|^s}\diff\rho(y).
$$
If we change the order of the integration in the right-hand side  we obtain
$$
\int_{\R^2}e^{\f{u}{\kappa}-V}\diff x\leq \f{1}{M}\int_{\R^2}\diff\rho(y)\int_{\R^2}e^{-V(x)}\f{1}{|x-y|^s}\diff x.
$$
Let us also recall that we have that
\begin{eqnarray*}
\int_{\R^2}e^{-V(x)}\f{1}{|x-y|^s}\diff x&=&\int_{B_1(y)}e^{-V(x)}\f{1}{|x-y|^s}\diff x+\int_{B_1(y)^c}e^{-V(x)}\f{1}{|x-y|^s}\diff x\\
&\leq&\int_{B_1(y)}\f{1}{|x-y|^s}\diff x
   +\int_{\R^2}e^{-V(x)}\diff x.
\end{eqnarray*}
where the two integrals can both be bounded by a constant only depending on $s$ (hence on $\kappa$) and on $K_V$. 
Therefore, in the first part of the statement ($M<8\pi$ and $R(x) = + \infty$), we obtain
$$
\int_{\R^2}\rho u\diff x \le \kappa\int_{\R^2}\rho\log\rho +\kappa\int V\diff\rho+C, 
$$
and the statement follows by adding $(2-\kappa)\int \rho\log\rho$ to both sides,
using $\mathcal F[\rho]=\int\rho\log\rho-\frac12\int\rho u$.
\step[Critical mass, a bound on $B$ and a new bound on $A$]
We now move on to the second case, i.e. $M=8\pi$, where we also have to take into account the second factor. In this case we choose $V=0$. This requires a different estimate for both $A$ and $B$.

We look at $B$ and we note that the mass on the complementary of $B_{R(x)}(x)$ is $\varepsilon$. Thus, 
\begin{align*}
    B \le  \exp\left( \int_{(B_{R(x)}(x))^c}\log\left(\f{1}{R^{\f{1}{2\pi\kappa}}}\right)\rho(y)\diff y\right) \le \exp\left( \log\left(\f{1}{R(x)^{\f{1}{2\pi\kappa}}}\right) \varepsilon\right)  = \frac{1}{R(x)^{\frac{\varepsilon}{2\pi\kappa}}}.
\end{align*}
%
Observe that we have $\frac{\varepsilon}{2\pi\kappa}=\frac 4\kappa-s$. Furthermore, up to changing the value of $\varepsilon$, we can assume $\kappa>0$. Moreover, since $R(x)\ge R$ we observe that for $y\in B_{R(x)}(x)$, we have $R(x)\ge \max\left(R,|x-y|\right)$. We can now compute the second term on the right-hand side of \eqref{eq:Inequality for entropy} and we obtain 
\begin{align*}
\int_{\R^2}e^{\f{u}{\kappa}}\diff x 
&\lesssim \int_{\R^2} \frac{1}{R(x)^{\frac{4}{\kappa}-s}}\int_{B_{R(x)}(x)}\f{1}{|x-y|^s}\rho(y)\diff y\diff x\\
&=\int_{\R^2}\rho(y)\diff y\int_{\R^2}\frac{1}{R(x)^{\frac{4}{\kappa}-s}}\f{1}{|x-y|^s}\rho(y)\ind_{|x-y|\le R(x)}\diff x \\
&\leq \int_{\R^2}\rho(y)\diff y\int_{\R^2} \frac{1}{\max\left(R^{\frac{4}{\kappa}-s}, |x-y|^{\frac{4}{\kappa}-s}\right)}\f{1}{|x-y|^s}\diff x\\
&=\int_{\R^2}\rho(y)\diff y\int_{B_{R}(y)}\frac{1}{R^{\frac{4}{\kappa}-s}}\f{1}{|x-y|^s}\rho(y)\diff x  + \int_{\R^2}\rho(y)\diff\int_{ B_{R}(y)^c}\f{1}{|x-y|^{\frac{4}{\kappa}}}\diff x \\
&= C(\kappa) R^{2-\frac{4}{\kappa}}.
\end{align*}
In the last inequality, we used the fact that the total mass of $\rho$ is $8\pi$ and we explicitly computed the integrals in $x$.
%
%
%
Hence, the conclusion of this argument is that we have shown that for $2 - \f{\varepsilon}{4\pi}<\kappa<2$ there exists a  constant $C(\kappa)$ such that
$$
\int_{\R^2}\rho u\diff x \le \kappa\int_{\R^2}\rho\log\rho +C(\kappa)R^{2-\frac{4}{\kappa}}. 
$$
We can optimize this result. Let us consider the scaling $\rho_{\lambda}(x)=\lambda^{2}\rho(\lambda x)$ which preserves the mass. Then, a similar argument can be performed for $\rho_{\lambda}$ and we obtain that
$$
    \int_{\R^2}\rho_{\lambda} u_{\lambda}\diff x \le \kappa\int_{\R^2}\rho_\lambda\log\rho_\lambda +C(\kappa)\left(\frac{R}{\lambda}\right)^{2-\frac{4}{\kappa}}. 
$$
We notice that we have
\begin{equation*}
    \int_{\R^2}\rho_{\lambda} u_{\lambda} = \int_{\R^2}\rho u + 32\pi \log\lambda\quad \mbox{ and }\quad\int_{\R^2}\rho_{\lambda}\log\rho_{\lambda}= \int_{\R^2}\rho\log\rho + 16\pi \log\lambda.  
\end{equation*}
Therefore, we have
$$
    \int_{\R^2}\rho u\diff x \le \kappa\int_{\R^2}\rho\log\rho +C(\kappa)\left(\frac{R}{\lambda}\right)^{2-\frac{4}{\kappa}}+(\kappa-2) 16\pi\log\lambda. 
$$
We can optimize  in $\lambda$. We choose $\lambda=R$ and this yields that there exists a new constant $C(\kappa)$ such that 
$$
    \int_{\R^2}\rho u\diff x \le \kappa\int_{\R^2}\rho\log\rho +C(\kappa)(1+\log R)
$$
Once more, we use the definition of $\mathcal F$ and we obtain the desired result. \qedhere
\end{enumeratesteps}
\end{proof}

\subsubsection{Estimates on the \texorpdfstring{$L^m$}{Lm} norm for \texorpdfstring{$d>2$}{d>2}}\label{sec:Estimates Lm}

We now turn to the higher-dimensional ($d>2$) critical Keller--Segel system. Again, we rely on a family of compact sets, as in the two-dimensional case, and seek a replacement for the $L \log L$ bound — i.e., an $\|\rho\|_{L^m}$-bound, where $m = 2 - \frac{2}{d}$. In order to obtain such a control in the critical regime, we exploit a recent stability result for the Lane--Emden inequality due to Carlen, Lewin, Lieb, and Seiringer~\cite{Carlen_Lewin_Lieb_Seiringer24}, which provides an HLS-type estimate even at critical mass.

\medskip

\noindent

\begin{theorem}[Stability of the Lane--Emden inequality {\cite{Carlen_Lewin_Lieb_Seiringer24}}]\label{thm:stability_lane_emden}
Let $m = 2 - \frac{2}{d}$ and $M_c$ be the critical mass in dimension $d>2$. There exists a constant $C(d) \ge 0$ such that for all nonnegative $\rho\in L^m(\R^d)$ with $\int \rho = M_c$, one has
\[
\mathcal{F}[\rho]
  = 
\f{1}{\,m-1\,}\,\int_{\R^d}\rho^m \,\diff x 
  -  
\frac12 \int_{\R^d}\rho\,u \,\diff x
  \ge 
C\,\inf_{\ell}\int_{\R^d}\left(\rho^{m/2} - \ell^{m/2}\right)^2\,\diff x,
\]
where the infimum $\ell$ is taken over all the minimizes of the Lane--Emden-type inequality~\eqref{eq:HLS} with mass $M_c$.  
\end{theorem}

As proved for instance in~\cite[Proposition 3.5]{Blanchet_Carrillo_Laurencot09}, the minimizers of the Lane-Emden inequality are given by the two parameters family 
$$
\ell_{\lambda,x_0}(x)= \lambda^d \bar{\ell}(\lambda(x-x_0))
$$
where $\bar{\ell}$ satisfies for some $R>0$ the problem
$$
\begin{cases}
&\f{m}{m-1}\Delta \bar{\ell}^{m-1} +\bar{\ell}=0 \quad  \text{in $B_{R}(0)$}, \\
&\bar{\ell}=0 \quad  \text{on $\partial B_{R}(0)$}.
\end{cases}
$$
Up to changing $\lambda$ we assume $\|\bar{\ell}\|_{L^m}=1$.

Here, the key step is that the free energy $\mathcal{F}[\rho]$ controls some notion of distance to the Lane--Emden profile $\ell$. Combining this stability inequality and the fact that $\rho$ is not a  Dirac yields the bound we need.

\begin{theorem}\label{thm:new_proof_hls_multid}
Let $d>2$ and $m = 2 - \frac{2}{d}$.  Suppose $K\subset \mathcal{M}(\R^d)$ is a set of nonnegative measures with total mass $M_c$, which is compact under narrow convergence and does not contain any  Dirac masses. Then, for any constant $C_1>0$, there is a constant $C=C(K)>0$ such that
\[
    \sup\left\{\|\rho\|_{L^m} : \rho\in K, \, \mathcal{F}[\rho]\le C_1\right\}  \le  C.
\]
\end{theorem}

\begin{proof}
Argue by contradiction.  If no such $C$ existed, we could build a sequence $\rho_n\in K$ with
\[
\|\rho_n\|_{L^m}  \to  +\infty,
\quad
\mathcal{F}[\rho_n] \le  C_1.
\]
By Theorem~\ref{thm:stability_lane_emden}, for each $\rho_n$, there exists a Lane--Emden optimizer $\ell_n$ such that
\[
\int_{\R^d}\left|\rho_n^{\,m/2} - \ell_n^{\,m/2}\right|^2 \le C',
\]
for some $C'>0$ independent of $n$. 

We know each $\ell_n$ is of the form $\ell_{\lambda_n,x_n}(x)$ (i.e.\ a scaling and translation of a fixed Lane--Emden solution $\bar{\ell}$).  If $\{\lambda_n\}$ were bounded, then $\|\ell_n\|_{L^m}$ would also be bounded. That would imply that $\|\rho_n\|_{L^{m}}$ is bounded. We then use the elementary inequality 
$$
    a^p\le 2(a^{p/2}-b^{p/2})^2 + 3 b^p,
$$
that we apply to $a=\|\rho_n\|_{L^{m}}$ and $b=\| \ell_n\|_{L^{m}}$.  Hence we can consider $\lambda_n\to\infty$. Set $\tilde{\rho}_n= T_n\sharp \rho_n$ where $T_n(x)=\lambda_n(x-x_n)$. Then, if we perform a change of variables in the inequality
$$
    \int_{\R^d}|\rho_n^{m/2}-\ell_n^{m/2}|^2\le C, 
$$
we obtain
$$
    \int_{\R^d}|\tilde{\rho}_n^{m/2}-\bar{\ell}^{m/2}|^2\to  0 
$$
as $n\to+\infty$ and $\lambda_n^{d(1-m)}\to 0$. Observe the functional 
$$
    J[\rho]\coloneqq \int_{\R^d}|\rho^{m/2}-\bar{\ell}^{m/2}|^2\diff x = \int_{\R^d}\rho^m- 2\int_{\R^d}\rho^{m/2}\bar{\ell}^{m/2}+\int_{\R^d}\bar{\ell}^{m} 
$$
is convex with respect to $\rho$ as $m>1$ and $m/2 < 1$, so in particular it is lower semi continuous with respect to the weak convergence in $L^m(\R^d)$. Since $$
\int_{\R^d}|\tilde{\rho}_n^{m/2}-\bar{\ell}^{m/2}|^2\to  0 
$$
as $n\to+\infty$ we deduce that up to a subsequence not relabeled, $\tilde{\rho}_n\rightharpoonup\tilde{\rho}$ weakly in $L^m$ and by lower semi continuity:
$J(\tilde{\rho})=0$. Therefore $\tilde{\rho}=\bar{\ell}$. Since the mass of $\bar{\ell}$ is $M_c$, we deduce that in fact $\tilde{\rho}_n$ is tight and narrowly converges to $\bar{\ell}$. Moreover, we know that the set $\{\rho_n\}_n$ is tight since the sequence belongs to $K$. From Lemma~\ref{lem:contradiction_ Dirac_mass} we conclude that, up to a subsequence, $\rho_{n}$ converges narrowly  to a  Dirac mass. This is in contradiction with the fact that $K$ is a compact set for the narrow convergence with no  Dirac masses. 
Combining all steps proves that the solution remains in a series of compact sets $K[\rho_{kT_0}]$, each prevent Dirac blow-up. From this combination of steps we also prove a uniform $L\log L$ bound in dimension $2$ and a uniform $L^m$ bound when $d>2$. Similarly than in the subcritical mass case, we obtain a global Li--Yau estimate for $\rho$ at critical mass. 
\end{proof}

This completes the proof of Theorems~\ref{thm:Li--Yau} and~\ref{thm: Aronson--B\'enilan} in the critical mass case.

\section{Justification of the Li--Yau and  Aronson--B\'enilan estimates}\label{sec:justification}   
In the formal derivation of the Li--Yau and Aronson--B\'enilan estimates for the Keller--Segel system, we used a comparison principle for parabolic equations. This requires some attention to be carried out. Let us explain how this is usually done. Suppose you have a quantity $f=f(t,x)$ (here $f=\Delta v$, for our goals) and you know that $f$ satisfies a parabolic differential inequality of the form
$$\partial_t f\geq A:D^2f+w\cdot \nabla f+\omega(f),$$
where $A$ is a matrix field, depending on the variables $(t,x)$, which is semi-positive definite (often the second-order term $A:D^2f$ takes the form $a(t,x)\Delta f$) and $w$ is a vector field. The goal is then to compare $f$ to the solution of the ODE $g'=\omega(g)$, with an initial datum $g(0)$ such that $f\geq g(0)$, and obtain $f(t,x)\geq g(t)$ everywhere. 

This can be done by an abstract parabolic comparison result (we refer for instance to \cite[Theorem 9.7]{Lieberman96}), which exploits the fact that $f$ and $g$ are respectively a super-solution and a sub-solution of the same equation (the second one being independent of $x$), and then the order of the initial datum is preserved. Yet, this procedure requires some assumptions, and in particular it is required that the vector field $w$ should be bounded, or satisfy some bounds, roughly the same which are needed in order to guarantee global existence of solutions to the ODE $y'=w(t,y)$. This is easy to understand: the comparison principle implies uniqueness of the solution of the PDE with given initial datum, and if we consider the first-order case, i.e. $A=0$, it is clear that the PDE $\partial_t f=w\cdot\nabla f$ only has uniqueness on the set of points $(t,x)$ in space-time which can be reached by trajectories of the ODE $y'=w(t,y)$ starting from $t=0$. Other boundedness assumptions (on $f$, on $A$\dots) are also usually required in order to carry on the comparison principle, but the one on $w$ is quite striking.

Another possibility is to redo the proof by hands considering the quantity $F(t) \coloneqq \min_x f(t,x)$ and looking at the ODE satisfied by $F$. This works under the following assumptions: assume that $f$ is smooth and that, at least locally in time, there exists a compact set $K$ such that for every $t$ the function $f$ attains its global minimum on $K$ (locally in time means that every $t_0$ is contained in an open interval $I\ni t_0$ such that this holds for $t\in I$, the set $K$ possibly depending on $I$). This strategy is explicitly used, for instance, in \cite{David_Santambrogio23}. It works as follows: first we write $F(t)=\min_{x\in K} f(t,x)$, then we observe that $F$ is locally Lipschitz continuous on $I$ as an infimum of Lipschitz functions of $t$, indexed with $x$, which are equi-Lipschitz because of the smoothness of $f$ and the compactness of $K$. Hence, $F$ is differentiable for a.e. $t$ and the envelop theorem states that for every $t_0$ where $F$ is differentiable, if $x_0\in K$ is a point realizing the minimal value of $f(t,\cdot)$, then $F'(t_0)=\partial_t f(t_0,x_0)$. Then, it is possible to use the minimality of $x_0$ in order to get rid of the terms $A(t_0,x_0):D^2f(t_0,x_0)$ (which is nonnegative) and $w(t_0,x_0)\cdot \nabla f(t_0,x_0)$ (which vanishes, independently of any bound on $w$) and to obtain $F'\geq \omega(F)$.

Using the first strategy (i.e. a blackbox comparison principle) has an extra difficulty in our case: we do not really have an inequality only involving $f$ and its derivatives, but some terms (that we called $Q(u)$) explictly involve $\min f$. It means that we have an inequality of the form $\partial_t f\geq A:D^2f+w\cdot \nabla f+\omega(f,\min f),$ and that we want to compare to the solutions of the ODE $g'=\omega(g,g)$.  Now that we have listed the main technical obstructions, we want to explain what is usually done in the proof of the Aronson--B\'enilan estimate for the porous medium equation, as done in \cite{Vazquez06}: the initial datum is usually approximated with initial data which are bounded from below by a strictly positive constant; this bound is preserved in the evolution and makes the porous medium equation uniformly parabolic; this allows to obtain boundedness of all the derivatives of the solution, and hence of the vector field $w$ appearing in the parabolic inequality satisfied by $f=\Delta p$ for $p = v + u$, and of all the other relevant terms. This strategy is impossible to perform in our case, as for us the finiteness of the mass of $\rho_0$ on the whole space is crucial (not only it has to be finite, but also subcritical or at most critical). 

For all these reasons, we will explain here in this section how to obtain an approximation argument which guarantees that $\Delta v$ is smooth and attains its minimum on a compact set, without increasing the mass of the initial datum.

First, we approximate the initial datum $\rho_0$ with bounded and smooth initial data. From bounded $\rho_0$ it is well-known, see for instance \cite{Carrillo_Santambrogio18}, that at least for short time an $L^\infty$ solution to the Keller--Segel system exists. The existence time could depend on the initial $L^\infty$ norm but this is not an issue: indeed, our goal is to obtain bounds on $\delta(t)$, and if we manage to do so, we also obtain bounds on $\|\rho_t\|_{L^\infty}$, which means that the $L^\infty$ norm does not deteriorate and that the interval on which we have such a boundedness property can be extended.
Next, we need to guarantee all the other bounds hold in order to perform our analysis. More precisely we need the following key points:
\begin{itemize}
\item At least for $d>2$, we need strict positivity of the solution, in order to have local uniform parabolicity and ensure smoothness. We will indeed guarantee more than this, i.e. uniform upper and lower bounds on $\rho$ of the form $\rho(t,x)\approx (1+|x|)^{-\beta}$ for some exponent $\beta>d$ (this condition is necessary for the finiteness of the mass).
\item We consider $p$ the pressure of the diffusive term, i.e. $p = v + u$ which corresponds with $p = \log \rho$ if $d = 2$ or $p = \frac{m}{m-1} \rho^{m-1}$ if $d > 2$. We need $\Delta p(t,x)\to 0$ as $|x|\to\infty$ uniformly in $t$, and here as well we will prove better than this, and more precisely we will prove $|\nabla p(t,x)|\leq C(1+|x|)^{-(d-1)}$ and $|D^2 p(t,x)|\leq C(1+|x|)^{-\beta_1}$ for another exponent $\beta_1$ very close to $d$ 
\end{itemize}

The main statement of this section is the following.

\begin{proposition}\label{prop:lipschitz-porous_main}
    Any initial datum with finite mass $M$ and finite free energy can be approximated with initial data $\rho_0$ such that:
    \begin{itemize}
        \item The mass of $\rho_0$ is $M$ and its free energy approximates the one of the target initial datum.
        \item The evolution $\rho$ stemming from $\rho_0$ is stricly positive and $C^\infty$ on any interval $[0,T]$ on which it is bounded, and it satisfies uniformly in $t$
        \begin{itemize}
            \item $\rho_t\approx (1+|x|)^{-\beta}$, where $\beta$ is an exponent as close as we want to $d$,
            \item $|\nabla p_t|\lesssim (1+|x|)^{-(d-1)}$,
            \item $|D^2 p_t|\lesssim (1+|x|)^{-\beta_1}$, where $\beta_1$ is another exponent as close as we want to $d$.
             \end{itemize}
    \item For any $t\in [0,T]$ the function $\Delta v[\rho_t]$ attains its minimum over $\R^d$ and the minimum points are contained in the same compact set $K\subset\R^d$.
    \end{itemize}
\end{proposition}

\begin{proof}
It is absolutely standard to approximate the desired initial datum with bounded and smooth initial data $\rho_0$, preserving its mass and approximating its free energy. Moreover, we can impose strict positivity and the behavior at infinity of $\rho_0$. We choose this behavior to be a polynomial decay, that we take of the form $\rho_0\approx (1+|x|)^{-\beta}$. Because of integrability, we need $\beta>d$, but we choose $\beta$ very close to $d$. 

We then assume extra conditions on the initial datum $\rho_0$, but express them in terms of the pressure $p_0$. 
More precisely, for $d>2$, calling $p_0$ the pressure associated with $\rho_0$ (i.e. $p_0=\frac{m}{m-1}\rho_0^{1-2/d}$), we have $p_0\approx (1+|x|)^{-\sigma}$ where $\sigma=\frac{d-2}{d}\beta>d-2 $. We also impose other conditions, compatible with this one: we require $|\nabla p_0(x)|\leq C(1+|x|)^{-\sigma-1}$ and $|D^2 p_0(x)|\leq C(1+|x|)^{-\sigma-2}$. 

For $d=2$ the expressions are slightly different: we also impose a polynomial decay on $\rho_0$, of the form $\rho_0\approx (1+|x|)^{-(2+\varepsilon)}$, but now we have $p_0=\log\rho_0$ and we require $|\nabla p_0(x)|\leq C(1+|x|)^{-1}$, $|D^2 p_0(x)|\leq C(1+|x|)^{-2}$. 

Proposition~\ref{prop:lipschitz-porous_main} and  Proposition~\ref{prop:upper_lower_bound_pressure}  prove that, on any time interval $[0,T]$ where $\rho$ is bounded, than both $\rho(t,x)$ and $D^2p(t,x)$ tend uniformly to $0$ as $|x|\to\infty$. In particular, their minimal value is attained on a same bounded compact set. This allows to apply the strategy described above and in particular the maximum principle so as to obtain a bound on $\delta(t)$. In particular, the $L^\infty$ norm of $\rho$ does not blow up and the same argument can be extended in time.
\end{proof}

\begin{proposition}[Upper and lower bound on the pressure]\label{prop:upper_lower_bound_pressure}
Let $\rho_{0}$ be an initial condition as in \Cref{prop:lipschitz-porous_main} 
for~\eqref{eq:KS}.  Let $\rho:[0,T]\times\R^d\to(0,\infty)$ be a classical
solution of the Keller--Segel system. Assume that $\rho_0$ satisfies 
$$
C_1(1+|x|)^{-\beta}\le \rho_0(x)\le C_2
(1+|x|)^{-\beta}
$$
for two positive constants $C_1,C_2$. Then there exist two constants $\tilde C_1,\tilde C_2$, depending on $C_1$, $C_2$, $d$, $\beta$, $T$ and $\|\rho\|_{L^\infty}$, such that $\rho$ satisfies for all $t,x$
$$
\tilde C_1(1+|x|)^{-\beta}\le \rho(t,x)\le \tilde C_2(1+|x|)^{-\beta}.
$$
\end{proposition}

\begin{proof}
We start from the case $d>2$ and we use the equation satisfied by $p$:
\begin{equation}\label{eq:eq_for_pressure}
\partial_t p = |\nabla p|^2 + (m-1)p\Delta p + \nabla p\cdot\nabla u + (m-1) p \rho.   
\end{equation}
Let $\eta:\R^d\to \R$ be a function depending on the space variable $x$ only. Multiplying~\eqref{eq:eq_for_pressure} by $\eta$ we obtain
\begin{equation}\label{eq:eq_of_p_eta}
\partial_t(p\eta) - |\nabla p|^2\eta - (m-1)\eta p \Delta p + \eta\nabla p \cdot\nabla u  + (m-1)p\eta\Delta u=0,
\end{equation}
which can be rewritten as

\begin{multline*}
\partial_t(p\eta) = \frac{|\nabla(p\eta)-p\nabla \eta|^2}{\eta} + (m-1)p\Delta (p\eta) - 2(m-1)p\frac{\nabla(\eta p)-p\nabla\eta}{\eta}\cdot \nabla \eta  - (m-1)p^2\Delta \eta\\
- \nabla(\eta p)\cdot\nabla u + p\nabla \eta \cdot\nabla u  -(m-1)p\eta \Delta u.
\end{multline*}
We choose $\eta$ with $|\Delta\eta|\lesssim \eta$ and $|\nabla\eta|\lesssim \eta $  as well as $\eta\geq 1$ (thinking of $\eta$ as being polynomial). Afterwards we set $f=p\eta$ and we use  the boundedness of $p$, of $\nabla u$ (which is a consequence of \Cref{boundsDuDDu}) and of $\Delta u$ (which equals $-\rho$) to obtain 
$$
    \partial_t f \leq |\nabla f|^2+w\cdot \nabla f+Cf+ (m-1) p \Delta f,\quad \partial_t f \geq w\cdot \nabla f+Cf+ (m-1) p \Delta f,
$$
where the vector field $w$ is given by $w=-2mp\nabla\log\eta-\nabla u$ and is thus bounded. The constant $C$ depends on the $L^\infty$ bounds on $\rho$ and $p$ and on the bounds on the ratios $|\Delta\eta|/\eta$ and $|\nabla\eta|/\eta$. The function $f$ is also bounded from below (since it is nonnegative) and from above if we also assume $\eta$ to be bounded (but please observe that the bounds we write will not depend on $\|\eta\|_{L^\infty}$).

A standard parabolic comparison principle together with a simple use of the Gr\"onwall Lemma provides bounds on $f$ in terms of its bounds at time $t=0$. Taking $\eta=(1+|x|)^{\sigma}$ proves the desired lower bound on $\rho$. For the upper bound we first use a truncated version of the same $\eta$, i.e. $\tilde \eta=K\arctan(\frac 1K \eta)$. We obtain $p(t,x)\tilde\eta(x)\leq e^{Ct}\sup ( p_0\tilde\eta ) \leq \sup ( p_0\eta )$ and we then take the limit $K\to\infty$, thus obtaining the upper bound on $\rho$.

The case $d=2$ is slightly different as we have $p=\log\rho$ and we lose boundedness from below of the pressure. We then perform a similar analysis on the equation on $\rho$ instead of $p$, that can be written as
$$
    \partial_t\rho=\Delta\rho-\nabla\rho\cdot\nabla u +\rho^2
$$
and, multiplying times $\eta$, becomes
$$
    \partial_t(\rho\eta)=\Delta(\rho\eta)-2\nabla(\rho\eta)\cdot\frac{\nabla\eta}{\eta}+2\rho\frac{|\nabla\eta|^2}{\eta}-\rho\Delta\eta-\nabla(\rho\eta)\cdot\nabla u +\rho\nabla\eta\cdot\nabla u+\rho^2\eta.
$$
Using the same bounds as before and setting $f=\rho\eta$ we then obtain
$$\partial_t f= \Delta f+\nabla f\cdot w+c(t,x)f,$$
where $w=-2\nabla \log\eta-\nabla u$ and $c$ is a function bounded by a constant only depending on the $L^\infty$ bounds on $\rho$ and on the ratios $|\Delta\eta|/\eta$ and $|\nabla\eta|/\eta$. The conclusion is the same as before taking this time $\eta=(1+|x|)^\beta$ or a truncation of it.
\end{proof}

\begin{lemma}\label{boundsDuDDu}
If $\rho$ is a probability density belonging to $L^\infty(\R^d)$, then for $u=\Gamma*\rho$ we have $\nabla u\in L^\infty(\R^d)$. If moreover $\rho$ satisfies $\rho(x)\leq C(1+|x|)^{-\beta}$ for $\beta>d$, then we have $|\nabla u(x)|\leq C(1+|x|)^{-(d-1)}$.

If $\rho$ is also $C^{0,\alpha}(\R^d)$ for $\alpha\in(0,1]$, then $D^2u $ is bounded and if on each ball $B(x_0,R)$ with $|x_0| = 2R$ the values of $\rho$ are bounded by $C(1+R)^{-\beta}$ and the $C^{0,\alpha}$ seminorm by $C(1+R)^{-(d+\alpha)}$, then on the ball $B(x_0,R/2)$ we have $|D^2u|\leq C(1+|x|)^{-d}$.
\end{lemma}
\begin{proof}
We have
$$
    |\nabla u(x)|\leq C\int_{\R^d}\frac{\rho(x-y)}{|y|^{d-1}}dy,
$$
hence, by separating the integral on $B_r$ and $B_r^c$ we obtain
$$
    |\nabla u(x)|\leq C\|\rho\|_{L^\infty(B(x,r))}\int_{B_r}\frac{1}{|y|^{d-1}}dy+\frac{C}{r^{d-1}}.
$$
In order to prove $\nabla u\in L^\infty$ it is enough to choose $r=1$ and observe that the integral of $|y|^{1-d}$ converges on each bounded set.  Its value on the ball $B_r$ equals $\int_{B_r}\frac{1}{|y|^{d-1}}dy=Cr$ where $C$ is a dimensional constant. If we use $\rho(x)\leq C(1+|x|)^{-\beta}$ we then just need to take $r=|x|/2$, so that $ \|\rho\|_{L^\infty(B(x,r))} \leq C(1+|x|)^{-\beta}$. Hence, we obtain
$$|\nabla u(x)|\leq C(1+|x|)^{-\beta}|x|+C|x|^{1-d}.$$
The first term in the right-hand side is smaller than $O(|x|^{1-d})$, which provides the claim.

We now go on with the second-order estimates. Let us call $u''$ a combination of second derivatives of $u$ which is either $u_{ij}$ for $i\neq j$ or $u_{ii}-u_{jj}$ also for $i\neq j$. We can write

\[
u''(x) 
  = 
C \int_{\R^d} \frac{\rho(x-y)}{| y|^d}c(y)\,\diff y,
\]
where $c$ is a function defined on $\R^d\setminus\{0\}$ which only depends on the direction of $y$ (it is $0$-homogeoneous) and has zero average on each sphere (more precisely, we have $c(y)=\frac{y_iy_j}{|y|^2}$ in the case $u''= u_{ij}$ and $c(y)=\frac{y_i^2-y_j^2}{|y|^2}$ in the case $u''=u_{ii}-u_{jj}$). Because of this zero-average property we can also write
\[
u''(x) 
  = 
C \int_{\R^d} \frac{\rho(x-y)-\rho(x)}{| y|^d}c(y)\,\diff y.
\]
We then fix $r\leq R$ and separate the integral into three parts: $B_r$, $B_R\setminus B_r$, and $B_R^c$. We have
\[
|u''(x)|\leq C[\rho]_{C^{0,\alpha}(B_r)}\int_{B_r}\frac{|y|^\alpha}{|y|^d}dy+ C[\rho]_{L^\infty(B_R)}\int_{B_R\setminus B_r}\frac{1}{|y|^d}dy+\frac{C}{R^d}
\]
For the proof of $D^2u\in L^\infty$ it is enough to choose $r=R=1$ and observe that we have $\int_{B_r}\frac{|y|^\alpha}{|y|^d}dy=Cr^\alpha<+\infty$. This proves that all the mixed derivatives and all the differences of pure second derivatives are bounded. Using the fact that the average of the pure second derivatives is $\frac{\Delta u}{d}=-\frac{\rho}{d}$, we obtain the result for the whole Hessian. 

Using now the more refined assumption on the behavior of $\rho$ on each ball of the form $B(x_0,R)$ with $|x_0| = 2R$ we choose $r=R/2$ and obtain
$$|u''(x_0)|\leq C(1+R)^{-(d+\alpha)}R^\alpha+C(1+R)^{-\beta}+CR^{-d}\leq C(1+R)^{-d},$$
where we used $\int_{B_R\setminus B_r}\frac{1}{|y|^d}dy=C\log(R/r)$.  This proves the desired result (again, using the Laplacian to estimate the pure second derivatives) for $|x_0|\geq 1$. For $|x_0|\leq 1$ it is enough to use the boundedness of $D^2u$ which we just proved.
\end{proof}

Finally, we show that we can provide suitable decays for the gradient and the Hessian of the pressure for $t>0$. We rely on the following result from classical parabolic theory, which we state in a convenient form.

\begin{lemma}\label{lem:regularity_lieberman}
Let $d\ge 1$, $T>0$ and $R>0$.  
Assume that $u\in C^{2,1} \left([0,T]\times B_R \right)$ solves
\begin{equation}\label{eq:mainPDE}
    \partial_t u  =  a(t,x,u)\Delta u  +  b \left(t,x,u,\nabla u \right),
\end{equation}
with smooth nonnegative initial data $u_0$, a diffusion coefficient $a$ bounded above and below by positive constants, and with a nonlinear term $b$ that is at most quadratic in the gradient variable, i.e. $|b(t,x,u,\xi)|\leq C(1+|\xi|^2)$. Also assume Lipschitz dependence of $a$ and $b$ in terms of $u$. Then we have the following facts:
\begin{itemize}
\item if $\nabla u_0$ is bounded then $\nabla u_t$ is also bounded uniformly on $B_{R/2}$ for every $t\in [0,T]$;
\item if $D^2 u_0$ is bounded then $D^2 u_t$ is also bounded uniformly on $B_{R/2}$ for every $t\in [0,T]$.
\end{itemize}
All the above bounds on $B_{R/2}$ (on the 
$L^\infty$ norms of the gradient and of the Hessian) depend on the corresponding bound on $u_0$, on the bounds on $a$ and $b$ (for the last one, also on the Lipschitz bounds of $a$ and $b$ in terms of $x$), on the $L^\infty$ norm of $u$, and on $R$ and $T$. 
These bounds degenerate as $R\to 0$ or $T\to\infty$. 
\end{lemma}

This lemma can be inferred from the theory for parabolic partial differential equations, see for instance~\cite{LSU68, Amann90, Lieberman96}. For the particular case we mention in \Cref{lem:regularity_lieberman} where we use $u \in C^{2,1} ([0,T] \times B_R)$ we can refer to \cite[Theorem 2.1]{Wen_Fan_Asiri_Alzahrani_El-Dessoky_Kuang17} which extends the theory in \cite{LSU68}. 

Instead, we apply this fact to precise estimates on the behavior of $\rho$ and $p$ at infinity. 

\begin{proposition}\label{prop:lipschitz-porous}
Let $\rho:[0,T]\times\mathbb R^{d}\to[0,\infty)$ be a bounded solution of the Keller--Segel system with smooth nonnegative initial density $\rho_0$ such that
\[
     \rho_0(x)
     \sim|x|^{-\beta},
     \quad\text{as }|x|\to\infty,
\]
for an exponent $\beta>d$ but close to $d$. If $d>2$ we denote by $\sigma:=\beta\frac{d-2}{d}$, so that $\sigma$ is larger than $d-2$ but close to it. In this case, setting as usual $p_0=\rho_0^{1-2/d}$, we have $p_0(x)\sim |x|^{-\sigma},$ and we also assume
$$|\nabla p_0(x)|\le C(1+|x|)^{-\sigma-1}, \quad |D^2p_0(x)|\leq C(1+ |x|)^{-\sigma-2}.$$
Then, there exist a constant $C$ such that
\begin{equation}\label{eq:MainGradDecay}
     |\nabla p(t,x)|
      \le 
     C(1+|x|)^{-(d-1)},\quad |D^2p(t,x)|\leq   C(1+|x|)^{\sigma+2-2d}
     \quad
     \mbox{for all }x\in\mathbb R^{d},\,t\in[0,T].
\end{equation}

For $d=2$, we use instead $p_0=\log\rho_0$ and we assume 
$$|\nabla p_0(x)|\le C(1+|x|)^{-1}, \quad |D^2p_0(x)|\leq C(1+ |x|)^{-2}.$$
Then, there exist a constant $C$ such that
\begin{equation}\label{eq:MainGradDecay2}
     |\nabla p(t,x)|
      \le 
     C(1+|x|)^{-1},\quad |D^2p(t,x)|\leq   C(1+|x|)^{-2}
     \quad
     \mbox{for all }x\in\mathbb R^{2},\,t\in[0,T].
\end{equation}
\end{proposition}

\begin{proof}
We start from the case $d>2$. By Proposition~\ref{prop:upper_lower_bound_pressure}, we deduce $\rho(x)\sim |x|^{-\beta}$ and $p(x)\sim |x|^{-\beta\left(1-\frac{2}{d}\right)}$ as $|x|\to \infty$. If we consider $u=\Gamma*\rho$ we deduce 
$$
|Du|(x)\lesssim|x|^{1-d} \quad \text{as $|x|\to \infty$}.
$$
uniformly in $t$.

We now consider the equation satisfied by $p$, which is \eqref{eq:eq_for_pressure}, that we write in the form
$$\partial_t p=a\Delta p+f(p)+|\nabla p|^2-w\cdot \nabla p,$$
where $a=(m-1)p$, $f(p)=p\rho$, and $w=\nabla u$.
We consider this equation in a domain $[0,T]\times B_{0}$, where $B_{0}$ is a ball of radius $R_0$ centered at a point $x_0$ with $|x_0|=2R_0\gg 1 $. We then take a suitable scaling of the form 
$$\tilde p(t,x)=Mp(Kt,Rx)$$
and we look at the equation satisfied by $\tilde p$  on a domain of the form $[0,T/K]\times \tilde B$, where $\tilde B$ is a ball of radius $R_0/R$. The equation satisfied by $\tilde p$ is 
$$\partial_t \tilde p=\frac{K}{R^2}\tilde a\Delta \tilde p+g(\tilde p)+\frac{K}{MR^2}|\nabla \tilde p|^2-\tilde w\cdot \nabla \tilde p.$$
Here $\tilde a$ is proportional to $\tilde p$ and the vector field $\tilde w$ is given by $\tilde w(t,x)=\frac KR\nabla u(Kt,Rx)$.

We call $\bar P_0$ a typical value of $p$ on the ball $B_0$, i.e. $\bar P_0=O(R_0^{-\sigma})$, with $\sigma=\beta(1-2/d)$.  We then choose 
$$M\bar P_0=1, \qquad R=\bar P_0R_0^{d-1},\qquad K\bar P_0=R^2.$$
Note that the second condition implies $R=O(R_0^{-\sigma+d-1}),$ and the exponent $-\sigma+d-1$ is very close to $1$ (but smaller), so that we have $R=o(R_0)$. 

With this choice we obtain at the same time that $\tilde p$ is bounded, and that the coefficient in front of $\Delta\tilde p$ is of order 1 and it is a (uniformly) Lipschitz function of $\tilde p$. Moreover the coefficient in front of $|\nabla \tilde p|^2$ is equal to $1$ and $\tilde w$ satisfies $|\tilde w|\leq C\frac{K}{R}R_0^{-(d-1)}=C$. The function $g$ is given by $g(\tilde p)=MKp^{\frac{2d-2}{d-2}}=M^{-\frac{d}{d-2}}K\tilde p^{\frac{2d-2}{d-2}}$. The Lipshitz constant of $g$ in terms of $\tilde p$ is then of the order of $M^{-\frac{d}{d-2}}K=\bar P_0^{\frac{d}{d-2}}K=R^2\bar P_0^{\frac{2}{d-2}}\leq R_0^2-2\frac{\sigma}{d-2}\leq C$.

Moreover, the ball $\tilde B$ is of radius $R_0/R$, bounded below by a constant and $K=R^2\bar P_0^{-1}$ is also bounded below by a constant, hence the domain where we consider $\tilde p$ satisfies upper bounds on its time length and lower bounds on its spatial size.

We now look at the initial bounds on $\tilde p$ at time $t=0$, and in particular at its gradient. We have $|\nabla\tilde p(0,x)|=MR|\nabla p(0,Rx)|\lesssim MR\bar P_0/R_0\lesssim M\bar P_0=1$. From Lemma~\ref{lem:regularity_lieberman}  we deduce a uniform bound on $|\nabla \tilde p|$.

This bound translates into $|\nabla p(t,x)|\leq \frac{C}{MR}=\frac{C\bar P_0}{R}=O(R_0^{-(d-1)})$ (note that this bound is slightly worse than the one that we had for $t=0$, where we had $|\nabla p|\leq \bar P_0/R_0=R_0^{-\sigma-1}$). 

We now use Lemma \ref{boundsDuDDu} with $\alpha=1$ in order to obtain $|D^2u(x)|\leq C(1+|x|^d)^{-1}$. To do this, we need to prove $|\nabla \rho(t,x)|\leq (1+|x|^{d+1})^{-1}$, but we have
$$|\nabla \rho(t,x)|=Cp(t,x)^{\frac{2}{d-2}}|\nabla p(t,x)|=C\rho(t,x)^{\frac 2d}|\nabla p(t,x)|\leq C(1+|x|^\beta)^{-\frac 2d}(1+|x|^{d-1})^{-1}$$
and the condition $\beta>d$ provides the desired estimate. 

We then look again at the equation satisfied by $\tilde p$. Looking at the vector field $\tilde w$ and recalling 
$\tilde w(t,x)=\frac KR\nabla u(Kt,Rx)$ we have 
$$D\tilde w(t,x)=K D^2u(Kt,Rx)$$
and hence $\Lip(\tilde w)\leq K R_0^{-d}=R^2\bar P_0^{-1}R_0^{-d}=O(R_0^{-2\sigma+2d-2+\sigma-d}).$ The exponent $-2\sigma+2d-2+\sigma-d$ equals $-\sigma+d-2<0$ and hence we have $\Lip(\tilde w)\leq C$.

We can thus apply the second part of Lemma~\ref{lem:regularity_lieberman}  in order to obtain uniform bounds on $D^2\tilde p$. This requires, of course, to control $D^2\tilde p$ at time $t=0$. We then use $|D^2\tilde p(0,x)|=MR^2|D^2 p(0,Rx)|\leq MR^2\bar P_0/R_0^2\leq M\bar P_0=1$.

From Lemma~\ref{lem:regularity_lieberman} we deduce a uniform bound on $|D^2 \tilde p|$ which, again, translates into the estimate $|D^2 p(t,x)|\leq \frac{C}{MR^2}=\frac{C\bar P_0}{R^2}=O(R_0^{-\sigma+2\sigma-2(d-1)})$. The exponent is now $\sigma-2(d-1)$, which, using $\sigma>d-2$, gives an exponent slightly worse than $-d$, but it can taken as close as we want to it. 

We now consider the case $d=2$. In this case the equation satisfied by $p=\log \rho$ is 
$$
    \partial_t p =\Delta p+\rho+|\nabla p|^2-\nabla u\cdot \nabla p.
$$
We perform a similar scaling as before, but we choose now $M=1$, $R=R_0$, and $K=R^2$.
We then obtain
$$\partial_t\tilde p =\Delta\tilde p+g(\tilde\rho)+|\nabla\tilde p|^2-\tilde w\cdot \nabla \tilde p, $$
where $\tilde w=\frac KR \nabla u(Kt,Rx)=R\nabla u(Kt,Rx)$.
With our choices for the scaling parameters we have $|\tilde w|\leq C$. Moreover, the oscillation of $p$ on the ball $B_0=B(x_0,R_0)$ is bounded by a universal constant (depending on $\beta$). The same bounds holds of course for $\tilde p$.
The function $g(\tilde p)=K\exp(\tilde p)$ is bounded by $K(1+R_0^{\beta})^{-1}=O(R_0^{2-\beta})\leq C$ and the Lipschitz constant of $g$ in terms of $\tilde p$ is also proportional to $K\max \exp(\tilde p)\leq C$.

Again we want to apply the results of Lemma~\ref{lem:regularity_lieberman} and obtain a uniform Lipschitz bound on $\tilde p$. To do this we only need to consider the initial Lipschitz bound, which is given by $\frac KR\,\Lip(p(t=0), B_0)=R_0\, \Lip(p(t=0), B_0)\leq C$. The function $\tilde p$ is then uniformly Lipschitz in space, which means that $\nabla p(t,x)$ for $x\in B(x_0,R_0/2)$ is bounded by $CR_0^{-1}$.

From this, using $\rho=\exp{(p)}$, we deduce that on the same ball $B(x_0,R_0)$ we have $\nabla \rho=\rho\nabla p\leq R_0^{-1-\beta}$ and we can apply Lemma~\ref{boundsDuDDu}  to obtain $D^2u(t,x)\leq C(1+|x|^2)^{-1}$. We then obtain a uniform Lipschitz bound in space for $\tilde w$. Together with the initial bound $D^2p(t=0)\leq R_0^{-2}$ we conclude a uniform bound for $D^2\tilde p(t,x)$ and hence $D^2 p(t,x)\leq C(1+|x|)^{-2}$.
\qedhere 
\end{proof}

\appendix

\section*{Acknowledgements}

This work was supported by the European Union via the ERC AdG 101054420 EYAWKAJKOS project. AFJ was also partially funded via the ERC AdG 883363 Nonlocal-CPD (Nonlocal PDEs for Complex Particle Dynamics: Phase Transitions, Patterns and Synchronization). The authors are thankful to Pierre-Damien Thizy (Lyon 1) for useful discussions about the Lane--Emden equation.

\section*{AI Statement}

This paper has been fully written by humans. However, the proof of Proposition 5.10 (which is the most concrete improvement w.r.t. a previously posted version of the same paper) has been obtained with the help of LLM tools, which suggested, in interaction with some of the authors, the use of the function $z$, which turned out to be the missing point in the argument. 

\FloatBarrier

{
\addcontentsline{toc}{section}{References}
	\small
	\bibliographystyle{abbrv}
	\bibliography{references}	
}

\end{document}